\documentclass[11pt]{article}
\usepackage{amsfonts}
\usepackage{amsmath}
\usepackage{amssymb}
\usepackage{latexsym}

\usepackage{diagrams}

\usepackage[all]{xy}

\newarrow{Embed} >---> 
\newarrow{Onto} ----{>>}

\newtheorem{thm}{Theorem}[section]
\newtheorem{cor}[thm]{Corollary}
\newtheorem{lemma}[thm]{Lemma}
\newtheorem{prop}[thm]{Proposition}

\newtheorem{exa}[thm]{Example}

\newtheorem{rem}[thm]{Remark}

\newtheorem{claim}[thm]{Claim}

\def\qed{{\hspace{2mm}{\hfill \small $\Box$}}}
 \def\pf{{\noindent{\bf Proof.\hspace{2mm}}}}

\def\C{{\mathbb C}}

\def\H{{\mathbb H}}
\def\Z{{\mathbb Z}}

\begin{document}

\title{On character varieties, sets of discrete characters, and non-zero degree maps}

\author{Michel Boileau\footnote{{ Partially supported
by the Laboratoire de Math\'ematiques \'Emile Picard, UMR CNRS 5580}}   
\hspace{.2mm} and Steve Boyer\footnote{{ Partially supported by NSERC and FQRNT}}}

\date{} 

\maketitle

\section {Introduction}

\subsection{General introduction} \label{gen}

Character variety methods have proven an essential tool for the investigation of problems in low-dimensional topology and have been instrumental in the resolution of many well-known problems. In this paper we use them to study homomorphisms between the fundamental groups of $3$-manifolds, in particular those induced by non-zero degree maps. We assume throughout that our manifolds are compact, connected, orientable, and $3$-dimensional. A {\it knot manifold} is a compact, connected, irreducible, orientable $3$-manifold whose boundary is an incompressible torus. We shall restrict our attention, for the most part, to {\it small} knot manifolds, that is, those which contain no closed essential surfaces. This is a simplifying hypothesis and though many of the results discussed in the paper extend to the general setting, we will not discuss them. A small knot manifold is atoroidal and Haken, hence it is either hyperbolic or admits a Seifert fibred structure with base orbifold of the form $D^2(p, q)$ for some integers $p, q \geq 2$. 

We call a homomorphism $\varphi: \Gamma_1 \to \Gamma_2$ between two groups a {\it virtual epimorphism} if its image is of finite index in $\Gamma_2$. For instance, a non-zero degree map between manifolds induces a virtual epimorphism on the level of fundamental groups. The first part of the paper investigates virtual epimorphisms between the fundamental groups of small knot manifolds. We will see that the existence of such homomorphisms places constraints on the algebraic decomposition of a knot manifold's $PSL_2(\mathbb C)$-character variety and, as a consequence, we will determine a priori bounds on the number of virtual epimorphisms between the fundamental groups of small knot manifolds with a fixed domain. 
This work yields minimality results which will be applied to illustrate the results of the second part of the paper. There we fix a small knot manifold $M$ and investigate various sets of characters of representations $\rho: \pi_1(M) \to PSL_2(\mathbb C)$ whose images are discrete. It turns out that the topology of these sets is intimately related to the algebraic structure of the $PSL_2(\mathbb C)$-character variety of $M$ as well as dominations of manifolds by $M$ and its Dehn fillings. In particular, we apply our results to study families of non-zero degree maps $f_n: M(\alpha_n) \to V_n$ where $M(\alpha_n)$ is the $\alpha_n$-Dehn filling $M$ and $V_n$ is either a hyperbolic manifold or $\widetilde{SL_2}$ manifold. Using this, the existence of infinite families of small, closed, connected, orientable manifolds which do not admit non-zero degree maps, other than homeomorphisms, to any hyperbolic manifold, or even manifolds which are either reducible, Haken, or admit a geometric structure, is determined.

In the remainder of the introduction we give a more detailed description of our results and the organization of the paper. Here is some notation and terminology we shall use.

Throughout, $\Gamma$ will denote a finitely generated group. We call a homomorphism $\rho: \Gamma \to PSL_2(\mathbb C)$ discrete, non-elementary, torsion free, abelian, etc. if its image has this property. If $\chi_\rho \in X_{PSL_2}(M)$ is the character of $\rho$ we will call it discrete, non-elementary, torsion free, abelian, etc. if each representation $\rho': \pi_1(M) \to PSL_2(\mathbb C)$ with $\chi_{\rho'} = \chi_\rho$ has this property. For instance we can unambiguously refer to a character as being either irreducible, non-elementary, or torsion free. 

A {\it slope} on the boundary of a knot manifold $M$ is a $\partial M$-isotopy class of essential simple closed curves. Slopes correspond bijectively with $\pm$ pairs of primitive elements of $H_1(\partial M)$ in the obvious way. The {\it longitudinal slope} on $\partial M$ is the unique slope $\lambda_M$ having the property that it represents a torsion element of $H_1(M)$. When $M$ is the exterior of a knot $K$ in a closed $3$-manifold $W$, there is a unique slope $\mu_K$ on $\partial M$, called the {\it meridinal slope}, which is homologically trivial in a tubular neighbouhood of the knot. If $W$ is a $\mathbb Z$-homology $3$-sphere, then $\mu_K$ and $\lambda_M$ are dual in the sense that the homology classes they carry form a basis for $H_1(\partial M)$. 

Each slope $\alpha$ on $\partial M$ determines an element of $\pi_1(M)$ well-defined up to conjugation and taking inverse. We will sometimes use this connection to evaluate a representation on a slope, but only in a context where the statement being made is independent of the choice of element of $\pi_1(M)$. For instance we may say that $\rho(\alpha) \in PSL_2(\mathbb C)$ is parabolic, or loxodromic, or trivial. 

A representation $\rho: \pi_1(M) \to PSL_2(\mathbb C)$ is {\it peripherally nontrivial} if $\rho(\pi_1(\partial M))$ does not equal $\{\pm I\}$. When $M$ is small, there are only finitely many characters of representations which are not peripherally non-trivial. Indeed, there are only finitely many characters $\chi_\rho$ for which $\rho(\pi_1(\partial M))$ is trivial or a parabolic subgroup of $PSL_2(\mathbb C)$ (cf. Corollary \ref{smallcharactervariety}). Thus, apart from finitely many exceptions, a discrete, torsion-free character is the character of a representation $\rho$ for which there is a unique slope $\alpha$ on $\partial M$ such that $\rho(\alpha) = \pm I$. In this case we call $\alpha$ the \emph{slope of $\rho$}.

A \emph{hyperbolic manifold} is one whose interior admits a complete, finite volume, hyperbolic structure. A closed manifold which admits an $\widetilde{SL_2}$ structure is called an $\widetilde{SL_2}$ manifold. Similarly we will refer to closed $Nil$ manifolds, $\mathbb E^3$ manifolds, etc. Two families of manifolds we will focus on are the family ${\cal H}$ of hyperbolic $3$-manifolds and the family ${\cal M}$ of $3$-manifolds which are either reducible, Haken, or admit a geometric structure. According to Thurston's Geometrization Conjecture, which has been claimed by Perelman, ${\cal M}$ is the set of all compact, connected, orientable $3$-manifolds. We shall assume this holds below. 

We say that $M$ \emph{dominates} $N$, written $M \geq N$, if there is a continuous, proper map from $M$ to $N$ of non-zero degree. Moreover, if $N$ is not homeomorphic to $M$ we say that $M$ \emph{strictly dominates} $N$. The relation $\geq$ is a partial order when restricted to manifolds in ${\cal M}$ which are aspherical but are neither torus (semi) bundles or Seifert manifolds with zero Euler number \cite{Wan1}, \cite{Wan2}. This partial order is far from well-understood, even when restricted to hyperbolic 3-manifolds. 

A knot manifold is {\it minimal} if the only knot manifold it dominates is itself. (Note that each knot manifold dominates $S^1 \times D^2$, but also that the latter is not a knot manifold.) For example, using elementary $3$-manifold topology we can see that the total space of a punctured torus bundle is minimal if and only if its monodromy is not a proper power (see  [BWa1, Prop. 2.6]).  A closed, connected, orientable $3$-manifold is {\it minimal} if the only manifold it dominates is one with finite fundamental group. (It is easy to see that every closed, connected, orientable $3$-manifold dominates each such manifold with finite fundamental group.) A manifold $V$ is ${\cal H}$-{\it minimal} if the only manifold in ${\cal H}$ it dominates is itself. Note that we do not require that $V \in {\cal H}$. Similarly we can define the notion of an ${\cal M}$-{\it minimal} manifold. If the reader prefers not to assume the geometrization theorem of Perelman, then in what follows, a closed, minimal $3$-manifold should be taken to mean a closed, ${\cal M}$-minimal manifold. 

\subsection{Virtual epimorphisms of knot manifold groups and domination} 

Let $X_0$ be an algebraic component of $X_{PSL_2}(\Gamma)$, the $PSL_2(\mathbb C)$-character variety of $\Gamma$ (see \S \ref{varieties}). In \S \ref{subvarieties} we note that $X_0$ determines a normal subgroup $\hbox{Ker}(X_0) \subset \Gamma$ with the property that for the generic character $\chi_\rho \in X_0$, $\hbox{Ker}(X_0) = \hbox{kernel}(\rho)$. For a small knot manifold $M$ there is a closely related normal subgroup $K_M(X_0)$ (\S \ref{seminorm}). Define ${\cal I}_M$ to be the set of isomorphism classes of groups $\pi_1(M) / K_M(X_0)$ where $X_0$ ranges over the set of algebraic components which contain an irreducible character. In Theorem \ref{epibound} we show that the cardinal $|{\cal I}_M|$ gives an upper bound for number of isomorphism classes of 
groups $\pi_1(N)$ where $N$ is a small knot manifold for which there is an epimorphism $\varphi: \pi_1(M) \to  \pi_1(N)$. For epimorphisms induced by non-zero degree maps we obtain an interesting 
refinement. 

The set of algebraic components of $X_{PSL_2}(M)$, the $PSL_2(\mathbb C)$-character variety of the fundamental group of a small knot manifold $M$, are partitioned into two types - those whose Culler-Shalen seminorms are norms and those which are not (see \S \ref{culler-shalen}). Let ${\cal N}_M$ be the set of isomorphism classes of groups $\pi_1(M) / K_M(X_0)$ where $X_0$ ranges over the set of algebraic components whose associated Culler-Shalen seminorm is a norm.  

We will consider two such maps $f_j: M \to N_j$ ($j = 1,2$) between knot manifolds to be {\it equivalent} if there is a homeomorphism $h: N_1 \to N_2$ such that the following diagram commutes up to homotopy: 
$$\xymatrix@R=20pt@C=70pt
{& (N_1, \partial N_1)\ar[dd]^{h}  \\ 
(M, \partial M)\ar[ur]^{f_1}\ar[dr]_{f_2}   & \\ 
& (N_2, \partial N_2)}  $$
The following theorem is proven in \S \ref{seminorm} (see Theorem \ref{dombound}).  
\begin{thm}   \label{domboundintro} 
Let $M$ be a small knot manifold. \\ 
$(1)$ The number of equivalence classes of $\pi_1$-surjective non-zero degree maps $M \to N$ is bounded above by $|{\cal I}_M|$. More precisely, \\ 
\indent $(a)$  The number of equivalence classes of $\pi_1$-surjective non-zero degree maps $M \to N$ \\ \indent \hspace{.55cm} where $N$ is hyperbolic is bounded above by $|{\cal N}_M|$. \\ 
\indent $(b)$ The number of equivalence classes of $\pi_1$-surjective non-zero degree maps $M \to N$  \\ \indent \hspace{.55cm} where $N$ is Seifert is bounded above by $|{\cal I}_M| - |{\cal N}_M|$. \\ 
$(2)(a)$ The number of equivalence classes of non-zero degree maps from $M$ to a hyperbolic  \\ \indent \hspace{.45cm} manifold is bounded above by a constant depending only on $X_{PSL_2}(M)$. \\
\indent $(b)$ The number of homeomorphism classes of Seifert fibred manifolds dominated by $M$  \\ \indent \hspace{.45cm} is bounded above by a constant depending only on $X_{PSL_2}(M)$.
\end{thm}
This result can be used to give many examples of minimal and ${\cal H}$-minimal knot exteriors. For instance it is shown in Example \ref{twistpretzel} that the exterior of the $(-2,3,n)$ pretzel knot is ${\cal H}$-minimal while those of twist knots or $(-2,3,n)$ pretzel knots ($n \not \equiv 0$ (mod $3$)) are minimal.  On the other hand, the theorem's usefulness is tempered by the difficulty, in general, of determining the algebraic decomposition of $X_{PSL_2}(M)$. Indeed, there only handful families of knot manifolds whose character varieties have been explicitly determined. For instance, little definite information is known on the number of algebraic components of the character varieties of two-bridge knot exteriors. We can, nevertheless, use character variety methods to study dominations by such manifolds. In \S \ref{rigidity} we associate to each small knot exterior $M$ a function $d_M: \pi_1(M) \to \mathbb Z$ with the property that if $\varphi: \pi_1(M) \to \pi_1(N)$ is a virtual epimorphism, then $d_N(\varphi(\gamma)) \leq d_M(\gamma)$ for all $\gamma \in \pi_1(M)$. Further, if $\gamma$ is not rigid (\S \ref{rigidity}), then $d_N(\varphi(\gamma)) = d_M(\gamma)$ implies that $\varphi$ is 1-1. This leads to our next result (see Theorem \ref{homomorphismsequence}).  

\begin{thm}  \label{intro:domseqbound} 
Let $M$ be a small knot manifold and consider a sequence of virtual epimorphisms    
$$\pi_1(M) \stackrel{\varphi_1}{\longrightarrow} \pi_1(N_1) \stackrel{\varphi_2}{\longrightarrow}  \cdots \stackrel{\varphi_n}{\longrightarrow}  \pi_1(N_n)$$ 
none of which is injective. If $N_i$ is small for each $i$, then $n \leq d_M(\gamma)$ for each totally non-rigid element $\gamma \in \pi_1(M)$. Moreover, if $n = d_M(\gamma)$, then $N_n$ is a twisted $I$-bundle over the Klein bottle. 
\end{thm}
The non-injectivity of $\varphi_i$ is a necessary condition since small knot manifolds which are Seifert fibred admit self-coverings of arbitrarily large degree. 

Precise calculations of $d_M(\gamma)$ can be made for various families of knot manifolds. In the case where $M$ is the exterior of a two-bridge knot we obtain the following explicit bounds. (See Theorems \ref{2bridgebound} and \ref{2bridgedegree1bound}.)

\begin{thm} \label{twobridgedomboundintro}
Let $M_{p/q}$ denote the exterior of the $(\frac{p}{q})$ two-bridge knot and consider a sequence of virtual epimorphisms none of which is injective: 
$$\pi_1(M_{p/q}) \stackrel{\varphi_1}{\longrightarrow} \pi_1(N_1) \stackrel{\varphi_2}{\longrightarrow}  \cdots \stackrel{\varphi_n}{\longrightarrow}  \pi_1(N_n)$$ 
$(1)$ If $N_i$ is small for each $i$, then $n < \frac{p-1}{2}$. \\ 
$(2)$ If each $\varphi_i$ is induced by a non-zero degree map, then $n+1$ is bounded above by the number of distinct divisors of $p$. \\
$(3)$ If each $\varphi_i$ is induced by a degree one map, then $n+1$ is bounded above by the number of distinct prime divisors of $p$.
\end{thm}
Theorem \ref{twobridgedomboundintro} immediately yields an infinite family of minimal two-bridge knot exteriors. 

\begin{cor} $\;$ \\ 
$(1)$ If $p$ is prime, then $M_{p/q}$ is minimal if and only if it is hyperbolic  $($i.e. $q \not \equiv \pm 1 \hbox{ $($mod $p$ $)$}$.$)$ \\
$(2)$ If $p$ is a prime power, any degree one map $M_{p/q} \to N$, $N$ a knot manifold, is homotopic to a homeomorphism.  
\end{cor}
Two-bridge knot exteriors are not minimal in general. For instance, T. Ohtsuki, R. Riley, and M. Sakuma have given a systematic construction of degree one maps between two-bridge knot exteriors.

\subsection{Families of discrete characters and domination}

The fundamental group of a small knot manifold $M$ admits many discrete, non-elementary representations with values in $PSL_2(\mathbb C)$. For instance when $M$ is hyperbolic, its holonomy representation is discrete and non-elementary, as are the holonomy representations of the hyperbolic Dehn fillings of $M$. Similarly, when $M$ is Seifert fibred but not a twisted $I$-bundle over the Klein bottle, a holonomy representation of its base orbifold is discrete and non-elementary, as are those of the  base orbifolds of the generic Dehn filling of $M$. One of the problems we investigate in the second part of the paper is to determine to what extent these are the only discrete non-elementary representations of $\pi_1(M)$.  

Set 
$$D(M) = \{ \chi_\rho \in X_{PSL_2}(M): \rho \hbox{ is discrete and non-elementary}\}.$$ 
Classic work of J\o rgensen and Marden shows that $D(M)$ is closed in $X_{PSL_2}(M)$ (see \S \ref{convkleinhyp}). Their results combine with the work of Culler and Shalen on ideal points of curves of $PSL_2(\mathbb C)$-characters to show that if $D(X_0) = D(M) \cap X_0$ is not compact for some component $X_0$ of $X_{PSL_2}(M)$, there is a connected essential surface $S$ in $M$ such that the restriction of each character in $X_0$ to $\pi_1(S)$ is elementary. It follows that $S$ is an annulus if $X_0$ contains the character of a faithful representation. In particular, $M$ is not hyperbolic. (Morgan and Shalen used this approach to prove that if the set of discrete faithful characters of the fundamental group of a compact $3$-manifold is not compact, then the group splits non-trivially along a virtually abelian subgroup \cite{MS3}.) We use these ideas to construct various infinite families of small hyperbolic knot exteriors $M$ for which $D(M)$ is compact (see \S \ref{unbounded}). 

To each representation $\rho: \pi_1(M) \to  PSL_2(\C)$ is associated a volume \emph{vol($\rho$)} $\in \mathbb R$ defined by taking any pseudo-developing map from the universal cover $\tilde M$ into $\H^3$ and integrating the pull-back of the hyperbolic volume form on a fundamental domain of $M$ (see \cite{Dun}, \cite{Fra} for more details). This value depends only on the character of $\rho$ so it makes sense to talk of the volume of a character. Moreover, the associated volume function $vol: X_{PSL_2}(M) \to \mathbb R$ is continuous (indeed analytic). Here is a natural way to construct representations with non-zero volume. 

Let $V$ be a compact, connected, orientable, hyperbolic manifold with holonomy representation $\rho_V: \pi_1(V) \to PSL_2(\mathbb C)$. If $V$ is a knot manifold, $f: M \to V$ a non-zero degree map, and  $\rho = \rho_V \circ f_\#$, then $vol(\rho) = \hbox{degree}(f) vol(V) \ne 0$. Similarly, if $V$ is closed, $M(\alpha)$ is a Dehn filling of $M$, $f: M(\alpha) \to V$ a non-zero degree map, and $\rho$ the composition $\pi_1(M) \to \pi_1(M(\alpha))\stackrel{f_\#}{\longrightarrow} \pi_1(V) \stackrel{\rho_V}{\longrightarrow} PSL_2(\mathbb C)$, then $vol(\rho) = vol(\rho_V \circ f_\#) = \hbox{degree}(f) vol(V) \ne 0$. Note that each of these non-zero volume representations  is discrete and torsion free. 

There is a converse to this construction. The image of a discrete, torsion free, non-zero volume representation $\rho$ is the fundamental group of an element $V$ of ${\cal H}$ which is either a knot manifold or closed depending on whether or not $\rho|\pi_1(\partial M)$ is injective (Lemma \ref{standardimage}) In the former case it is easy to see that $\rho$ is induced, as above, by a map $M \to V$. Since $0 \ne vol(\rho) = \hbox{degree}(f) vol(V)$ so $\hbox{degree}(f) \ne 0$. In the latter case, we use the fact that the abelian subgroups of $\pi_1(V)$ are cyclic to see that there is at least one slope $\alpha$ on $\partial M$ such that $\rho(\alpha) = \pm I$. It follows that $\rho$ is obtained, as above, from a map $f: M(\alpha) \to V$ with $\hbox{degree}(f) \ne 0$. Thus discrete, non-elementary, torsion free, non-zero volume representations of the fundamental group of a knot exterior $M$ correspond to non-zero degree maps of $M$ or its Dehn fillings to a hyperbolic manifold. 

A {\it principal component} $X_0$ of the $PSL_2(\mathbb C)$-character variety of a finitely generated group $\Gamma$ is a component which contains the character of a discrete, faithful, irreducible representation of $\Gamma / Z(\Gamma)$, where $Z(\Gamma)$ denotes the centre of $\Gamma$. Our next result is a combination of Theorem \ref{convord} and Lemma \ref{projectiveconvergence}. 

\begin{thm}\label{thm:characters} 
Let $M$ be a small hyperbolic knot manifold, $X_0$ a component of $X_{PSL_2}(M)$, and suppose that $\{\chi_n\} \subset X_0$ is a sequence of distinct characters of non-zero volume representations $\rho_n: \pi_1(M) \to PSL_2(\mathbb C)$ with image  a torsion-free, cocompact, discrete group $\Gamma_n$. For $n \gg 0$, let $\alpha_n$ be the slope of $\rho_n$. Then up to taking a subsequence, one of the following two possibilities arises: \\
$(a)$  the slopes $\alpha_n$ converge projectively to the class of a boundary slope of $M$; or \\
$(b)$ $\lim \chi_n$ exists and is the character of  a discrete, non-elementary, torsion free, non-zero volume representation $\rho_0$ such that: \\ 
\indent {$(i)$} $\rho_0|\pi_1(\partial M)$ is 1-1 and $\rho_0(\pi_1(M))$ is a finite index subgroup of the fundamental group \\\indent \hspace{4mm} of a $1$-cusped hyperbolic manifold $V$. \\
\indent {$(ii)$} there are slopes $\beta_n$ on $\partial V$ such that for each $n$ the fundamental group of the Dehn  \\ \indent \hspace{4mm} filled manifold $V(\beta_n)$ is isomorphic to $\Gamma_n$ and the character $\chi_n$ is induced by the  \\\indent \hspace{4mm} composition $\pi_1(M) \to \rho_0(\pi_1(M)) \hookrightarrow \pi_1(V) \to \pi_1(V(\beta_n)) \cong \Gamma_n$. \\
\indent {$(iii)$} $X_0 = \rho_0^*(Y_0)$ for a principal component $Y_0$ of $X_{PSL_2}(V)$
\end{thm} 
In certain circumstances we can guarantee that conclusion (b) of the theorem holds. For instance, this is the case when $M$ is hyperbolic and the characters $\chi_n$ lie on a principal component $X_0$ of the $PSL_2(\mathbb C)$-character variety of $\pi_1(M)$ (Corollary \ref{principalcompact}). More generally, it is ruled out if we suppose that the characters $\chi_n$ lie on a curve component $X_0$ of $X_{PSL_2}(M)$ such that one of the following two conditions holds: \\
\indent (a) for each ideal point $x_0$ of $X_0$ there are a component $S_0$ of an essential surface  assoc- \\ \indent \hspace{5.5mm} iated to $x_0$ and a character $\chi \in X_0$ such that $\chi|\pi_1(S_0)$ is non-elementary; or \\
\indent (b) the Culler-Shalen seminorm of $X_0$ is a norm and each ideal point of $X_0$ has an assoc- \\ \indent \hspace{5.5mm} iated essential surface $S_0$ with $|\partial S_0| \leq 2$. \\ 
See Corollary \ref{nonelementarycompact} for the justification of case (a) and Corollary \ref{lessthan2} for that of case (b). The following is a consequence of Theorem \ref{thm:characters} (see Corollary \ref{discetecharsprincipal}.)

\begin{cor}\label{cor:finite volume}
Let $M$ be a small hyperbolic knot manifold. Then all but finitely many of the discrete, non-zero volume characters on a principal curve $X_0$ of the $PSL_2(\mathbb C)$-character variety of $\pi_1(M)$ are induced by the complete hyperbolic structure on the interior of $M$ or by Dehn fillings of manifolds finitely covered by $M$.
\end{cor} 

One of our principal motives for investigating families of discrete characters was to address the following question posed by Shicheng Wang: If there are non-zero degree maps between infinitely many distinct Dehn fillings of two knot manifolds $M$ and $N$, are they induced by a non-zero degree map $M \to N$? Here is a version of Theorem \ref{thm:characters} for non-zero degree maps which provides a partial answer (see \S \ref{dominationhyperbolic}).

\begin{thm}\label{thm:domination} 
Let $M$ be a small hyperbolic knot manifold and suppose that there is a slope $\alpha_0$ on $\partial M$ such that the Dehn filled manifold $M(\alpha_0)$ does not  dominate a hyperbolic manifold. Let $\{\alpha_n\}_{n \geq 1}$ be a sequence of distinct slopes on $\partial M$ which do not subconverge projectively to a boundary slope. If there are dominations $f_n: M(\alpha_n) \geq  V_n$ where $V_n$ is a hyperbolic manifold, then there exist a compact hyperbolic manifold $V_0$ with a domination $f: M \geq V_0$, a subsequence $\{j\}$ of $\{n\}$, and slopes $\beta_j$ on $\partial V_0$ such that: \\
\indent {$(i)$} For each $j$, $V_0(\beta_j) \cong V_j$. \\
\indent {$(ii)$} The following diagrams are commutative up to homotopy : 
$$\begin{array}{ccc} 
M & \stackrel{f}{\longrightarrow} & V_0 \\
\downarrow & & \downarrow \\
M(\alpha_j) & \stackrel{f_j}{\longrightarrow} & V_j \cong V_0(\beta_j) 
\end{array}$$ 
Thus infinitely many of the dominations $f_n: M(\alpha_n) \geq  V_n$ are induced by $f$. If we assume further that the dominations $f_n: M(\alpha_n) \geq V_n$ are strict, then $f: M \geq V_0$ is strict as well.
\end{thm} 

The only known example of closed hyperbolic ${\cal H}$-minimal $3$-manifold is $\frac12$ surgery on the figure eight knot \cite{RWZ}. The following consequences of Theorem \ref{thm:domination} show that closed ${\cal H}$-minimal manifolds are actually quite plentiful. 

We will denote the projective space of $H_1(M; \mathbb R)$ by $\mathbb P(H_1(\partial M; \mathbb R))$ and the class of non-zero element $\beta \in H_1(\partial M; \mathbb R)$ by $[\beta]$. The following is part (1) of Theorem \ref{minimal hyperbolic}. 
 
\begin{thm}\label{thm:h-minimal}
Let $M$ be a small, hyperbolic ${\cal H}$-minimal knot manifold and suppose that there is a slope $\alpha_0$ on $\partial M$ such that the Dehn filled manifold $M(\alpha_0)$ does not dominate any closed hyperbolic manifold. If $U \subset \mathbb P(H_1(\partial M; \mathbb R))$ is the union of disjoint closed arc neighbourhoods of the finite set of boundary slopes of $M$, then $\mathbb P(H_1(\partial M; \mathbb R)) \setminus U$ contains only finitely many projective classes of slopes $\alpha$ such that $M(\alpha)$ is not ${\cal H}$-minimal. In particular, $M$ admits infinitely many ${\cal H}$-minimal Dehn fillings. 
\end{thm}
This theorem applies to many hyperbolic knot manifolds. For instance, it applies to punctured torus bundles whose monodromies are pseudo-Anosov and not proper powers, or the exterior of the $(-2,3,n)$ pretzel ($n \ne 1,3,5$). 

The meridinal slope of a knot in the $3$-sphere whose exterior is small is never a boundary slope (see Theorem 2.0.3 of \cite{CGLS}). Thus Theorem \ref{thm:h-minimal} implies:

\begin{cor}\label{cor:integer surgery}
Let $M$ be the exterior of a small hyperbolic knot in $S^3$ and let $\mu,\lambda \in H_1(\partial M)$ represent the meridinal and longuitudinal slope respectively. If $M$ is ${\cal H}$-minimal, then for all but finitely many $n \in \Z$, the Dehn filled manifold $M(n\mu + \lambda)$ is ${\cal H}$-minimal. 
\end{cor}
For certain two-bridge knot exteriors we can say more (see \S \ref{sec:h-minimal}). 
 
\begin{cor}\label{cor:2-bridge}
Let $M$ be  the exterior of a hyperbolic $\frac{p}{q}$ two-bridge knot with $p$ prime. Then all but finitely many Dehn fillings of $M$ yield ${\cal H}$-minimal manifolds. 
\end{cor} 

In order to construct families of closed minimal manifolds it is necessary to prove a version of Theorem \ref{thm:characters} for discrete representations to $PSL_2(\mathbb R)$. Set 
$$D(M; \mathbb R) = \{ \chi_\rho \in D(M) : \rho \hbox{ has image in } PSL_2(\mathbb R)\}$$ 
and $D(X_0; \mathbb R) = D(M; \mathbb R) \cap X_0$ where $X_0$ is a subvariety of $X_{PSL_2}(M)$. In section \S \ref{realconvergent} we prove a result on 
convergent sequences of characters in $D(M; \mathbb R)$ whose topological 
interpretation is investigated in \S \ref{psl2domination}. An example of the sort of result we obtain is the next theorem (see Corollary \ref{cor:s-domination}).  

\begin{thm} \label{cor:s-dominationintro} 
Let  $M$ will be a small hyperbolic knot manifold with $H_1(M) \cong \mathbb Z$, $\{\alpha_n\}$ a sequence of distinct slopes on $\partial M$, and $\{\chi_n\} \subset D(M; \mathbb R)$ a sequence of characters of representations $\rho_n$ such that $\rho_n(\alpha_n) = \pm I$ for all $n$. If there are infinitely many distinct characters $\chi_n$ and the sequence $\{\chi_n\}$ subconverges to a character $\chi_{\rho_0}$ such that $\rho_0(\lambda_M) \ne \pm I$, then $M$ strictly dominates a Seifert manifold with incompressible boundary. 
\end{thm} 
If we remove the condition that $H_1(M) \cong \mathbb Z$ from the hypotheses of the theorem, we can still construct a non-zero degree map from $M$ to a Seifert orbifold with incompressible boundary. We cannot rule out, though, the possibility that the manifold underlying the orbifold is not a solid torus.  

This last theorem can be used to construct infinite families of closed minimal manifolds. For instance, we have the following theorem (see Theorem \ref{finitesl2domination}). 

\begin{thm} \label{finitesl2dominationintro}
Suppose that $M$ is a small ${\cal H}$-minimal hyperbolic knot manifold which has the following properties:\begin{description} 
\vspace{-.45cm} \item[\hspace{1.5mm} {\rm(a)}]  There is a slope $\alpha_0$ on $\partial M$ such that $M(\alpha_0)$ is ${\cal H}$-minimal.
\vspace{-.35cm}  \item[\hspace{1.5mm} {\rm(b)}]  For each norm curve  $X_0 \subset X_{PSL_2}(M)$ and for each essential surface $S$ associated to an ideal point of $X_0$ there is a character $\chi_\rho \in X_0$ which restricts to a strictly irreducible character on $\pi_1(S)$. 
\vspace{-.35cm} \item[\hspace{1.5mm} {\rm(c)}]  There is no surjective homomorphism from $\pi_1(M)$ onto a Euclidean triangle group.
\vspace{-.35cm} \item[\hspace{1.5mm} {\rm(d)}]  There is no epimorphism $\rho: \pi_1(M) \to \Delta(p,q,r) \subset PSL_2(\mathbb R)$ such that $\rho(\pi_{1} (\partial M))$ is elliptic or trivial.
\end{description}
\vspace{-.45cm}  Then all but finitely many Dehn fillings $M(\alpha)$ yield a minimal manifold.
\end{thm}
As a consequence, we will prove our next result in Corollary \ref{cor:s-minimal2}.

\begin{cor}\label{prop:g-minimalintro}
If $M$ is the exterior of a hyperbolic twist knot, then all but finitely many Dehn filling $M(\alpha)$ yield a minimal manifold. 
\end{cor}

Our final results show that quite general hypotheses on a minimal knot exterior imply that it admits infinitely many minimal Dehn fillings. (See Theorems \ref{uinfinite} and \ref{minginfinite}.) 

\begin{thm} \label{uinfiniteintro} 
Let $M$ be an ${\cal H}$-minimal, small, hyperbolic knot manifold and suppose 
that  $H_1(M) \cong \mathbb Z \oplus T$ where $T$ is torsion prime to $6$ and 
$H_1(\partial M) \to H_1(M)/T \cong \mathbb Z$ is surjective. Suppose as well 
that  \\ 
\indent $(a)$ there is a slope $\alpha_0$ on $\partial M$ such that 
$\pi_1(M(\alpha_0))$ admits no homomorphism onto a \\ \indent \hspace{5mm} 
non-elementary Kleinian group or a Euclidean triangle group, and  \\ 
\indent $(b)$ either \\
\indent \hspace{5mm} $(i)$ there is no discrete, non-elementary representation $\rho \in R_{PSL_2(\mathbb R)}(M)$ such that \\ \indent \hspace{10mm} $\rho(\pi_1(M))$ is isomorphic to a free product of two non-trivial cyclic groups and  \\ \indent \hspace{10mm} $\rho(\pi_1(\partial M))$ is parabolic, or \\ \indent \hspace{5mm} $(ii)$ $T = \{0\}$, $M$ is  minimal and there is no representation 
$\rho: \pi_1(M(\lambda_M)) \to$  \\ \indent \hspace{10mm} $PSL_2(\mathbb R)$ such that 
$\rho(\pi_1(M(\lambda_M)))$ is a free product of two 
non-trivial  \\ \indent \hspace{10mm} cyclic groups and $\rho(\pi_1(\partial M))$ is parabolic. \\ 
Then there are infinitely many slopes $\alpha$ on $\partial M$ such that 
$M(\alpha)$ is minimal.  
\end{thm} 
When $T \neq \{0\}$, this corollary applies to the exterior of many knots in 
lens spaces. For instance, to the knot manifold obtained by $m$-Dehn surgery 
on one component of the right-hand Whitehead link, $m$ odd. See Example 
\ref{lensspace}. When $T = \{0\}$ it applies to the exterior of many knots in 
the $3$-sphere. 
For instance it is remarked in Example \ref{s3ginfiniteegs} that if the 
Alexander polynomial of a knot manifold $M$ with $H_1(M) \cong \mathbb Z$ is 
not divisible by the Alexander polynomial of a non-trivial torus knot, there 
is no homomorphism of $\pi_1(M)$ onto the free product of two non-trivial 
finite cyclic groups. Thus if $M$ is minimal, small, and hyperbolic, there are 
infinitely many slopes $\alpha$ on $\partial M$ such that $M(\alpha)$ is 
minimal. The proof of the following corollary is also discussed in 
this example.

\begin{cor} \label{intro:2bridgefillingminimal}
Let $M$ be the exterior of a $\frac{p}{q}$ two-bridge knot with $p$ is prime and $q \not \equiv \pm 1$ $($mod $p$$)$, or of a  $(-2, 3, n)$ pretzel with $n \not \equiv 0$ $($mod $3$$)$. Then there are infinitely many slopes $\alpha$ on $\partial M$ such that $M(\alpha)$ is minimal. 
\end{cor}

\subsection{Organization of the paper and acknowledgements} 

The basic properties of $PSL_2(\mathbb C)$-character varieties and Culler-Shalen theory are described in \S \ref{varieties}. The main result of this section states that the morphism $\varphi^*: X_{PSL_2}(\Gamma_2) \to X_{PSL_2}(\Gamma_1)$ induced by a virtual epimorphism $\varphi: \Gamma_1 \to \Gamma_2$ is generically 1-1 when restricted to the union of the positive dimensional algebraic components of $X_{PSL_2}(\Gamma_2)$ which contain a strictly irreducible character (Corollary \ref{1-1}). Section 3 deals with virtual epimorphisms of the fundamental groups of knot manifolds and contains the proofs of Theorems \ref{domboundintro}, \ref{intro:domseqbound}, and \ref{twobridgedomboundintro}. In \S \ref{sec:characters} we begin our study of families of discrete characters. In particular, unbounded sequences of such characters are studied in \S \ref{unbounded} and convergent ones in \S \ref{convergent}. This leads to the proofs of Theorems \ref{thm:characters} and \ref{thm:domination}. We discuss ${\cal H}$-minimal Dehn filling in \S \ref{sec:h-minimal}, including the proof of Theorem \ref{thm:h-minimal}. The last two sections of the paper deal with families of discrete $PSL_2(\mathbb R)$-characters (\S \ref{psl2rcharacters}) and constructing minmal manifolds (\S \ref{gminimal}). The proofs of Theorems \ref{cor:s-dominationintro}, \ref{finitesl2dominationintro}, \ref{uinfiniteintro} and Corollaries \ref{prop:g-minimalintro}, \ref{intro:2bridgefillingminimal} are found here. Finally, three appendices are included which deal, respectively, with the following topics: a smoothness criterion for dihedral characters; peripheral values of representations of the fundamental groups of twist knot exteriors; the bending construction. 

The authors wish to thank Shicheng Wang for many stimulating conversations 
and for bringing to their attention his question about Dehn fillings and non-zero degree maps.

\section{Varieties of $PSL_2(\mathbb C)$-characters} \label{varieties}

\subsection{Generalities} \label{generalities}

In what follows we shall refer to the elements of $PSL_2(\mathbb C)$
as matrices. Denote by ${\cal D}$ the abelian subgroup of $PSL_2(\mathbb C)$ consisting of diagonal matrices and by ${\cal N}$ the subgroup consisting of those matrices which are either
diagonal or have diagonal coefficients $0$. Note that ${\cal D}$ has index $2$ in ${\cal N}$ and any element in ${\cal N} \setminus {\cal D}$ has order $2$. Further, the centre of ${\cal N}$, which we will denote by $Z({\cal N})$, is isomorphic to $\mathbb Z/2$ and is generated by  $\pm \left(\begin{array}{cc} i & 0 \\ 0 & -i \end{array} \right)$. 

The action of $SL_2(\mathbb C)$ on $\mathbb C^2$ descends to one of $PSL_2(\mathbb C)$ on $\mathbb CP^1$. We call a representation $\rho$ with values in $PSL_2(\mathbb C)$ {\it irreducible} if the associated action on $\mathbb CP^1$ is fixed point free, otherwise we call
it {\it reducible}. We call it {\it strictly irreducible} if the action has no invariant subset in $\mathbb CP^1$ with fewer than three points. Note that 
\vspace{-.35cm} 
\begin{itemize}
\item $\rho$ is reducible if and only if it is conjugate to a representation
whose image consists of upper-triangular
matrices.
\vspace{-.15cm} 
\item $\rho$ is conjugate to a representation with image in ${\cal D}$ if
and only if the action on $\mathbb
CP^1$ has at least two fixed points. It is conjugate into ${\cal
N}$ if and only if it leaves a two point subset of $\mathbb CP^1$ invariant.
\vspace{-.15cm}
\item $\rho$ is is strictly irreducible if and only if it is irreducible but not conjugate into ${\cal N}$.  
\vspace{-.15cm} 
\item if $\rho$ is irreducible and $A \in PSL_2(\mathbb C)$ satisfies $A \rho A^{-1}
= \rho$, then either $A = \pm I$ or up to conjugation, $A = \pm \left(\begin{array}{cc} i & 0 \\ 0 & -i \end{array} \right)$ and $\rho$ conjugates into ${\cal N}$. Thus if $\rho$ is strictly irreducible, then $A = \pm I$.
\end{itemize}
\vspace{-.25cm} 

The action of $PSL_2(\mathbb C)$ on $\mathbb CP^1 = S^2_\infty$ extends over $\mathbb H^3$ yielding an identification $PSL_2(\mathbb C) = \hbox{Isom}_+(\mathbb H^3)$. A representation is called {\it elementary} if the associated action on $\overline{\mathbb H}^3$ has a finite orbit. Equivalently, the representation is reducible or conjugates to one with image in either $SO(3) = PSU(2)$ or ${\cal N}$. 

Let $\Gamma$ be a finitely generated group. The set $R_{PSL_2}(\Gamma)$ of representations of $\Gamma$ with values in $PSL_2(\mathbb C)$ admits the structure of a $\mathbb C$-affine algebraic set [LM] called the {\it $PSL_2(\mathbb C)$-representation variety} of $\Gamma$. The action of $PSL_2(\mathbb C)$ on $R_{PSL_2}(\Gamma)$ determines an algebro-geometric quotient $X_{PSL_2}(\Gamma)$ whose coordinate ring is $\mathbb C[R_{PSL_2}(\Gamma)]^{PSL_2(\mathbb C)}$ and a regular map $t: R_{PSL_2}(\Gamma) \to X_{PSL_2}(\Gamma)$ [LM]. This quotient is called the {\it $PSL_2(\mathbb C)$-character variety} of $\Gamma$. For $\rho \in R_{PSL_2}(\Gamma)$, we denote $t(\rho)$ by $\chi_\rho$ and refer to it as the {\it character} of $\rho$. If $\chi_{\rho_1} = \chi_{\rho_2}$ and $\rho_1$ is irreducible, then $\rho_1$ and $\rho_2$ are conjugate representations. We can therefore call a character $\chi_\rho$ reducible, irreducible, or strictly irreducible if $\rho$ has that property. Each reducible character is the character of a diagonal representation, that is, one with image in ${\cal D}$. The property of an irreducible representation being conjugate into $SO(3)$ is also determined by its character (see Proposition III.1.1 of \cite{MS1} for instance) and so if $\chi_{\rho_1} = \chi_{\rho_2}$ and $\rho_1$ is elementary, then so is $\rho_2$. In this case we call the character elementary. 

When $\Gamma$ is the fundamental group of a path-connected space $Y$, we write $R_{PSL_2}(Y)$ for $R_{PSL_2}(\pi_1(Y))$, $X_{PSL_2}(Y)$ for $X_{PSL_2}(\pi_1(Y))$, and refer to them respectively as the $PSL_2(\mathbb C)$-representation variety of $Y$ and $PSL_2(\mathbb C)$-character variety of $Y$. 

Each $\gamma \in \Gamma$ determines an element $f_\gamma$ of the
coordinate ring $\mathbb C [X_{PSL_2}(\Gamma)]$ according to the formula 
$$f_\gamma(\chi_\rho) = \mbox{trace}(\rho(\gamma))^2 - 4.$$

A homomorphism $\varphi: \Gamma_1 \to \Gamma_2$ determines morphisms $\varphi^*: R_{PSL_2}(\Gamma_2) \to R_{PSL_2}(\Gamma_1), \; \rho \mapsto \rho \circ \varphi$ 
and $\varphi^*: X_{PSL_2}(\Gamma_2) \to X_{PSL_2}(\Gamma_1), \; \chi_\rho \mapsto \chi_{\rho \circ \varphi}$. For $\gamma \in \Gamma_1$ and $\chi_\rho \in X_{PSL_2}(\Gamma_2)$ we have 
$$f_\gamma(\varphi^*(\chi_\rho)) = f_{\varphi(\gamma)}(\chi_\rho). \eqno{(2.1.1)}$$

We end this section with a useful observation

\begin{lemma} \label{closed} 
If the image of $\varphi: \Gamma_1 \to \Gamma_2$ is of finite index $n$ in $\Gamma_2$, then $\varphi^*:  X_{PSL_2}(\Gamma_2) \to X_{PSL_2}(\Gamma_1)$ is a closed map with respect to the Zariski topology. 
\end{lemma}

\pf Let $X_0$ be a Zariski closed subset of $X_{PSL_2}(\Gamma_2)$ and let $Y_0 = \overline{\varphi^*(X_0)}$. If $\bar X_0, \bar Y_0$ are projective closures of $X_0, Y_0$, then $\varphi^*$ determines a surjective projective morphism $\bar \varphi^*: \bar X_0 \to \bar Y_0$. Let $y_0 \in Y_0$ and choose $x_0 \in \bar X_0$ such that $\bar \varphi^*(x_0) = y_0$, and a projective curve $C \subseteq \bar X_0$ which contains $x_0$. Set $C_0 = C \cap C_0$ and note that if $x_0 \not \in C_0$, there is some $\gamma \in \pi_1(M)$ such that $f_\gamma(x_0) = \infty$ (cf. Theorem 2.1.1 of \cite{CS}). For $A, B \in SL_2(\mathbb C)$ we have $\hbox{trace}(AB) + \hbox{trace}(A^{-1} B) = \hbox{trace}(A) \hbox{trace}(B)$, and this identity can be used inductively to show that $f_{\gamma^n}$ is a degree $|n|$ polynomial in $f_{\gamma}$. In particular, $f_{\gamma^n}(x_0) = \infty$. On the other hand there is some $\delta \in \Gamma_1$ such that $\varphi(\delta) = \gamma^n$. Then $f_\delta(y_0) = f_\delta(\bar \varphi^*(x_0)) = f_{\varphi(\delta)}(x_0) = f_{\gamma^n}(x_0) = \infty$. But this contradicts the fact that $y_0 \in Y_0$. Thus $x_0 \in C_0 \subseteq X_0$ and so $\varphi^*$ is onto $Y_0$. 
\qed 

\subsection{Subvarieties of $X_{PSL_2}(\Gamma)$} \label{subvarieties} 

The set of reducible representations $R_{PSL_2}^{red}(\Gamma) \subseteq R_{PSL_2}(\Gamma)$ ($\Gamma$ a finitely generated group) is a closed algebraic subset (cf. the proof of Corollary 1.4.5 of \cite{CS}). The sets $R_{SO(3)}(\Gamma), R_{{\cal D}}(\Gamma)$, and $R_{{\cal N}}(\Gamma)$ of representations of $\Gamma$ with values in $SO(3), {\cal D}$, and ${\cal N}$ are also closed algebraic subsets of $R_{PSL_2}(\Gamma)$. A similar statement holds for their images $X_{SO(3)}(\Gamma), X_{PSL_2}^{red}(\Gamma)$, and $X_{{\cal N}}(\Gamma)$ in $X_{PSL_2}(\Gamma)$. In particular, the set $X_{Elem}(\Gamma) = X_{SO(3)}(\Gamma) \cup X_{{\cal N}}(\Gamma)$ of characters is Zariski closed in $X_{PSL_2}(\Gamma)$.

A subvariety $X_0$ of $X_{PSL_2}(\Gamma)$ is called {\it non-trivial} if it contains the
character of an irreducible representation. It is called {\it strictly non-trivial} if it contains the character of a strictly irreducible representation. The property of being (strictly) irreducible is open so that the generic character of a (strictly) non-trivial subvariety of $X_{PSL_2}(\Gamma)$ is (strictly) irreducible. Let $X_+^{irr}(\Gamma)$ denote the union of the positive dimensional non-trivial components of $X_{PSL_2}(\Gamma)$ and $X_+^{str}(\Gamma)$ the union of its positive dimensional strictly non-trivial components.

For each non-trivial subvariety $X_0$ of $X_{PSL_2}(\Gamma)$ there is a subvariety $R_{X_0}$ of $R_{PSL_2}(\Gamma)$ uniquely determined by the condition that it is conjugation invariant and $t(R_{X_0}) = X_0$ (cf. Lemma 4.1 of [BZ1]). We define the \emph{kernel} of $X_0$ to be the normal subgroup of $\Gamma$ given by  
$$\hbox{Ker}(X_0) = \bigcap \limits_{\rho \in R_{X_0}} \hbox{kernel}(\rho).$$
For instance $\hbox{Ker}(X_0) = \{1\}$ if $R_{X_0}$ contains an injective representation. 

\begin{lemma}  \label{ker} 
Let $X_0$ be a non-trivial subvariety of $X_{PSL_2}(\Gamma)$. \\ 
$(1)$ There is a subset $V$ of $R_{X_0}$ which is a countable union of proper, closed, conjugation invariant algebraic subsets of $R_{X_0}$ such that for $\rho \in R_{X_0} \setminus V$, $\hbox{kernel}(\rho) = \hbox{Ker}(X_0)$. \\
$(2)$ If $\varphi: \Gamma_1 \to \Gamma_2$ is a homomorphism and $X_0$ is a subvariety of $X_{PSL_2}(\Gamma_2)$ such that $Y_0 = \overline{\varphi^*(X_0)}$ is non-trivial, then $\hbox{Ker}(Y_0) = \varphi^{-1}(\hbox{Ker}(X_0))$. In particular, $\hbox{kernel}(\varphi) \subseteq \hbox{Ker}(Y_0)$.  
\end{lemma}

\pf (1) For each $\gamma \in \pi_1(M)$ set $V_\gamma = \{\rho \in R_{X_0}
\; | \; \rho(\gamma) = \pm I\}$. Then  $V_\gamma$ is a closed,
conjugation invariant algebraic subset of $R_{X_0}$. It is clear that
$\gamma \in \hbox{Ker}(X_0)$ if and only if $V_\gamma  = R_{X_0}$.  Set $V = \bigcup_{\gamma \not \in \hbox{Ker}(X_0)} V_\gamma$ and observe that $\rho \in R_{X_0} \setminus V$ if and only if $\rho(\gamma) \ne \pm I$
for each $\gamma \not \in \hbox{Ker}(X_0)$. In particular, $\hbox{kernel}(\rho) =
\hbox{Ker}(X_0)$ for such $\rho$. This proves (1). 

(2) Now $\varphi^*(X_0) = t(\varphi^*(R_{X_0})) \subseteq t(\overline{\varphi^*(R_{X_0})}) \subseteq \overline{t(\varphi^*(R_{X_0}))} = \overline{\varphi^*(X_0)} = Y_0$ and since $\overline{\varphi^*(R_{X_0})}$ is closed and conjugation invariant in $R_{PSL_2}(\Gamma_1)$, Theorem 3.3.5(iv) of \cite{Ne} implies that $t(\overline{\varphi^*(R_{X_0})})$ is Zariski closed in $X_{PSL_2}(\Gamma_1)$. It follows that $R_{Y_0} = \overline{\varphi^*(R_{X_0})}$. Hence noting that $\varphi^*(\rho)(\gamma) = \rho(\varphi(\gamma)) = \pm I$ whenever $\gamma \in \varphi^{-1}(\hbox{Ker}(X_0))$ and $\rho \in R_{X_0}$, it follows that $\rho'(\gamma) = \pm I$ for all $\rho' \in R_{Y_0}$. In other words, $\gamma \in \hbox{Ker}(Y_0)$. Conversely if $\gamma \in \hbox{Ker}(Y_0)$ and $\rho \in R_{X_0}$, then $\rho(\varphi(\gamma)) = \varphi^*(\rho)(\gamma) = \pm I$. Thus $\gamma \in \varphi^{-1}(\hbox{Ker}(X_0))$. This proves (2). 
\qed

We call a component $X_0$ of $X_+^{irr}(\Gamma)$ {\it principal} if it contains the character of a discrete, faithful,  irreducible representation of $\Gamma/Z(\Gamma)$ where $Z(\Gamma)$ denotes the centre of $\Gamma$. It is clear that $\hbox{Ker}(X_0) \subseteq Z(\Gamma)$. 

\begin{lemma} \label{kerprincipal} $\;$ \\
$(1)$ If $X_0$ is a principal component of $X_+^{irr}(\Gamma)$, then $\hbox{Ker}(X_0) = Z(\Gamma)$. \\
$(2)$ If $\varphi: \Gamma_1 \to \Gamma_2$ is a homomorphism and $X_0$ is a subvariety of $X_{PSL_2}(\Gamma_2)$ such that $\overline{\varphi^*(X_0)}$ is principal, then $\hbox{kernel}(\varphi) \subseteq Z(\Gamma_1)$. 
\end{lemma}

\pf (1) It suffices to show that $Z(\Gamma) \subseteq \hbox{Ker}(X_0)$. To that end we note that if $\rho \in R_{PSL_2}(\Gamma)$ is irreducible, then every element in $\rho(Z(\Gamma))$ has order $1$ or $2$. In particular for $\gamma \in Z(\Gamma)$ and $\rho \in R_{X_0}$ we have $f_\gamma(\chi_\rho) \in \{0, -4\}$. Hence $f_\gamma|X_0$ is constant and since it vanishes at a discrete faithful character, it is identically zero. Thus $\rho(\gamma) = \pm I$ for all $\rho \in R_{X_0}$ and therefore
$Z(\Gamma) \subseteq \hbox{Ker}(X_0)$. 

Part (2) follows from part (1) and part (2) of the previous lemma. 
\qed 

\subsection{Restriction} \label{restriction} 

If $\varphi: \Gamma_1 \to \Gamma$ is surjective, it is easy to see that $\varphi^*: X_{PSL_2}(\Gamma) \to X_{PSL_2}(\Gamma_1)$ is injective. The goal of this section is to show that a similar conclusion is true for virtual epimorphisms $\varphi: \Gamma_1 \to \Gamma$ as long as we restrict $\varphi^*$ to $X_{+}^{str}(\Gamma)$. 

Let $D_n, T_{12}, O_{24}$ denote, respectively, the dihedral group of order $2n$, the tetrahedral group of order $12$, and the octahedral group of order $24$. Set 
$${\cal K} = \{ \pm I, \pm \left(\begin{array}{cc} i & 0 \\ 0 & -i \end{array} \right), \pm \left(\begin{array}{cc} 0 & 1 \\ -1 & 0 \end{array} \right), \pm \left(\begin{array}{cc} 0 & i \\ i & 0 \end{array} \right)\}  \subset {\cal N}$$ 
and observe that ${\cal K} \cong D_2$. It is well-known and easy to verify that the centraliser $Z_{PSL_2}({\cal K})$ of ${\cal K}$ in $PSL_2(\mathbb C)$ is ${\cal K}$ and its normaliser $N_{PSL_2}({\cal K})$ is isomorphic to $O_{24}$. The only other subgroups of $PSL_2(\mathbb C)$ which contain ${\cal K}$ as a normal subgroup are the (unique) subgroup of ${\cal N}$ isomorphic to $D_4$ and the (unique) index $2$ subgroup of $N_{PSL_2}({\cal K})$, which is isomorphic to $T_{12}$. 

Given any homomorphisms $\rho \in R_{\cal N}(\Gamma)$ and $\varphi: \Gamma \to Z({\cal N})$, the equation
$$\rho'(\gamma) = \varphi(\gamma) \rho(\gamma)$$ 
defines an element $\rho' \in R_{{\cal N}}(\Gamma)$ which will denote by $\varphi \cdot \rho$. 

\begin{lemma} \label{determined}
Let $\Gamma_0$ be a normal subgroup of a finitely generated group $\Gamma$ and suppose that $\rho_1, \rho_2 \in R_{PSL_2}(\Gamma)$ restrict to the same irreducible
representation $\rho_0 \in R_{PSL_2}(\Gamma_0)$.  Then either \\ 
$(a)$ $\rho_1 = \rho_2$, or \\
$(b)$ after a similtaneous conjugation of $\rho_1$ and $\rho_2$ we have $\rho_j(\Gamma) \subset {\cal N}$ for $j = 1, 2$ and there is a homomorphism $\varphi:  \Gamma/\Gamma_0 \to Z({\cal N})$ such that $\rho_2 = \varphi \cdot \rho_1$, or \\
$(c)$ $\rho_0(\Gamma_0) \cong D_2$ and $\rho_1(\Gamma) = \rho_2(\Gamma)$ is isomorphic to one of $D_2, D_4, T_{12},$ or $O_{24}$. There are only finitely many orbits in $R_{PSL_2}(\Gamma)$ for which this case arises. 
\end{lemma}

\pf Fix $\gamma \in \Gamma$ and set $\rho_j(\gamma) = A_j$ for $j = 1, 2$. Then for $\gamma_0
\in \Gamma_0$ we have $A_1
\rho_0(\gamma_0) A_1^{-1} = \rho_0(\gamma \gamma_0 \gamma^{-1}) = A_2
\rho_0(\gamma_0) A_2^{-1}$. Thus $A_2^{-1} A_1 \in Z_{PSL_2}(\rho_0(\Gamma_0))$, the centraliser of $\rho_0(\Gamma_0)$ in $PSL_2(\mathbb C)$. The irreducibility of $\rho_0$ then implies that either $A_2 = A_1$ or, after conjugation, $\rho_0(\Gamma_0) \subset {\cal N}$ and $A_2^{-1} A_1 = \pm \left(\begin{array}{cc} i & 0 \\ 0 & -i \end{array} \right)$. If the former occurs for each $\gamma \in \Gamma$ then $\rho_1 = \rho_2$ and so we are in case (a). Suppose then that the latter arises for some $\gamma \in \Gamma$. 

Similtaneously conjugate $\rho_1$ and $\rho_2$ so that $\rho_0(\Gamma_0) \subset {\cal N}$. If some $B \in \rho_0(\Gamma_0)$ has order larger than $2$, then $B \in {\cal D}$ and for any matrix $A \in \rho_j(\Gamma)$ (either $j$) we have $A B A^{-1} \in {\cal D}$. Thus $A B A^{-1} = B^\epsilon$ for some $\epsilon = \pm 1$ and therefore $A \in {\cal N}$. It follows that $\rho_j(\Gamma) \subset {\cal N}$ for both $j$. From the previous paragraph we know that if $\gamma \in \Gamma$ and $A_j = \rho_j(\gamma)$, then 
$\pm (A_2^{-1} A_1)^2 = \pm I$. Moreover, $\pm (A_2^{-1} A_1) B (A_2^{-1} A_1)^{-1} = B$ and our restrictions on $B$ then imply that $A_2^{-1} A_1 \in {\cal D}$. Thus $\rho_2(\gamma)^{-1}\rho_1(\gamma) \in \{\pm I, \pm \left(\begin{array}{cc} i & 0 \\ 0 & -i \end{array} \right)\} = Z_{PSL_2}({\cal N})$. It follows that $\gamma \mapsto \rho_1(\gamma)^{-1}\rho_2(\gamma)$ induces a homomorphism $\varphi: \Gamma/\Gamma_0 \to Z_{PSL_2}({\cal N})$ such that $\rho_2 = \varphi \cdot \rho_1$. This is case (b). 

Suppose then that each $\pm I \ne B \in \rho_0(\Gamma_0) \subset {\cal N}$ has order $2$. Then $\rho_0(\Gamma_0) = {\cal K}$ and is a normal subgroup of  $\rho_j(\Gamma)$ ($j = 1, 2$). Hence by our remarks above we know that $\rho_j(\Gamma)$ is a subgroup of $N_{PSL_2}({\cal K}) \cong O_{24}$ and is isomorphic to one of $D_2, D_4, T_{12}$, and $O_{24}$. Conjugation induces an exact sequence $1 \to {\cal K} \to \rho_j(\Gamma) \to Aut({\cal K}) \cong S_3$ and we saw above that for each $\gamma \in \Gamma$, $\rho_1(\gamma)^{-1} \rho_2(\gamma) \in Z_{PSL_2}({\cal N}) \subset {\cal K}$. Therefore the images of $\rho_1(\Gamma)$ and $\rho_2(\Gamma)$ in $Aut({\cal K})$ coincide and so $|\rho_1(\Gamma)| = |\rho_2(\Gamma)|$. Since ${\cal K} \subseteq \rho_j(\Gamma)$, we actually have $\rho_1(\Gamma) = \rho_2(\Gamma)$. Finally, since two finite subgroups of $PSL_2(\mathbb C)$ which are abstractly isomorphic are conjugate in $PSL_2(\mathbb C)$ and $\Gamma$ is finitely generated, there are only finitely many orbits in $R_{PSL_2}(\Gamma)$ of representations with image either $D_2, D_4, T_{12},$ or $O_{24}$. This is case (c). 
\qed

\begin{cor} \label{str11}
Let $\Gamma_0$ be a finitely generated normal subgroup of a finitely generated group $\Gamma$ and suppose that $\rho_1, \rho_2 \in R_{PSL_2}(\Gamma)$ are strictly irreducible with images different from $T_{12}$ and $O_{24}$. Suppose further that $\chi_{\rho_1|\Gamma_0} = \chi_{\rho_2|\Gamma_0}$ and is irreducible. Then $\chi_{\rho_1} = \chi_{\rho_2}$. 
\qed 
\end{cor}

\begin{prop} \label{inj}
Let  $X_0$ be a positive dimensional non-trivial subvariety of $X_{PSL_2}(\Gamma)$ and $\Gamma_0$ a finitely generated normal subgroup of $\Gamma$. Then one of the following three situations arises. \\
$(a)$ $R_{X_0} \subset R_{PSL_2}(\Gamma/\Gamma_0)$. That is,
$\rho(\Gamma_0) = \{\pm I\}$
for each $\rho \in R_{X_0}$. \\
$(b)$ The restriction to $\Gamma_0$ of each $\rho \in R_{X_0}$ is
conjugate into ${\cal N}$ and
$X_0 \subset X_{{\cal N}}(\Gamma)$. \\
$(c)$ The restriction to $\Gamma_0$ of the generic $\rho \in R_{X_0}$ is
strictly irreducible.
Moreover, there is a Zariski open, conjugation invariant, connected subset $U$ of $R_{X_0}$ such that if $\rho_1, \rho_2 \in U$ and their restrictions to $\Gamma_0$ have
the same characters, then $\chi_{\rho_1} = \chi_{\rho_2}$.
\end{prop}

\pf We shall suppose that conclusion (a) of the lemma does not occur.
In particular, there is some $\rho_0 \in R_{X_0}$ such that $\rho_0(\Gamma_0) \ne \{\pm I\}$. Then there is a Zariski open, conjugation invariant, connected subset $U_0 \subset R_{X_0}$ consisting of irreducible representations such that $\rho(\Gamma_0) \ne \{\pm I\}$ for all $\rho \in U_0$. As there are only finitely many conjugacy classes of representations with image isomorphic to $D_2$, we also suppose that $\rho(\Gamma) \not \cong D_2$ for $\rho \in U_0$.

As a first case, assume that the restriction to $\Gamma_0$ of each
$\rho \in U_0$ is reducible. If $\rho(\Gamma_0)$ is not diagonalisable
for some $\rho \in U_0$, it has a
unique fixed point in $\mathbb CP^1$, and so the fact that
$\rho(\Gamma_0)$ is normal in $\rho(\Gamma)$
implies that $\rho$ is reducible, which is not the case. Thus
$\rho(\Gamma_0)$ is diagonalisable for $\rho \in
U_0$. It follows that  for such $\rho$, $\rho(\Gamma)$ leaves invariant a two element subset of $\mathbb C P^1$ and therefore conjugates to a representation with image in ${\cal N}$. Since we have assumed that $\rho(\Gamma) \not \cong D_2$, there is a unique index $2$ subgroup $\Gamma_\rho$ of $\Gamma$ such that $\rho|\Gamma_\rho$ is diagonalisable while $\rho(\Gamma_\rho) \not \subset \mathbb Z/2$. Fix generators $\gamma_0, \gamma_1, \ldots , \gamma_n$ of $\Gamma$ so that $\rho(\gamma_j)^2 = \pm I$ if and only if $j = 0$. Then $\Gamma_\rho$ is generated by 
$\gamma_0^2, \gamma_1, \ldots , \gamma_n, \gamma_1', \ldots , \gamma_n'$ where $\gamma_j' = \gamma_0 \gamma_j \gamma_0^{-1}$. For all $\rho' \in U_0$ close enough to $\rho$ and $j \geq 1$ we have $\rho'(\gamma_j)^2 \ne \pm I$ and therefore $\rho'(\gamma_j) \in {\cal D}$. Hence $\Gamma_\rho \subseteq \Gamma_{\rho'} $ for $\rho'$ in an open neighbourhood of $\rho$. It follows that $\Gamma_\rho$ is independent of $\rho \in U_0$. Denote this common subgroup by $\Gamma_1$. Each $\rho_0 \in R_{X_0}$ is a limit of representations $\rho_n \in U_0$. It follows that $\rho_0|\Gamma_1$ is reducible and as above we see that $\rho_0$ is either reducible or conjugates into ${\cal N}$. In either case, $\chi_{\rho_0} \in X_{{\cal N}}(\Gamma)$. Thus conclusion (b) holds.

Next assume that the restriction to $\Gamma_0$ of the generic $\rho \in R_{X_0}$
is irreducible but not strictly irreducible. For such $\rho$, after a conjugation we may suppose that $\rho(\Gamma_0) \subset {\cal N}$. If $\rho(\Gamma_0) \ne {\cal K}$ then $\rho(\Gamma) \subset {\cal N}$, while if $\rho(\Gamma_0) = {\cal K}$ then $\rho(\Gamma)$ is isomorphic to one of ${\cal K}, D_4, T_{12}$ or $O_{24}$ (c.f. the proof of Lemma \ref{determined}). There are only finitely many characters of representations $\rho \in R_{PSL_2}(\Gamma)$ with image $T_{12}$ or $O_{24}$, so the generic character in $X_0$ lies in $X_{\cal N}(\Gamma)$. It follows that conclusion (b) holds.

Finally suppose that the restriction to $\Gamma_0$ of some $\rho \in R_{X_0}$
is strictly irreducible. Then there is a Zariski open subset $U \subseteq U_0$ such that for each $\rho \in U$, $\rho|\Gamma_0$ is strictly irreducible and the image of $\rho$ is neither $T_{12}$ nor $O_{24}$. If $\rho_1, \rho_2 \in U$ restrict to representations of $\Gamma_0$ with the same character, then Lemma \ref{determined} shows that $\chi_{\rho_1} = \chi_{\rho_2}$. This is conclusion (c).
\qed

\begin{cor} \label{1-1}
Let $\varphi: \Gamma_1 \to \Gamma$ be a virtual epimorphism. Then $Y := \varphi^*(X_+^{str}(\Gamma))$ is a Zariski closed subset of $X_+^{str}(\Gamma_1)$. Further, $\varphi^*$ sends distinct components of $X_+^{str}(\Gamma)$ to Zariski dense subsets of distinct components of $Y$ and is generically one-to-one on each of these components.  
\end{cor}

\pf It suffices to prove the result when $\Gamma_1$ is a finite index subgroup of $\Gamma$ and $\varphi$ is the inclusion. 

Fix a finite index subgroup $\Gamma_0 \subseteq \Gamma_1$ which is normal in
$\Gamma$ and a component $X_0$ of  $X_+^{str}(\Gamma)$. Since $X_0$ is positive dimensional, it cannot be contained in $X_{PSL_2}(\Gamma/ \Gamma_0)$. Thus the third option in Proposition \ref{inj} must arise with respect to the restriction map $R_{PSL_2}(\Gamma) \to R_{PSL_2}(\Gamma_0)$. In particular, the induced map $X_0 \to X_{PSL_2}(\Gamma_0)$ is generically $1-1$ and so its image is a positive dimensional subvariety of $X_{PSL_2}^{str}(\Gamma_0)$. This image is also closed by 
Lemma \ref{closed}. Corollary \ref{str11} implies that distinct components of $X_+^{str}(\Gamma)$ are sent to distinct subvarieties of $X_{PSL_2}^{str}(\Gamma_0)$. Finally, since the restriction $X_+^{str}(\Gamma) \to X_{PSL_2}(\Gamma_0)$ factors through the map $X_+^{str}(\Gamma) \to X_{PSL_2}(\Gamma_1)$, the conclusions of the corollary hold. 
\qed

\subsection{Culler-Shalen theory}\label{culler-shalen}

In this section, $M$ will denote compact, connected, orientable,
irreducible $3$-manifold whose boundary is a torus. 

Any complex affine curve $C$ admits an affine desingularisation
$C^\nu \stackrel{\nu}{\longrightarrow} C$ where $\nu$ is surjective and
regular. Moreover, the smooth
projective model $\tilde C$ of $C$ is obtained by adding a finite number of
ideal points to $C^\nu$. Thus
$\tilde C = C^\nu \cup {\cal I}(C)$ where ${\cal I}(C)$ is the set of ideal
points of $C$. There are natural
identifications between the function fields of $C, C^\nu,$ and $\tilde C$.
Thus to each $f \in \mathbb C(C)$ we
have corresponding $f^\nu \in \mathbb C(C^\nu) = \mathbb
C(C)$ and $\tilde f \in \mathbb C(\tilde
C) = \mathbb C(C)$ where $f^\nu = f \circ \nu = \tilde f|C^\nu$. 

The proof of the
following basic result can be found in
Proposition 3.2.1 of \cite{CS}.

\begin{prop} {{\rm (Thurston)}} \label{posdim}
Let $V$ be a compact orientable $3$-manifold and $\rho \in R_{PSL_2}(V)$ an irreducible representation such that $\rho(\pi_1(T)) \ne \{\pm I\}$ for each toral boundary component $T$ of $V$. Then the dimension of any component $X_0$ of $X_{PSL_2}(V)$ which contains $\chi_{\rho}$ is at least $3t - \frac{3}{2}\chi(\partial V)$ where $t$ is the number of toral boundary components of $V$. 
\qed
\end{prop}

Recall that each $\gamma \in \pi_1(M)$ determines an element $f_\gamma$ of the
coordinate ring $\mathbb C [X_{PSL_2}(M)]$ satisfying
$$f_\gamma(\chi_\rho) = \mbox{trace}(\rho(\gamma))^2 - 4$$
where $\rho \in R_{PSL_2}(M)$. Each $\alpha \in H_1(\partial M) = \pi_1(\partial M)$ defines an element of $\pi_1(M)$ well-defined up to conjugation and therefore determines an element $f_\alpha \in \mathbb C [X_{PSL_2}(M)]$. Similarly each slope $\alpha$ on $\partial M$ determines an element of $\pi_1(M)$ well-defined up to conjugation and taking inverse, and so defines $f_\alpha \in \mathbb C[X_{PSL_2}(M)]$.

To each curve $X_0$ in $X_{PSL_2}(M)$ we associate the function 
$$d_{X_0}: \pi_1(M) \to \mathbb Z, \; d_{X_0}(\gamma) =  \hbox{degree}(f_\gamma: X_0 \to \mathbb C). $$
Standard trace identities imply that for $n \in \mathbb Z$, 
$$d_{X_0}(\gamma^n) = |n| d_{X_0}(\gamma).$$
More generally, it was shown in\cite{CGLS} that $d_{X_0}$ has nice properties when restricted to abelian subgroups of $\pi_1(M)$. For instance, when restricted to 
$\pi_1(\partial M)$ it gives rise to a {\it Culler-Shalen seminorm} 
$$\| \cdot \|_{X_0}: H_1(\partial M; \mathbb R) \to [0, \infty)$$
where for each $\alpha \in H_1(\partial M) = \pi_1(\partial M)$
we have $\|\alpha\|_{X_0} =  d_{X_0}(\alpha)$. 

We say that a curve $X_0 \subset X_{PSL_2}(M)$ is a {\it norm curve} if $\|\cdot\|_{X_0}$ is a norm. If $\|\cdot\|_{X_0} \ne 0$, though it is not a norm, there is a primitive element $\beta \in H_1(\partial M)$ well-defined up to sign such that $\|\beta\|_{X_0} = 0$. In this case we say that $X_0$ is a {\it $\beta$-curve}. 

For $x \in \tilde X_0$ and $\gamma \in \pi_1(M)$, we denote by
$Z_x(\tilde f_\gamma), \Pi_x(\tilde f_\gamma)$ the multiplicity of $x$
as a zero, respectively pole, of $\tilde
f_\gamma$. From the definition of $\|\cdot\|_{X_0}$ we see that for each $\alpha \in H_1(\partial M)$ we have
$$\|\alpha\|_{X_0} = \sum_{x \in \tilde X_0} Z_x(\tilde f_\alpha) =
\sum_{x \in {\cal I}(X_0)} \Pi_x(\tilde f_\alpha). \eqno{(2.4.1)}$$
When $M$ is hyperbolic, there is an essentially canonical choice of curve
$X_0 \subset X_{PSL_2}(M)$ characterized by the fact that it contains the
character of a discrete, faithful representation. In this case, it is known \cite{CGLS} that $\|\cdot\|_{X_0}$ is a norm.

Consider a curve $X_0$ in $X_{PSL_2}(M)$. We say that a sequence of characters
$\chi_n \in X_0$ {\it converges to an ideal point} $x_0 \in \tilde X_0$ if
there are a sequence $\{x_n\}$ in $X_0^\nu \subset \tilde X_0$ and an ideal point $x_0 \in {\cal I}(X_0)$ such
that $\nu(x_n) = \chi_n$ for all $n$ and $\lim_n x_n = x_0$. 

For a path-connected space $X$, a representation $\rho \in R_{PSL_2}(X)$, a
path-connected
subspace $Q$ of a space $X$ with inclusion map $i: Q \to X$, set
$$\rho^Q := \rho\circ i_\#: \pi_1(Q) \to PSL_2(\mathbb C).$$
Since $\rho_Q$ is determined up to conjugation, there is a well-defined
$$\chi_\rho^Q = \chi_{\rho^Q}.$$

\begin{prop} \label{limitsreps} {\rm (\cite{CS})}
Suppose that $X_0$ is a curve in $X_{PSL_2}(M)$ and $\rho_n \in R_{X_0}
\subset R_{PSL_2}(M)$
is a sequence of representations whose characters $\chi_n$ converge to an
ideal point $x_0$ of $\tilde X_0$.
Then there is an essential surface $S \subset M$ whose complementary
components $A_1, A_2, \ldots , A_n$
satisfy the following properties. \\
$(a)$ For each $i$, the characters $\chi_n^{A_i}$ converge to a character
$\chi_0^{A_i}$.
Thus if $S_j$ is a component of $S$, then $\chi_0^{S_j} := \lim_n
\chi_n^{S_j} \in X_{PSL_2}(S_j)$ exists.
Further, $\chi_0^{S_j}$ is reducible. \\
$(b)$ For each $i$, there are conjugates $\sigma_n^{A_i}$ of $\rho_n^{A_i}$
which converge
to a representation $\sigma_0^{A_i} \in R_{PSL_2}(A_i)$ for which
$\chi_{\sigma_0^{A_i}} = \chi_0^{A_i}$.
\qed
\end{prop}

A representation $\sigma_0^{A_i} \in R_{PSL_2}(A_i)$ obtained as a limit of
some conjugates of
$\rho_n^{A_i}$ is said to be a {\it limiting representation associated to
the sequence $\{\rho_n\}$}.

Any essential surface $S \subset M$ as described in Proposition \ref{limitsreps} is said to be {\it associated} to the ideal point $x_0$. 

\begin{prop} \label{idealvalue} {\rm (\cite{CS}, \cite{CGLS}, \cite{CCGLS})}
Let $x_0$ be an ideal point of a curve $X_0$ in $X_{PSL_2}(M)$.  There is at least one primitive class $\alpha \in H_1(\partial M)$ such that $\tilde f_\alpha(x_0) \in \mathbb C$. Further, \\
$(1)$ if there is exactly one such class (up to sign), then it is
a boundary class and any surface $S$ associated to $x_0$ has non-empty boundary of
slope $\alpha$. Further, $\tilde f_\alpha(x_0) = (\lambda - \lambda^{-1})^2$ where $\lambda$ is a root of unity,  \\
$(2)$ if there are rationally independent classes $\alpha, \beta \in
H_1(\partial M)$ such that $\tilde f_\alpha(x_0), \tilde f_\beta(x_0) \in \mathbb C$, then $\tilde f_\gamma(x_0) \in \mathbb C$ for each $\gamma \in H_1(\partial M)$ and the surface $S$ can be chosen to be closed. 
\qed
\end{prop}

\begin{cor} \label{smallcharactervariety} 
Suppose that $M$ is a small knot exterior. \\
$(1)$ If $X_0$ is a non-trivial component of $X_{PSL_2}(M)$ and $x_0$ an ideal point of $X_0$, there is a primitive class $\alpha \in H_1(\partial M)$ such that $\Pi_{x_0}(\tilde f_{\alpha}) > 0$. Thus $\|\cdot \|_{X_0}$ is either a norm curve or a $\beta$-curve for some primitive $\beta \in H_1(\partial M)$. \\
$(2)$ If $\alpha \in H_1(\partial M)$ is a slope such that $X_{PSL_2}(M(\alpha))$ is infinite, then $\alpha$ is a boundary slope.
\end{cor}

\pf The first statement follows immediately from the previous proposition. For the second, assume that $X_{PSL_2}(M(\alpha))$ is infinite and choose a curve $X_0 \subset X_{PSL_2}(M(\alpha)) \subset X_{PSL_2}(M)$. Since $M$ is small, any essential surface in $M$ associated to an ideal point $x_0$ has boundary. Moreover, since $X_0 \subset X_{PSL_2}(M(\alpha))$, $\tilde f_\alpha(x_0) = 0$. Part (1) of the corollary shows that $\alpha$ is the unique slope with this property and therefore part (1) of the previous proposition shows that it is a boundary slope. 
\qed 

\begin{cor} \label{normcondition}
Let $X_0$ be a curve in $X_{PSL_2}(M)$ containing the characters of two discrete representations $\rho_1, \rho_2$ such that $\rho_j(\pi_1(\partial M))$ contains a non-trivial loxodromic element of $PSL_2(\mathbb C)$. If there are rationally independent classes $\alpha_1, \alpha_2 \in H_1(\partial M)$ such that $\rho_j(\alpha_j) = \pm I$ for $j = 1, 2$, then $\|\cdot\|_{X_0}$ is a norm. 
\end{cor}

\pf Our hypotheses imply that $\rho_j(\pi_1(\partial M)) \cong \mathbb Z \oplus \mathbb Z/c_j$ for some $c_j \geq 1$ and that any element of infinite order in this group is loxodromic. In particular this is the case for any element of $H_1(\partial M)$ which is rationally independent of $\alpha_j$ and therefore $f_{\alpha_1}(\chi_{\rho_2}) = |\mbox{trace}(\rho_2(\alpha_1))|^2 - 4 \ne 0$. Since $f_{\alpha_1}(\chi_{\rho_1}) = 0$, $f_{\alpha_1}|X_0$ is not constant. If there is a primitive class $\beta \in H_1(\partial M)$ such that $f_\beta|X_0$ is constant, then for some $j$, $\alpha_j$ and $\beta$ are rationally independent and so $\rho_j(\beta)$ is loxodromic. It follows that $f_\beta \equiv f_\beta(\chi_j) = (\lambda - \lambda^{-1})^2$ where $\lambda$ is not a root of unity. In particular $\tilde f_\beta$ takes on this value at each ideal point of $X_0$. Proposition \ref{idealvalue} now shows that $\tilde f_{\alpha_1}(x_0) \in \mathbb C$ for each ideal point $x_0$ of $X_0$. But this impossible as it would imply that $f_{\alpha_1}|X_0$ is constant. Thus $\|\cdot\|_{X_0}$ is a norm. 
\qed

\section{Dominations between small knot manifolds} \label{sec:knot exterior}

\subsection{Bounds on dominations between small knot manifolds } \label{seminorm} 

Let $M$ be a small knot manifold and denote by $T_1(M)$ the torsion subgroup of $H_1(M)$ and $F_1(M)  = H_1(M)/T_1(M) \cong \mathbb Z$ its free part. If $K_M = \hbox{kernel}(\pi_1(M) \to H_1(M) \to F_1(M))$, then
$$Z(\pi_1(M)) \cap K_M = \{1\} $$ 
where $Z(\pi_1(M))$ denotes the centre of $\pi_1(M)$. This is obvious when $M$ is hyperbolic. When it is Seifert fibred, it admits a Seifert structure with base orbifold $D^2(p,q)$ and so $Z(\pi_1(M))$ is generated by the class of the fibre. Since $\partial M \ne \emptyset$, $M$ admits a horizontal surface, and as $D^2(p,q)$ is orientable, the surface is non-separating. Thus the fibre class has infinite order in $H_1(M)$ which yields the desired conclusion. (See \cite{Wal1} for more details.) 

Given a component $X_0$ of $X_+^{irr}(M)$, set 
$$K_M(X_0) = \hbox{Ker}(X_0) \cap K_M.$$

Recall that a component $X_0 \subseteq X_+^{irr}(M)$ is called {\it principal} if it contains the character of an irreducible, discrete faithful representation of $\pi_1(M)/Z(\pi_1(M))$. The character varieties of small knot manifolds have dimension $1$ [CCGLS] and since they are either hyperbolic or Seifert fibred, they contain at least one principal component. 
Moreover, such a component is contained in $X_+^{str}(M)$ unless $M$ is a twisted $I$-bundle over the Klein bottle. 

For a small knot manifold define
$${\cal I}_M = \{\pi_1(M) / K_M(Y_0) : Y_0 \hbox{ an algebraic component of } X_{+}^{irr}(M)\} / \hbox{isomorphism}$$ 
and note 
$$|{\cal I}_M| \leq \# \hbox{ algebraic components of } X_{+}^{irr}(M).$$ 

\begin{thm} \label{epibound}
Let $M$ and $N$ be small knot manifolds and $X_0$ a principal component of $X_{PSL_2}(N)$. If $\varphi: \pi_1(M) \to \pi_1(N)$ is an epimorphism and $Y_0 = \varphi^*(X_0) \subseteq X_+^{irr}(M)$, then $\hbox{kernel}(\varphi) = K_M(Y_0)$. Thus $|{\cal I}_M|$ is an upper bound for the number of isomorphism classes of groups $\pi_1(N)$ where $N$ is a small knot manifold for which there is an epimorphism $\varphi: \pi_1(M) \to  \pi_1(N)$. 
\end{thm}

\pf By Lemmas \ref{ker} and \ref{kerprincipal} we have $\hbox{Ker}(Y_0) = \varphi^{-1}(\hbox{Ker}(X_0)) = \varphi^{-1}(Z(\pi_1(N)))$. Since $\varphi$ is surjective, it induces an isomorphism $F_1(M) \to F_1(N)$ and therefore $\varphi^{-1}(K_N) = K_M$. It follows that $K_N(Y_0) = \varphi^{-1}(Z(\pi_1(N))) \cap \varphi^{-1}(K_N) =  \varphi^{-1}(Z(\pi_1(N)) \cap K_N) =  \varphi^{-1}(1) =  \hbox{kernel}(\varphi)$. Hence $\pi_1(N) \cong \pi_1(M) / K_M(Y_0)$ represents an element of ${\cal I}_M$. This completes the proof. 
\qed

\begin{cor} \label{cor:canonical curve}
Suppose that $M$ is a small knot manifolds such that $X_+^{irr}(M)$ contains only principal components. Then any non-zero degree map $f: M \to N$ is homotopic to a cover. In particular, if $M$ covers no orientable manifold other than itself, it is minimal. 
\end{cor}

\pf Since $\hbox{degree}(f) \ne 0$, $f_\#(\pi_1(M))$ has finite index in $\pi_1(N)$ and so consideration of the cover $\tilde N \to N$ corresponding to $\hbox{image}(f_\#)$ and $\tilde f: M \to N$, the lift of $f$, we can suppose, without loss of generality, that $f_\#$ is surjective. Theorem \ref{epibound} implies that $\hbox{kernel}(f_\#) = K_M(Y_0)$ for a principal component $Y_0$ of $X_{PSL_2}(M)$. Thus by Lemma \ref{kerprincipal} $\hbox{kernel}(f_\#) = Z(\pi_1(M)) \cap K_M = \{1\}$ and therefore $f_\#$ is an isomorphism which preserves the peripheral structure. Hence $f$ is homotopic to a homeomorphism \cite{Wal2}. This completes the proof.  
\qed

\begin{exa} \label{twistpretzel} 
{\rm The $PSL_2(\mathbb C)$-character variety of the exterior of a non-trivial twist knot or a $(-2, 3, n)$ pretzel knot, $n \not \equiv 0$ (mod $3$), has a unique non-trivial component ((\cite{Bu}, \cite{Mat}). Further, any such manifold covers no orientable manifold but itself as otherwise the cover would be regular \cite{GW} and so the knot would admit a free symmetry. But by \cite{GLM}  and \cite{BolZ} (see also \cite{Ha}), this eventuality only occurs in the case of the trefoil knot exterior where the result is readily verified. Thus we have two infinite families of minimal manifolds. It is interesting to note that when $n \equiv 0$ (mod $3$), the character variety of the exterior of the $(-2, 3, n)$ pretzel knot has precisely two non-trivial components, one principal and the other corresponding to a strict domination of the trefoil knot exterior. (This follows from the analysis in the section ``$r$-curves" of  \cite{Mat}.) Thus $M$ is ${\cal H}$-minimal in this case.}
\end{exa}

We say that non-zero degree maps $f_j: M \to N_j$ ($j=1,2$) to be {\it equivalent} if there is a homeomorphism $g: N_1 \to N_2$ such that $f_2 \simeq g \circ f_1$. 

Set 
$${\cal N}_M = \{\pi_1(M) / K_M(Y_0) : Y_0  \hbox{ a norm curve in } X_+^{irr}(M) \} / \hbox{isomorphism} \subseteq {\cal I}_M.$$

\begin{thm} \label{dombound} 
Let $M$ be a small knot manifold. \\ 
$(1)$ The number of equivalence classes of $\pi_1$-surjective non-zero degree maps $M \to N$ is bounded above by $|{\cal I}_M|$. More precisely, \\ 
\indent $(a)$  The number of equivalence classes of $\pi_1$-surjective non-zero degree maps $M \to N$ \\ \indent \hspace{.55cm} where $N$ is hyperbolic is bounded above by $|{\cal N}_M|$. \\ 
\indent $(b)$ The number of equivalence classes of $\pi_1$-surjective non-zero degree maps $M \to N$  \\ \indent \hspace{.55cm} where $N$ is Seifert is bounded above by $|{\cal I}_M| - |{\cal N}_M|$. \\ 
$(2)(a)$ The number of equivalence classes of non-zero degree maps from $M$ to a hyperbolic  \\ \indent \hspace{.45cm} manifold is bounded above by a constant depending only on $X_{PSL_2}(M)$. \\
\indent $(b)$ The number of homeomorphism classes of Seifert fibred manifolds dominated by $M$  \\ \indent \hspace{.45cm} is bounded above by a constant depending only on $X_{PSL_2}(M)$.
\end{thm}

\begin{rem}
{\rm Small Seifert knot manifolds have base orbifolds of the form $D^2(p,q)$ where $p, q \geq 2$. They are also surface bundles over the circle with periodic monodromies. If $F$ is the fibre and $h: F \to F$ the monodromy, then $D^2(p, q) = F/h$. It is clear that such manifolds admit self-covering maps of arbitrarily high degree. On the other hand, if $M$ is a hyperbolic knot manifold, it is well-known that the degree of a proper map $f: M \to N$ is bounded above by $vol(M)$. These contrasting facts are the root of the difference in the statements of parts (2)(a) and (b) of the theorem. }
\end{rem}

\noindent {\bf Proof of Theorem \ref{dombound}.} A standard transversality argument shows that if $M \to N$ is a non-zero degree map, then $N$ is small.

(1) Let $f: M \to N$ be a $\pi_1$-surjective non-zero degree map. According to Theorem \ref{epibound}, if $X_0$ is a principal component of $X_+^{irr}(N)$ and $Y_0 = f_\#^*(X_0)$, then $\hbox{kernel}(f_\#) = K_M(Y_0)$. Suppose that $f': M \to N'$ is another $\pi_1$-surjective non-zero degree map and $X_0'$ is a principal component of $X_+^{irr}(N')$ such that $Y_0 = (f')_\#^*(X'_0)$ and therefore $\hbox{kernel}((f')_\#) = K_M(Y_0)$. We claim that $f$ and $f'$ are equivalent. To see this, observe that by construction, there is an isomorphism $\varphi: \pi_1(N) \to \pi_1(N')$ such that $f'_\# = \varphi \circ f_\#$. Since $N$ and $N'$ are Haken, it suffices to prove that $\varphi(\pi_1(\partial N)) \subseteq \pi_1(\partial N')$ \cite{Wal2}. This is clear when $N$, and therefore $N'$, is hyperbolic. Suppose then that they are Seifert fibred manifolds with base orbifolds ${\cal B}, {\cal B}'$. Since $\varphi(Z(\pi_1(N))) = Z(\pi_1(N'))$, $\varphi$ induces an isomorphism $\bar \varphi: \pi_1({\cal B}) \to \pi_1({\cal B}')$. By hypothesis, the finite index subgroup $f_\#(\pi_1(\partial M))$ of $\pi_1(\partial N)$ is sent into $\pi_1(\partial N')$ by $\varphi$ and therefore some positive power of a generator of $\pi_1(\partial {\cal B}) \cong \mathbb Z$ is sent by $\bar \varphi$ into $\pi_1(\partial {\cal B}')$. As $\pi_1({\cal B})$ is isomorphic to a free product of two finite cyclic groups ($\mathbb Z/p * \mathbb Z/q$ if ${\cal B} = D^2(p,q)$), the centralizer of any element is cyclic and so we conclude that $\bar \varphi(\pi_1(\partial {\cal B})) \subseteq \pi_1(\partial {\cal B}')$. This fact together with the identity $\varphi(Z(\pi_1(N))) = Z(\pi_1(N'))$ imply that $\varphi$ preserves peripheral subgroups and so is induced by a homeomorphism $g: N \to N'$. Thus the number of equivalence classes of $\pi_1$-surjective non-zero degree maps is bounded above by $|{\cal I}_M|$. 

The proof of (1)(a) follows from the fact that if $N$ is hyperbolic, then a principal curve $X_0 \subset X_{PSL_2}(N)$ is a norm curve (cf. \S \ref{culler-shalen}). Thus if $f: M \to N$ has non-zero degree, then $f_\#|\pi_1(\partial M)$ is injective and so $Y_0 = f_\#^*(X_0)$ is also a norm curve. Part (1)(b) follows similarly. Simply note that a principal curve $X_0$ for a small Seifert knot manifold is never a norm curve since the class $\gamma \in \pi_1(N)$ of a regular fibre is sent to $\pm I$ by each $\rho \in R_{X_0}$ (Proposition \ref{kerprincipal}(1)). 

(2)(a) Consider a non-zero degree map $f: M \to N$ where $N$ is hyperbolic. Since $M$ is small and $f_\#(\pi_1(M))$ has finite index in $\pi_1(N)$, $M$ must also be hyperbolic. 
Now $f$ factors $M \stackrel{\tilde f}{\longrightarrow} \tilde N \stackrel{g}{\longrightarrow} N$ where $\tilde f$ is $\pi_1$-surjective and $g$ is a cover of degree at most $vol(M)$, a constant determined by $X_{PSL_2}(M)$. By part (1), there are only finitely many possibilities for $M \stackrel{\tilde f}{\longrightarrow} \tilde N$ up to equivalence. Hence as $vol(\tilde N) \leq vol(M)$, we are reduced to proving the following claim.

\begin{claim} \label{hypbound}
Given a hyperbolic $3$-manifold $W$ with $vol(W) \leq vol(M)$, the number of equivalence classes of covers $p: W \to V$ is bounded above by a constant depending only on $vol(M)$. 
\end{claim}

\pf A $d$-fold covering $p: W \to V$ is determined up to equivalence by a homomorphism $\rho: \pi_1(V) \to {\cal S}_d$, where ${\cal S}_d$ is the symmetric group on a set of $d$ elements. Further, $\rho$ induces a transitive action of $\pi_1(V)$ on the set for which $\pi_1(W)$ is the stabilizer of an element. Thus $K_\rho := \hbox{kernel}(\rho) \subseteq \pi_1(W)$ and has index at most $(d-1)!$. It follows that $K_\rho$ is the fundamental group of a finite cover $\hat  W \to W$ whose volume is bounded above by $(d-1)! vol(M)$. Therefore, the isometry group $Out(K_\rho)$ of $\hat W$ is a finite group of cardinality at most $24 (d-1)! vol(M)$ since by \cite{Mey} the volume of a cusped hyperbolic $3$-orbifold is $\geq 1/6$. Since $d \leq vol(W) \leq vol(M)$, $|Out(K_\rho)|$ is bounded by a constant depending only on $vol(M)$.

Let $H_V = \rho(\pi_1(V)) \subseteq {\cal S}_d$ and consider the exact sequence 
$$1 \to K_\rho \to \pi_1(V) \stackrel{\rho}{\longrightarrow} H_V \to 1 $$ 
This extension is determined by the associated homomorphism $H_V \to \hbox{Out}(K_\rho)$ since $K_\rho$ has trivial centre (see Chap. IV.6 of \cite{Brn}). The number of such homomorphisms is bounded above by $|Out(K_\rho)|^{d!}$, which in turn is bounded by a constant depending only on $vol(M)$. Since the number of possible groups $H_V$ is bounded above by $2^{d!}$, it remains to show that the number of possibilities for $K_\rho$ is bounded above by a constant depending only on $vol(M)$. But there is a universal constant C such that the rank of $\pi_1(W)$ is no more than $C vol(W) \leq C vol(M)$ \cite{Ad}. On the other hand, the number of normal subgroups of $\pi_1(W)$ of index at most $(d-1)!$ is bounded above by $((d-1)!)!)^{rank(\pi_1(W))}$. This completes the proof of the Claim and therefore of (2)(a). 
\qed (Claim \ref{hypbound})

(2)(b) Part (1)(b) shows that the number of $\pi_1$-surjective dominations of $M$ to a Seifert manifold is bounded above by $|{\cal I}_M| - |{\cal N}_M|$. Given such a domination $M \to \tilde N$, let $D^2(p,q)$ be the base orbifold of $\tilde N$ and note that  $p, q$ are determined by the associated curve in $X_+^{irr}(M)$. To complete the proof, it suffices to show that the number of homeomorphism types of knot manifolds $N$ finitely covered by $\tilde N$ is bounded above by a constant depending only on $p, q$. 

Let $\tilde N \to N$ be a cover and suppose $D^2(a,b)$ is the base orbifold of $N$. There is an induced orbifold cover $D^2(p,q) \to D^2(a,b)$ of some degree $d \geq 1$. An elementary calculation based on Euler characteristics shows that if $d > 1$, then up to permutation of $a, b$ either \\ 
\indent (a) $a = d = 2$ and $p = q = b$, or \\
\indent (b) $a = b = p = q = 2$. \\
A small Seifert knot manifold $N$ with base orbifold $D^2(a,b)$ is the union of two vertical solid tori along a vertical annulus where the solid tori are of fibred type $(a,r)$ and $(b,s)$ where $1 = \gcd(a,r) = \gcd(b,s)$. The homeomorphism type of $N$ is unchanged if we alter the gluing map by a homeomorphism which extends over either solid torus. Hence 
the number of homeomorphism types of small Seifert knot manifolds with base $D^2(a,b)$ is at most $\frac{ab}{4} \leq \frac{pq}{4}$. Thus we are done.
\qed 

\subsection{Rigidity in $\pi_1(M)$ and bounds on sequences of dominations} \label{rigidity} 

We assume that $M$ is a small knot manifold in this section. 

Call $\gamma \in \pi_1(M)$ {\it rigid} if $f_\gamma|X_0$ is constant for some principal curve $X_0$ of $X_{PSL_2}(M)$. Equivalently, $d_{X_0}(\gamma) = 0$ (cf. \S \ref{culler-shalen}). (This condition is independent of the choice of principal curve.) For instance, if a positive power of $\gamma \in \pi_1(M)$ is central, then $\gamma$ is rigid. We call $\gamma$ {\it non-rigid} otherwise. Finally we call $\gamma \in \pi_1(M)$ {\it totally non-rigid} if $f_\gamma|X_0$ is non-constant for all curves $X_0 \subseteq X_{+}^{irr}(M)$. 

\begin{lemma} \label{non-rigid} 
Let $M$ be a small knot manifold.  \\
$(1)$ If $M$ is hyperbolic, every non-trivial element of $\pi_1(\partial M)$ is non-rigid. \\
$(2)$ If $M$ is Seifert, an element of $\pi_1(M)$ is rigid if and only if some non-zero power of it is central. \\ 
$(3)$ If $\alpha \in H_1(\partial M) = \pi_1(\partial M) \subset \pi_1(M)$ is a slope which is not a boundary slope, then $\alpha$ is totally non-rigid. 
\end{lemma}

\pf  Part (1) is proved in Proposition 1.1.1 of \cite{CGLS}. 

(2) Suppose that $M$ is Seifert. Since it is small, its base orbifold is of the form ${\cal B} = D^2(p, q)$ for some integers $p, q \geq 2$. If no positive power of $\gamma \in \pi_1(M)$ is central, then $\gamma$ projects to an element $\bar \gamma \in \pi_1({\cal B}) \cong \mathbb Z/p * \mathbb Z/q$ of reduced length at least $2$ with respect to any generators $x$ of $\mathbb Z/p$ and $y$ of $\mathbb Z / q$. It follows as in the proof of Theorem 1 of \cite{BMS} that $f_{\bar \gamma}$ is non-constant on each non-trivial curve of $X_{PSL_2}(\mathbb Z/p * \mathbb Z/q)$. Thus $\gamma$ is non-rigid. 

(3) Let $X_0$ be a non-trivial curve in $X_{PSL_2}(M)$ and suppose that $f_\alpha|X_0$ is constant. Then for any ideal point $x$ of $X_0$, $f_\alpha (x) \in \mathbb C$. But this impossible as otherwise Proposition \ref{idealvalue} implies that either $M$ is large or $\alpha$ is a boundary slope. Thus $f_\alpha|X_0$ is not constant.
\qed  

\begin{lemma} \label{samedegree}
Suppose that $M$ and $N$ are small knot manifolds and $\varphi: \pi_1(M) \to
\pi_1(N)$ is a virtual epimorphism.  \\
$(1)$ $\varphi^*$ induces a birational isomorphism between $X_+^{str}(N)$ and a
union of algebraic components of $X_+^{str}(M)$. In particular if $Y_0$ is a component of $X_{+}^{str}(N)$, then $X_0 = \varphi^*(Y_0)$ is a component of $X_{+}^{str}(M)$ and for each $\gamma \in \pi_1(M)$ we have
$$d_{X_0}(\gamma) = d_{Y_0}(\varphi(\gamma)).$$  
$(2)$ If there is a principal component $X_0$ of $X_{PSL_2}(M)$ contained in $\varphi^*(X_{+}^{str}(N))$, then $\varphi$ is injective. 
\end{lemma}

\pf Parts (1) follows from the remark in the opening paragraph of this section and Corollary \ref{1-1}. We consider part (2) then.

Suppose that $\varphi^*(X_{+}^{str}(Y_0)) = X_0$ for some component $Y_0$ of $X_{+}^{str}(N)$ and principal component $X_0$ of $X_{PSL_2}(M)$. Lemma \ref{ker} (2) shows that $\hbox{kernel}(\varphi) \subseteq \hbox{Ker}(X_0) \subseteq Z(\pi_1(M))$. Thus if $\hbox{kernel}(\varphi) \ne \{1\}$, $M$ is Seifert fibred, and as $\pi_1(M)$ and $\pi_1(N)$ are torsion free, $\hbox{kernel}(\varphi) = Z(\pi_1(M)) \cong \mathbb Z$. But this is impossible as it would imply that $\pi_1(N)$ contains a subgroup isomorphic to $\pi_1(M)/Z(\pi_1(M))$, which is the free product of two finite cyclic groups. Thus $\hbox{kernel}(\varphi) = \{1\}$. 
\qed

We define the {\it strict degree} of an element $\gamma$ of the fundamental group of a small knot manifold $M$ to be the sum 
$$d_M (\gamma)= \sum_{\stackrel{components \; X_0 \; of}{X_+^{str}(M)}} d_{X_0}(\gamma).$$
Note that $d_M (\gamma) > 0$ if $\gamma$ is non-rigid as long as $M$ is not a twisted $I$-bundle over the Klein bottle. The following lemma is of use in this case.

\begin{lemma} \label{kleinbottle} 
Let $M, N$ be small manifolds and $\varphi: \pi_1(M) \to \pi_1(N)$ a virtual epimorphism. \\
$(1)$ If $N$ is a twisted $I$-bundle over the Klein bottle and $X_0$ the principal curve in $X_+^{irr}(N)$, then $\varphi^*(X_0)$ is a non-trivial curve in $X_{PSL_2}(M)$. \\
$(2)$ If $M$ is a twisted $I$-bundle over the Klein bottle then so is $N$ and $\varphi$ is injective.
\end{lemma}

\pf (1) Let $\tilde N \stackrel{g}{\longrightarrow} N$ be the cover corresponding to $\varphi(\pi_1(M))$. A finite cover of $N$ is either homeomorphic to $N$ or $S^1 \times S^1 \times I$ and so as $M$ has first Betti number $1$, $\tilde N$ is also a twisted $I$-bundle over the Klein bottle. The reader will then verify that $Y_0 = (g_\#)^*(X_0)$ is a principal curve for $\tilde N_i$. But if $\tilde \varphi: \pi_1(M) \to \pi_1(\tilde N)$ is the surjection induced by $\varphi$,  $\tilde \varphi^*(Y_0) = \varphi^*(X_0)$ is a non-trivial curve in $X_{PSL_2}(M)$. 

(2) If $M$ is a twisted $I$-bundle over the Klein bottle, then $\pi_1(N)$ has a finite index abelian subgroup and therefore is also a twisted $I$-bundle over the Klein bottle. As in the proof of (1), if $\tilde N \to N$ is the cover corresponding to $\varphi(\pi_1(M))$, then $\tilde N$ is also a twisted $I$-bundle over the Klein bottle. But the fundamental group of such a manifold is Hopfian so the induced epimorphism $\pi_1(M) \to \pi_1(\tilde N)$ is an isomorphism. Thus $\varphi$ is injective. 
\qed

\begin{thm}\label{degreeinject}
Let $\varphi: \pi_1(M) \to \pi_1(N)$ be a virtual epimorphism. \\
$(1)$ $d_N(\varphi(\gamma)) \leq d_M(\gamma)$ for all $\gamma \in \pi_1(M)$. \\
$(2)$ If $\gamma \in \pi_1(M)$ is not rigid and $d_N(\varphi(\gamma)) = d_M(\gamma)$, then $\varphi$ is injective.  
\end{thm}

\pf The first assertion is a consequence of part (1) of Lemma \ref{samedegree}. To prove the second, note that Lemma \ref{samedegree}(2) shows that we can suppose there is a principal component $X_0$ of $X_{PSL_2}(M)$ which is not contained in the Zariski closure of the image of $\varphi^*$. Lemma \ref{kleinbottle}(2) shows that we can also suppose that $M$ is not a twisted $I$-bundle over the Klein bottle. Thus as $\gamma$ is not rigid, $d_{X_0}(\gamma) > 0$ and therefore $d_M (\gamma) \geq d_N(\varphi(\gamma)) + d_{X_0}(\gamma)> d_N(\varphi(\gamma)) = d_M (\gamma)$, which is impossible. 
\qed 

\begin{rem}
{\rm Note that under the hypotheses of part (2) of the theorem, work of Waldhausen \cite{Wal2} implies that $\varphi$ is induced by a covering map $M \to N$ as long as it  preserves the peripheral subgroups of $\pi_1(M)$ and $\pi_1(N)$. This is automatically satisfied if $N$ is hyperbolic. }
\end{rem}

Our next result gives an a priori bound on the length of certain sequences of homomorphisms between the fundamental groups of small knot manifolds. 

\begin{thm} \label{homomorphismsequence}
Let $M$ be a small knot manifold and consider a sequence of homomorphisms    
$$\pi_1(M) \stackrel{\varphi_1}{\longrightarrow} \pi_1(N_1) \stackrel{\varphi_2}{\longrightarrow}  \cdots \stackrel{\varphi_n}{\longrightarrow}  \pi_1(N_n)$$ 
none of which is injective. If $N_i$ is small and $\varphi_i$ is a virtual epimorphism for each $i$, then $n \leq d_M(\gamma)$ for each totally non-rigid element $\gamma \in \pi_1(M)$. Moreover, if $n = d_M(\gamma)$ for some such $\gamma$, then $N_n$ is a twisted $I$-bundle over the Klein bottle. 
\end{thm}

\pf Set $\psi_i = \varphi_i \circ \cdots \circ \varphi_1$ and let $\gamma \in \pi_1(M)$ be totally non-rigid. If $\psi_i(\gamma)$ is rigid for some $1 \leq i \leq n$ and $X_0 \subset X_{+}^{irr}(N_i)$ is a principal curve, then $f_{\psi_i(\gamma)}|X_0$ is constant, or equivalently, $f_\gamma|\psi_i^*(X_0)$ is constant (Identity (2.1.1)). Since $\gamma$ is totally non-rigid, $\psi_i^*(X_0)$ cannot be a non-trivial curve and therefore $X_0 \not \subset X_+^{str}(N_i)$ (Corollary \ref{1-1}). Thus $N_i$ is a twisted $I$-bundle over the Klein bottle. But Lemma \ref{kleinbottle}(1) shows that this case does not arise under our assumptions. It follows that $\psi_i(\gamma)$ is non-rigid for $1 \leq i \leq n$. Moreover, Lemma \ref{kleinbottle}(2) shows that if $N_i$ is a twisted $I$-bundle over the Klein bottle for some $i$, then $i = n$. Theorem \ref{degreeinject} now implies that 
$$d_{M}(\gamma) >  d_{N_{1}}(\psi_1(\gamma)) > \cdots > d_{N_{i}}(\psi_i(\gamma)) \cdots > d_{N_{n}}(\psi_{n}(\gamma)) \geq 0$$
with $d_{N_{n}}(\psi_{n}(\gamma)) = 0$ if and only if $N_n$ is a twisted $I$-bundle over the Klein bottle. This completes the proof. 
\qed 

\subsection{Dominations by two-bridge knot exteriors} \label{twobridgedomination}

Consider relatively prime integers $p,q$ where $p \geq 1$ is odd and let $k_{p/q}$
denote the two-bridge knot corresponding to the rational number $p/q$. Thus the $2$-fold cover of $S^3$ branched over $k_{p/q}$ is the lens space $L(p,q)$. It is a theorem of Schubert [Sch] that $k_{p/q}$ is equivalent to $k_{p'/q'}$ if and only if $L(p,q)$ is homeomorphic to $L(p',q')$. The exterior $M_{p/q}$ of $k_{p/q}$ is known to be small [HT]. Moreover it is hyperbolic unless $q \equiv \pm 1$ (mod $p$), in which case it is a $(p,2)$ torus knot.

The proof of the following unpublished result of Tanguay is contained in Appendix A.

\begin{prop} \label{meridiandegree} {{\rm (\cite{Tan})}}
Let $M$ be the exterior of the two-bridge knot of type $p/q$. If $\mu \in \pi_1(M)$ is a meridinal class, then $d_M(\mu)= \frac{p-1}{2}.$
\end{prop}

As a consequence we deduce: 

\begin{thm}
Consider a sequence of homomorphisms    
$$\pi_1(M_{p/q}) \stackrel{\varphi_1}{\longrightarrow} \pi_1(N_1) \stackrel{\varphi_2}{\longrightarrow}  \cdots \stackrel{\varphi_n}{\longrightarrow}  \pi_1(N_n)$$ 
none of which is injective. If $N_i$ is small and $\varphi_i$ is a virtual epimorphism for each $i$, then $n < \frac{p-1}{2}$. 
\end{thm}

\pf The meridinal slope $\mu$ of a two-bridge knot is not a boundary slope \cite{HT} so  Lemma \ref{non-rigid}(3) shows that it is totally non-rigid in $\pi_1(M_{p/q})$. Theorem \ref{homomorphismsequence} then yields the inequality $n \leq d_{M_{p/q}}(\mu) = \frac{p-1}{2}$ with equality only if $N_n$ is a twisted $I$-bundle over the Klein bottle. We saw in the proof of Lemma \ref{kleinbottle}(1) that if $\tilde N_n \to N_n$ is the cover corresponding to the image of $\varphi_n \circ \varphi_{n-1} \circ \ldots \circ \varphi_1$, then $\tilde N_n$ is also 
a twisted $I$-bundle over the Klein bottle. But this is impossible since $H_1(M_{p/q})$ is cyclic while $H_1(\tilde N_n)$ is not. Thus $n < \frac{p-1}{2}$.  
\qed 

This result can be significantly strengthened if the homomorphisms are induced by non-zero degree maps. This is the goal of the remainder of this section. 

\begin{thm} \label{2bridgemini} 
Let $N$ be a knot manifold and $\varphi: \pi_1(M_{p/q}) \to \pi_1(N)$ a homomorphism such that the image $\varphi(\mu)$ of a meridian $\mu$ is peripheral. \\
$(1)$ If $\varphi$ is an epimorphism, then $N$ is homeomorphic to the exterior $M_{p'/q'}$ of a 2-bridge knot in $S^3$. Moreover either $M_{p'/q'} = M_{p/q}$ or $p = kp'$ with $k > 1$. \\
$(2)$ If $\varphi(\pi_1(M))$ is of finite index $d$ in $\pi_1(N)$, then either $d = 1$ and the conclusions of part $(1)$ hold, or $N$ is Seifert fibred and $\varphi$ factors through an epimorphism $\tilde \varphi: \pi_1(M_{p/q}) \to \pi_1(M_{p'})$ to which the conclusions of part $(1)$ apply. Further, $\gcd(2p',d) = 1$. 
\end{thm}

\pf $(1)$  Since $\mu$ normally generates $\pi_1(M_{p/q})$, $ \varphi(\mu)$ does the same for $\pi_1(N)$. In particular $\varphi(\mu) \ne 1$ so if $\mu'$ is the slope on $\partial  N$ corresponding to the projective class of $\varphi_*(\mu) \in H_1(\partial N)$, then the manifold $ W = N( \mu')$ obtained by Dehn filling $\partial  N$ along the slope $\mu'$  is a homotopy $3$-sphere. 

Let $k'$ be the core of the surgery in $W =  N(\mu')$ and let $\widehat{W}_2(k')$ be the $2$-fold cover of
$W$ branched over $k'$. There is an induced surjective homomorphism
$\mathbb Z/p \cong \pi_1(L(p,q)) \to \pi_1(\widehat{W}_2(k'))$ and so the latter
group is finite cyclic $\mathbb Z/p'$ with $p'$ dividing $p$. Since $\pi_1(M_{p/q})$ is generated by two elements, the same holds for $\pi_1(N)$, hence $k'$ is a  $2$-generator knot in the homotopy sphere $W$. It follows as in  \cite{Wed} that $k'$ is prime
and thus $N$ cannot contain an essential annulus with slope $\mu'$. Thus the
$2$-fold branched covering $\widehat{W}_2(k')$ of $k'$ is irreducible and by
the orbifold theorem ([BP], [BLP], [CHK])
$\widehat{W}_2(k')$ it is itself a lens space and the covering involution conjugates to an orthogonal involution. Therefore $W = N(\varphi(\mu)) \cong S^3$ and $k'$ is a
two-bridge knot. In other words, $(N(\varphi(\mu)), k') \cong (S^3, k_{p'/q'})$ for some integers $p' \geq 1, q'$ with $p'$ dividing $p$ and $q'$ coprime with $p'$. Property P for two-bridge knots \cite{Tak} implies that $\varphi(\mu) = \mu'$ is a meridian of $k_{p'/q'}$. 

According to Theorem \ref{degreeinject} and Proposition \ref{meridiandegree}, either $\varphi$ is an isomorphism or $\frac{p-1}{2} = d_{M_{p/q}}(\mu) > d_{M_{p'/q'}}(\mu') =
\frac{p'-1}{2}$. In the first case  $k_{p/q} = k_{p'/q'}$ while in the second  
$p = kp'$ with $k > 1$. This is the conclusion of (1). 

$(2)$ Let $\tilde N \to N$ be the cover corresponding to the image of $\varphi$
and define $\tilde \varphi: \pi_1(M) \to \pi_1(\tilde N)$ and $\psi: \pi_1(\tilde N) \to
\pi_1(N)$ in the obvious way.  Part $(1)$ implies that $\tilde N$ is homeomorphic to the exterior $M_{p'/q'}$ of a 2-bridge knot in $S^3$. Hence by \cite{GW}, the cover $\tilde N \to N$ is regular and cyclic. If $N$ is hyperbolic, $\tilde N = N$ since hyperbolic 2-bridge knot exteriors admit no free symmetries by \cite{GLM} (see also \cite{Ha}), and therefore we are in case (1). Otherwise $N$ is Seifert and so $\tilde N$ is the exterior $M_p'$ of $k_{p'}$, which is the $(p',2)$ torus knot. Thus the Seifert structure on $M_{p'}$ has base orbifold $D^2(2, p')$. Since $p'$ is odd, the proof of part (2)(b) of Theorem \ref{dombound} shows that the cover $\tilde N \to N$ induces a homeomorphism of the underlying orbifolds. Thus the cover is a degree $d$ unwinding of a regular fibre of $N$ and so if $F$ is the fibre and $h: F \to F$ the monodromy of the realization of $N$ as a surface bundle over the circle, then $h^d$ is the monodromy of the realization of $\tilde N$ as a surface bundle. The induced homeomorphism on the level of orbifolds is $F/h \to F/h^d$ and so $d$ must be coprime with the order of $h$. But since $F/h \cong D^2(2,p')$, this order is a multiple of $2p'$. 
\qed 

The following two results are immediate consequences of the previous theorem: 

\begin{thm} \label{2bridgebound}
Consider a sequence of non-zero degree maps     
$$M_{p_0/q_0} = N_0 \stackrel{f_1}{\longrightarrow} N_1 \stackrel{f_2}{\longrightarrow}  \cdots \stackrel{f_n}{\longrightarrow}  N_n$$ 
between knot manifolds, none of which is homotopic to a homeomorphism. If $N_{n-1}$ is hyperbolic, there are coprime pairs $p_j, q_j$ $(1 \leq j \leq n)$ such that $N_j = M_{p_j / q_j}$ $(1 \leq j \leq n-1)$, $N_n$ is finitely covered by some $M_{p_n / q_n}$, and $p_{j-1} = k_j p_j$ for some integer $k_j > 1$ $(1 \leq j \leq n)$.  Hence, $n + 1$ is bounded above by the number of distinct multiplicative factors of $p$. 
\qed 
\end{thm}

\begin{cor} \label{2bridgemini2}
If $p$ is an odd prime, then $M_{p/q}$ is minimal if and only if it is hyperbolic  $($i.e. $q \not \equiv \pm 1 \hbox{ $($mod p$)$}$$)$.  
\qed
\end{cor}

\noindent If we consider domination via degree-one maps (i.e. 1-domination), we obtain stronger results:

\begin{cor} 
Let $N$ be a knot manifold and $f: M_{p/q} \to N$ a degree-one map. Then either $f$ is homotopic to a homeomorphism or $N$ is a two-bridge knot exterior $M_{p'/q'}$ where $p = kp', k > 1,$ and $\gcd(k, p') = 1$. 
\end{cor}

\pf Since a degree-one map induces an epimorphism on the level of fundamemtal groups, case $(1)$ of the previous theorem shows that if $f$ is not homotopic to a homeomorphism, then $N = M_{p'/q'}$ where $p = kp', k > 1$. Moreover, $f$ induces a degree-one map $L(p,q) \to L(p',q')$ between the $2$-fold branched covers. By Corollary 6 of [RoW], there is an integer c such that $q' \equiv (\frac{p}{p'}) c^2 q$ (mod $p'$). In particular this implies that $\gcd(\frac{p}{p'}, p') = 1$.
\qed

\begin{thm} \label{2bridgedegree1bound} 
Consider a sequence of degree-one maps     
$$M_{p/q} \stackrel{f_1}{\longrightarrow} N_1 \stackrel{f_2}{\longrightarrow}  \cdots \stackrel{f_n}{\longrightarrow}  N_n$$ 
between knot manifolds, none of which is homotopic to a homeomorphism. Then $n + 1$ is bounded above by the number of distinct prime factors of $p$. 
\qed 
\end{thm}

\begin{cor} \label{1minimal} 
If $p$ is a prime power, the two-bridge knot exterior $M_{p/q}$ does not stictly $1$-dominate any knot manifold. 
\qed 
\end{cor}

\section{Sets of discrete $PSL_2(\mathbb C)$-characters}\label{sec:characters}

We investigate sequences of $PSL_2(\mathbb C)$ characters of representations of the fundamental groups of small knot manifolds whose images are discrete. This leads us in particular to proofs of Theorems \ref{thm:characters} and \ref{thm:domination}. Our analysis relies fundamentally on the convergence theory of Kleinian groups and hyperbolic $3$-manifolds. 

\subsection{Convergence of Kleinian groups and hyperbolic $3$-manifolds}  \label{convkleinhyp}

A metric space is {\it proper} if all of its closed and bounded subsets are compact. A sequence of proper pointed metric spaces $(X_n, x_n)$ is said to converge \emph{geometrically} to a metric space $(X_0, x_0)$ if for every $r > 0$, the sequence of compact metric balls $\{B_{X_n}(x_n; r)\}$ converges in the Gromov
bilipschitz topology to $B_{X_0}(x_0; r)$. (See chapter 3 of \cite{Gro} and also chapter E \cite{BeP}, chapter 7 \cite{MT}.) 

We recall the thick/thin decomposition of a complete, finite volume hyperbolic $3$-manifold $V$ (chapter D of \cite{BeP}): given a 
positive constant $0< \mu \leq \mu_0$, where $\mu_0$ is the Margulis constant, $V$ decomposes as $V_{[\mu, \infty)} \cup V_{(0,\mu]}$ such that:  

\indent $V_{[\mu, \infty)} =  \{ x \in V : \hbox{inj}(x) \geq \mu \}$  is the $\mu$-thick part of $V$ 

\indent  $V_{(0,\mu]} = \{ x \in V : \hbox{inj}(x) \leq \mu \}$ is the $\mu$-thin part of $V$.

\noindent For $\mu \leq \mu_0,$ each component of $\mu$-thin part of $V$ is either empty, or a geodesic neighborhood of a closed geodesic (a Margulis tube, homeomorphic to $S^1 \times D^2$) or a cusp with torus cross sections (homeomorphic to $T^2 \times [0, \infty)$).

Let $\{V_n\}$  be a sequence of pointed, closed, connected, orientable, hyperbolic $3$-manifolds whose volumes are bounded above. There is a sequence of base points $x_n \in (V_{n})_{[\mu_0, \infty)}$  such that some subsequence $\{(V_j, x_j)\}$ converges to a pointed, complete, finite volume, hyperbolic $3$-manifold $(V, x)$. In particular this implies that given $\varepsilon > 0$ and $0< \mu \leq \mu_0$, for $j \geq n_{0}(\varepsilon,\mu)$ the $\mu$-thick parts of $V_j$ and $V$ are $(1+ \varepsilon)$-bilipschitz homeomorphic. Moreover 
$$vol(V) = \lim_j vol(V_j),$$ 
(see chapter E of \cite{BeP}, theorem 7.9 of \cite{MT}). Further, if $V$ is closed, then $V = V_{j}$ for $j \gg 0$ and if $V$ is not closed, $V_j$ is obtained by Dehn filling $V$ for $j \gg 0$ (see chapter 5 of \cite{Thu}, chapter E of \cite{BeP}). By a Dehn filling of a complete, non-compact, finite volume hyperbolic $3$-manifold $V$ we mean a Dehn filling of some compact core $V_0$ of $V$.

In order to simplify the presentation, base points for fundamental groups and pointed  metric spaces will often be supressed from the notation. In particular we will say that a sequence $\{V_n\}$ of hyperbolic manifolds converges geometrically to a hyperbolic manifold $V$ if it does so under a suitable choice of base points.

We come now to the algebraic counterpart of this notion of geometric convergence. A good source on this topic is the paper \cite{JM} of J\o rgensen and Marden. The torsion-free case is dealt with in chapter 7 of \cite{MT}.

The {\it envelope} of a sequence $\{\Gamma_n\}$ of subgroups of $PSL_2(\mathbb C)$  is defined as 
$$Env(\{\Gamma_n\}) := \{\gamma = \lim_n \gamma_n : \gamma_n \in \Gamma_n \mbox{ for all } n\} \subset PSL_2(\mathbb C).$$
Clearly, $Env(\{\Gamma_n\})$ is a subgroup of $PSL_2(\mathbb C)$. 

A {\it Kleinian} group is a discrete subgroup of $PSL_2(\mathbb C)$. A {\it Fuchsian} group is a discrete subgroup of $PSL_2(\mathbb R)$. 

\begin{prop} \label{elemordiscr} {\rm (cf. Lemmas 3.2 and 3.6, \cite{JM} ) }$\;$ \\
$(1)$ If each $\Gamma_n$ is a non-elementary Kleinian group, then $Env(\{\Gamma_n\})$ is either elementary or discrete. \\
$(2)$ If $Env(\{\Gamma_n\})$ is non-elementary, then each $\gamma \in \overline{\cup (\Gamma_n \setminus \{\pm I\})}$ is either loxodromic, or parabolic, or elliptic of finite order $m \geq 2$. In the latter case, for any subsequence $\{\gamma_{n'}\}$ which converges to $\gamma$, $\gamma_{n'}$ has order $m$ for $n' \gg 0$.  
\qed 
\end{prop}

\begin{cor} \label{torsionfree} 
If each $\Gamma_n$ is a torsion-free non-elementary Kleinian group, then $Env(\{\Gamma_n\})$ is either abelian or discrete, and if it is non-abelian, each $\gamma \in \overline{\cup (\Gamma_n \setminus \{1\})}$ is either loxodromic or parabolic.  
\qed   
\end{cor}

A sequence $\{\Gamma_n\}$ of subgroups of $PSL_2(\mathbb C)$ is said to {\it converge geometrically} to a subgroup $\Gamma_0$ of $PSL_2(\mathbb C)$
if $\Gamma_0 = Env(\{\Gamma_{j}\})$ for every subsequence $\{ j \}$ of $\{ n \}$. The sequence $\{\Gamma_n\}$ is said to {\it converge algebraically} to
$\Gamma_0$ if there is a finitely generated group $\pi$ and representations $\rho_n \in R_{PSL_2}(\pi)$ ($n \geq 0$) such that $\Gamma_n = \rho_n(\pi)$ and $\lim_n \rho_n = \rho_0$. Note that if $\{\Gamma_n\}$ converges  algebraically to $\Gamma_0$ and geometrically to $\Gamma$, then
$\Gamma_0 \subseteq \Gamma \subseteq Env(\{\Gamma_{n}\})$. 

We record the following result for later use. Proofs in the torsion-free case can be found in Theorems 7.6, 7.7, 7.12, 7.13, and 7.14 of \cite{MT}. The general case can be dealt with using the results of \cite{JM}. 

\begin{prop} \label{convergence}  
Suppose that $\pi$ is a finitely generated group and $\rho_n: \pi \to PSL_2(\mathbb C)$ is a sequence which converges to $\rho_0 \in R_{PSL_2}(\pi)$. For $n \geq 0$ set $\Gamma_n = \rho_n(\pi)$ and suppose that for $n \geq 1$, $\Gamma_n$ is a non-elementary Kleinian group. Then \\
$(1)$ $\Gamma_0$ is a non-elementary Kleinian group. \\
$(2)$ for $n \gg 0$ there is a homomorphism $\theta_n: \Gamma_0 \to
\Gamma_{n}$ such that
$\rho_{n} = \theta_n \circ \rho$. Further, \\ \indent $\lim \theta_n =
1_{\Gamma_0}$. \\
$(3)$ there are a non-elementary Kleinian group $\Gamma$ containing
$\Gamma_0$ and a subsequence $\{j\}$ of  \\ \indent $\{n\}$ such that $\{\Gamma_{j}\}$ converges geometrically to $\Gamma$. Moreover, the homomorphisms $\theta_{j}$ of \\ \indent part $(2)$ extend to homomorphisms $\Gamma \to \Gamma_{j}$, which we continue to denote $\theta_{j}$, in such a  \\ \indent way that $\lim_j \theta_j =
1_{\Gamma}$. \\
$(4)$ the quotient spaces $\mathbb H^3/\Gamma_{j}$ converge geometrically to  $\mathbb
H^3/\Gamma$ 
\qed
\end{prop}

In the remainder of the paper we investigate sets of discrete characters and apply our results to study sequences of non-zero degree maps between closed manifolds. Let $M$ be a knot manifold and $X_0$ a subvariety of $X_{PSL_2}(M)$. Set 
$$D(X_0) = \{ \chi_\rho \in X_0 : \rho \hbox{ is discrete and non-elementary}\}$$
$$D^*(X_0) = \{ \chi_\rho \in D(X_0) : \rho \hbox{ is torsion free}\}$$ 
$$D_{0}^*(X_0) = \{ \chi_\rho \in D^*(X_0) : \rho \hbox{ has non-zero volume}\}.$$
Note that the image of any $\rho \in R_{PSL_2}(M)$ whose character is contained in $D_0^*(X_0)$ is the fundamental group of a complete hyperbolic $3$-manifold. 

\begin{lemma} \label{standardimage} 
Let $\chi_\rho \in D^*_0(X_0)$, $\Gamma = \rho(\pi_1(M))$, and $V = \mathbb H^3 / \Gamma$.  \\
$(1)$ If $\rho|\pi_1(\partial M)$ is injective, then a compact core $V_0$ of $V$ is a hyperbolic knot manifold and there is a proper non-zero degree map $f: M \to V_0$ such  that $f_\#: \pi_1(M) \to \pi_1(V_0) = \Gamma$ is conjugate to $\rho$.  \\ 
$(2)$ If $\rho|\pi_1(\partial M)$ is not injective, $V$ is closed and there is a slope $\alpha$ on $\partial M$ and a non-zero degree map $f: M(\alpha) \to V$ such  that the composition $\pi_1(M) \to \pi_1(M(\alpha)) \stackrel{f_\#}{\longrightarrow} \pi_1(V) = \Gamma$ is conjugate to $\rho$.  \\ 
$(3)$ If $v_0 > 0$ is the minimal volume for complete, connected, orientable, hyperbolic $3$-manifolds, then $|vol(\chi_\rho)| \geq v_0$. 
\end{lemma} 

\pf If $\rho|\pi_1(\partial M)$ is injective, there is a compact core $V_0$ of $V$ and a torus $T$ in $\partial V_0$ such that $\rho(\pi_1(\partial M)) \subset \pi_1(T)$. Thus there is a proper map $f: (M, \partial M) \to (V_0, T)$ realizing $\rho$. By the definition of the volume of a representation (\cite{Dun}), $|\hbox{degree}(f)| vol(V) = |vol(\rho)|  \ne 0$. In particular $|\hbox{degree}(f)| \ne 0$, which implies that $\partial V_0 = T$ and $|vol(\rho)| \geq vol(V) \geq v_0$. On the other hand, if $\hbox{kernel}(\rho| \pi_1(\partial M)) \ne \{\pm I\}$, $\rho$ factors $\pi_1(M) \to \pi_1(M(\alpha)) \stackrel{\bar \rho}{\longrightarrow} \pi_1(V) = \Gamma$ for some slope $\alpha$ since the image of $\rho$ is torsion free. There is a map $f: M(\alpha) \to V$ associated to the homomorphism $\bar \rho$ and Lemma 2.5.4 of \cite{Dun} implies that $vol(\bar \rho) = vol(\rho)$. Then $|\hbox{degree}(f)| vol(V) = |vol(\bar \rho)| = |vol(\rho)| \ne 0$ and again we see that $|\hbox{degree}(f)| \ne 0$ so that $V$ must be closed and $|vol(\rho)| \geq vol(V) \geq v_0$. This completes the proof. 
\qed 

Here is a simple application of the results of this section.

\begin{prop} \label{closedinX_0}
$D(X_0), D^*(X_0),$ and $D^*_0(X_0)$ are closed in $X_0$. 
\end{prop}

\pf Suppose that $\lim_n \chi_n = \chi_0 \in X_0$ where $\chi_n \in D(X_0)$ for all $n$. Proposition 1.4.4 of \cite{CS} (or Corollary 2.1 of \cite{CL}) shows that there are a subsequence $\{j\}$ of $\{n\}$ and a convergent sequence of representations $\{\rho_j\} \subset R_{X_0}$ such that $\chi_j = \chi_{\rho_j}$. Set $\rho_0 = \lim_j \rho_j$ and note that $\chi_0 = \chi_{\rho_0}$. Proposition \ref{convergence} implies that $\chi_{\rho_0} \in D(X_0)$. Moreover, if we assume that each $\chi_n \in D^*(X_0)$, then part (2) of Proposition \ref{convergence} implies that $\chi_{\rho_0} \in D^*(X_0)$. Thus $D(X_0)$ and $D^*(X_0)$ are closed in $X_0$. In particular, if $\chi_n \in D^*_0(X_0)$ for all $n$, then $\chi_0 \in D^*(X_0)$. From the previous lemma and the continuity of the volume function we have $|vol(\chi_0)| = \lim_n |vol(\chi_n)| \geq v_0 > 0$. Thus $\chi_0 \in D_0^*(X_0)$, which completes the proof. 
\qed

\subsection{Unbounded sequences of discrete $PSL_2(\mathbb C)$-characters} \label{unbounded} 

In this section $M$ will be a small knot manifold and $X_0$ a non-trivial component of $X_{PSL_2}(M)$. We are interested in the asymptotic behaviour of the sets $D(X_0)$ and $D^*(X_0)$. Consider a sequence $\{\chi_n\} \subset D(X_0)$ which converges to an ideal point $x_0$ of $X_0$. Fix $\rho_n \in R_{X_0}$ such that $\chi_n = \chi_{\rho_n}$ and let $\alpha_0$ be  the $\partial$-slope associated to $x_0$.  

\begin{lemma} \label{projectiveconvergence} 
Let $M$ be a small knot manifold, $X_0$ a curve component of $X_{PSL_2}(M)$, and $\{\chi_n\} \subset D(X_0)$ a sequence which converges to an ideal point $x_0$ of $X_0$. Fix $\rho_n \in R_{X_0}$ such that $\chi_n = \chi_{\rho_n}$ and let $\alpha_0$ be  the $\partial$-slope associated to $x_0$.  \\
$(1)$ For $n \gg 0$, $\hbox{kernel}(\rho_n|\pi_1(\partial M)) \cong \mathbb Z$ and $\rho_n(\pi_1(\partial M)) \cong \mathbb Z \oplus \mathbb Z/c_n$ where the $\mathbb Z$ factor is generated by a loxodromic element and $c_n \geq 1$.  \\ 
$(2)$ Let $\alpha_n \in H_1(\partial M)$ be the element, unique up to sign, which generates the kernel of $\rho_n|\pi_1(\partial M)$ $($ $n \gg0$ $)$. Then $\lim_n [\alpha_n] = [\alpha_0]$.  \\
$(3)$ If $[\alpha_n] \ne [\alpha_m]$ for some $m, n \gg 0$, then $X_0$ is a norm curve.  
\end{lemma}

\pf (1) Since $M$ is small, there is a slope $\beta$ on $\partial M$ such that $\tilde f_\beta$ has a pole at $x_0$. Thus $\rho_n(\beta)$ is loxodromic for large $n$ and so for such $n$, $\rho_n|\pi_1(\partial M))$ contains no parabolics. On the other hand, a discrete subgroup of $PSL_2(\mathbb C)$ isomorphic to $\mathbb Z^2$ contains parabolic matrices. Thus $\hbox{kernel}(\rho_n|\pi_1(\partial M)) \ne \{\pm I\}$, which implies (1). 

(2) Since $M$ is small, Proposition \ref{idealvalue} shows that  there is a unique slope $\alpha_0 \in H_1(\partial M)$ such that $\tilde f_{\alpha_{0}}(x_0) \in \mathbb C$. Further, $\alpha_0$ is a boundary slope, any surface $S$ associated to $x_0$ has non-empty boundary of slope $\alpha_0$, and if $\alpha_{0}^{*} \in H_1(\partial M)$ is a slope dual to $\alpha_0$ (i.e. $\alpha_0 \cdot \alpha_0^* = 1$), then $\lim_n f_{\alpha_{0}^{*}}(\chi_n)  = \infty$. We must show that $\lim_n [\alpha_n] = [\alpha_0]$. 

To that end set $\alpha_n = p_n \alpha_0 + q_n \alpha_{0}^{*}$. By construction,  
$(p_n, q_n) \ne (0,0)$ and for $n$ large, $\rho_n(\alpha_{0}^{*})$ is loxodromic. By choice of $\alpha_n$ we have $(\rho_n(\alpha_0))^{p_n} = (\rho_n(\alpha_{0}^{*}))^{-q_n}$ and therefore the minimal translation lengths $\ell(\rho_n(\alpha_0))$ and $\ell((\rho_n(\alpha_{0}^{*}))$ of $\rho_n(\alpha_0)$ and $\rho_n(\alpha_{0}^{*})$ satisfy: 
$$\vert p_n \vert \ell(\rho_n(\alpha_0)) = \vert q_n \vert \ell((\rho_n(\alpha_{0}^{*})) > 0.$$ 
If $\pm A \in PSL_2(\mathbb C)$, then $\ell(\pm A) = |\log(|\frac{\hbox{trace}(A)}{2} + \sqrt{(\frac{\hbox{trace}(A)}{2})^2 - 4}|)|$ and so our hypotheses imply that $\lim_n \ell((\rho_n(\alpha_{0}^{*})) = \infty$ while $\lim_n \ell(\rho_n(\alpha_0))$ is bounded. Thus $\lim_n \frac{q_n}{p_n} = 0$ or equivalently, $\alpha_n$ converge projectively to $[\alpha_0]$. 

(3) follows from (1) and Corollary \ref{normcondition}. 
\qed 

\begin{thm} \label{elementaryimage}
Let $M$ be a small knot manifold, $X_0$ a non-trivial component of $X_{PSL_2}(M)$, and $\{\chi_n\} \subset D(X_0)$ a sequence which converges to an ideal point $x_0$ of $X_0$. If $S_0$ is a component of an essential surface associated to $x_0$ and $i_\#: \pi_1(S_0) \to \pi_1(M)$ is the inclusion induced homomorphism, then either \\
\indent $(a)$ $\overline{i_\#^*(X_0)} \subset X_{{\cal N}}(S_0)$, or \\
\indent $(b)$ $i_\#^*(X_0) = \{\chi_\rho\}$ where $\rho(\pi_1(S_0))$ is either the tetrahedral group, the octahedral group, \\ \indent \hspace{4.5mm} or the icosahedral group.  
\end{thm}

\pf Fix $\rho_n \in R_{X_0}$ such that $\chi_n = \chi_{\rho_n}$ and let $S_0$ be a component of an essential surface $S$ in $M$ associated to $x_0$. Since $\chi_n|\pi_1(S_0)$ converges to a character $\chi_\sigma \in X_{PSL_2}(S_0)$ (Proposition \ref{limitsreps}), we can replace the $\rho_n$ by conjugate representations so that after passing to a subsequence $\{j\}$, we have $\lim \rho_j|\pi_1(S_0) = \sigma$ where $\chi = \chi_\sigma$ (see Proposition 1.4.4 of \cite{CS} or Corollary 2.1 of \cite{CL}). We also know that $\sigma$ is reducible (Proposition \ref{limitsreps}) and so by taking $j \gg 0$, $\rho_j^{S_0}$ is discrete and elementary (Proposition \ref{convergence}). A discrete elementary subgroup of $PSL_2(\mathbb C)$ is either reducible, conjugates into ${\cal N}$, or is isomorphic to a  {\it polyhedral group} (i.e. the tetrahedral group, the octahedral group, or the icosahedral group). Thus $i_\#^*(\chi_j)$ is contained in $X_{{\cal N}}(S_0)$ or $\rho_j$ has polyhedral image. This proves the lemma when $i_\#^*(X_0)$ is a single character. Suppose, on the other hand, that $\overline{i_\#^*(X_0)}$ is a curve $Y_0 \subset X_{PSL_2}(S_0)$. Then $i_\#^*: X_0 \to Y_0$ is finite-to-one and since there only finitely many characters of representations in $R_{PSL_2}(S_0)$ with image a polyhedral group, $Y_0 \cap X_{{\cal N}}(S_0)$ is infinite. But $X_{{\cal N}}(S_0)$ is Zariski closed in $X_{PSL_2}(S_0)$, and so it contains $Y_0$. 
\qed 

\begin{cor} \label{nonelementarycompact}
Let $M$ be a small knot manifold, $X_0$ a non-trivial component of $X_{PSL_2}(M)$. Suppose that for each ideal point $x_0$ of $X_0$ there are a component $S_0$ of an essential surface associated to $x_0$ and a character $\chi \in X_0$ such that $\chi|\pi_1(S_0)$ is non-elementary. Then $D(X_0), D^*(X_0)$, and $D_0^*(X_0)$ are compact. 
\end{cor}

\pf Theorem \ref{elementaryimage} shows that $D(X_0)$ does not accumulate to an ideal point of $X_0$. The result then follows from Proposition \ref{closedinX_0}. 
\qed

\begin{cor} \label{principalcompact} 
Let $M, N$ be small hyperbolic knot manifolds and suppose that $\varphi: \pi_1(M) \to \pi_1(N)$ is a virtual epimorphism. Fix a principal component $Y_0 \subset X_{PSL_2}(N)$ and set $X_0 = \varphi^*(Y_0)$. Then $D(X_0), D^*(X_0)$, and $D_0^*(X_0)$ are compact. In particular this is true for a principal component of $X_{PSL_2}(M)$. 
\end{cor}

\pf First suppose that $M = N$ and $\varphi$ is the identity. By Corollary \ref{nonelementarycompact} it suffices to show that for each connected essential surface $S_0$ in $M$, there is a character $\chi_\rho \in X_0$ such that $\rho(\pi_1(S_0))$ is non-elementary. Fix such a surface and note that $\pi_1(S_0)$ is a non-abelian free group since $M$ is small and hyperbolic. Moreover, since $X_0$ is principal, it contains the character of a discrete faithful representation $\rho_0$ of $\pi_1(M)$. Thus $\rho_0(\pi_1(S_0))$ is a discrete and free of rank at least $2$ and as such is non-elementary.  Thus $D(X_0), D^*(X_0)$, and $D_0^*(X_0)$ are compact.

Now consider the general case and let $\chi_\rho \in D(X_0)$. By Lemma \ref{closed} there is a $\chi_{\rho'} \in Y_0$ such that $\varphi^*(\chi_{\rho'}) = \chi_\rho$. Since $\chi_\rho$ is irreducible we can suppose that $\rho = \rho' \circ \varphi$. Then $\rho'$ is non-elementary and since the image of $\varphi$ has finite index in $\pi_1(N)$, it is also discrete. In other words, $\chi_{\rho'} \in D(Y_0)$ and so $\chi_\rho = \varphi^*(\chi_{\rho'}) \in \varphi^*(D(Y_0))$. Hence $D(X_0)$ is contained in the compact subset $\varphi^*(D(Y_0))$ of $X_0$.  
\qed 

\begin{exa} 
{\rm Corollary \ref{principalcompact} implies that the set of discrete, non-elementary characters in the $PSL_2(\mathbb C)$ character variety of the exterior of either a hyperbolic twist knot or the $(-2, 3, n)$ pretzel knot, $n \not \equiv 0$ (mod $3$) is compact (cf. Example \ref{twistpretzel}).}
\end{exa} 

\begin{rem} \label{pqnoncompact}
{\rm The corollary is false if we assume that $N$ is Seifert fibred but not a twisted $I$-bundle over the Klein bottle. Indeed, suppose that $N$ has base orbifold $D(p,q)$ where $p, q \geq 2, (p,q) \ne (2,2)$. Each pair $\pm I \ne A_0, B_0 \in {\cal D}$ such that $A_0^p = B_0^q = \pm I$ determines a curve $Y_0 \subset X_{PSL_2}(\mathbb Z/p * \mathbb Z/q)$ consisting of the characters of homomorphisms sending a generator of $\mathbb Z/p$ to $A_0$ and one of $\mathbb Z/q$ to a conjugate $B$ of $B_0$ (see Example 3.2, \cite{BZ1}). Further, if a sequence $\{B_n\}$ of such conjugates is chosen so that $\lim_n |\hbox{trace}(A_0 B_n)| = \infty$, the associated characters tend to the unique ideal point of $Y_0$ (Example 3.2, \cite{BZ1}). On the other hand, if $A_0, B_0 \in PSL_2(\mathbb R)$ are chosen to have extreme negative trace (page 293, \cite{Kn}), they generate a discrete group isomorphic to $\mathbb Z/p * \mathbb Z/q \cong \pi_1(D(p,q))$ as long as $|\hbox{trace}(A_0B)| \geq 2$ (Theorem 2.3 \cite{Kn}). In particular, they determine a principal component $Y_0 \subset X_{PSL_2}(\mathbb Z/p * \mathbb Z/q) = X_{PSL_2}(\pi_1(D(p,q)) \subset X_{PSL_2}(N)$ for which $D(Y_0)$ is non-compact. By hypothesis, $Y_0 \subset X_+^{str}(N)$, so $X_0 := \varphi^*(Y_0) \subset X_+^{str}(M)$ (Corollary \ref{1-1}) and by construction, $D(X_0)$ is non-compact.}
\end{rem}

\begin{thm} \label{torsionfreeconvideal}
Let $M$ be a small knot manifold, $X_0$ a norm curve component of $X_{PSL_2}(M)$, and $\{\chi_n\} \subset D^*(X_0)$ a sequence which converges to an ideal point $x_0$ of $X_0$. If $S$ is an essential surface associated to $x_0$ and $S_0$ a component of $S$, then $S_0$ is separating and there is a complementary component $A$ of $S_0$ such that $\rho(\pi_1(A)))$ is abelian for each $\rho \in R_{X_0}$. There is a subsequence $\{j\}$ of $\{n\}$ such that $\rho_j(\pi_1(A)))$ is cyclic for all $j$.   
\end{thm}

\pf Choose $\rho_n \in R_{X_0}$ such that $\chi_n = \chi_{\rho_n}$. By Lemma \ref{projectiveconvergence} we may suppose that $\rho_n(\pi_1(\partial M))$ is loxodromic and since $\rho_n$ is torsion free, the lemma imlies that there is a unique slope $\alpha_n$ on $\partial M$ satisfying $\rho_n(\alpha_n) = \pm I$. We may suppose that the $\alpha_n$ are distinct, since $X_0$ is a norm curve, and that none of them are boundary slopes \cite{Hat}. 

Since $M$ is small, $S_0$ has non-empty boundary of slope $\alpha_0$, say. Then by construction, $\rho_n(\alpha_0) \ne \pm I$ is loxodromic for $n \gg 0$. According to Theorem \ref{elementaryimage}, $\rho_n(\pi_1(S_0))$ is elementary and since it is discrete, torsion free, and contains a loxodromic ($n \gg 0$), it is a cyclic subgroup of $PSL_2(\mathbb C)$. In this case we can apply the bending construction to $\chi_n$ along $\pi_1(S_0)$ (Appendix \ref{bending}). We claim that for $n \gg 0$, the bending of $\chi_n$ along $\pi_1(S_0)$ is trivial. For such $n$, $\chi_{n} \in X_0$ is contained in a unique component of $X_{PSL_2}(M)$, and so if the claim is false $X_0$ is obtained by bending $\chi_{n}$ along $\pi_1(S_0)$. In particular $\rho_m(\alpha_0)$ is independent of $m$, at least up to conjugation. But then $f_{\alpha_0}|X_0$ is constant and so $X_0$ cannot be a norm curve, contrary to our hypotheses. 

Suppose that $S_0$ is non-separating and write $M = A/ \{S_0^+ = S_0^-\}$ where $A$ is the complementary component of $S_0$ in $M$ and $S_0^+ \sqcup S_0^- \subseteq \partial A$ are parallel copies of $S_0$. Note that $\pi_1(M)$ is generated by $\pi_1(A)$ and $\gamma$, a homotopy class represented by a loop which intersects $S_0$ once transversely. Fix $n \gg 0$ and note that since $\rho_n$ cannot be bent non-trivially along $\pi_1(S_0)$, either $\rho_n(\pi_1(M)) \subset {\cal N}$ or $\rho_n$ is reducible (Lemma \ref{constbendnonsep}). As neither of these possibilities is satisfied in our situation, $S_0$ must be separating. Hence if $M = A \cup_{S_0} B$ where $A$ and $B$ are the complementary components of $S_0$ in $M$, the fact that for large $n$ the bending of $\chi_n$ along $\pi_1(S_0)$ is trivial, at least one of $\rho_n^A, \rho_n^B$ has cyclic image. Further, this image is trivial if the image of $\rho_n^{S_0}$ is trivial (Lemma \ref{constconj}). By passing to a subsequence and possibly exchanging $A$ and $B$, we can assume that for $n \gg 0$, $\rho_n^A$ has cyclic image and $\rho_n(\pi_1(A)) = \{\pm I\}$ if $\rho_n(\pi_1(S_0)) = \{\pm I\}$. 

Let ${\cal O}(\rho_n)$ denote the $PSL_2(\mathbb C)$ orbit of $\rho_n$. Since $\cup_{m \geq n} {\cal O}(\rho_n)$ is Zariski dense in $R_{X_0}$, $n \geq 1$, the previous paragraph shows that $\rho(\pi_1(S_0))$ is abelian for each $\rho \in R_{X_0}$. 
\qed

\begin{cor} \label{lessthan2}
Let $M$ be a small knot manifold and $X_0$ a norm curve of $X_{PSL_2}(M)$. Then $D^*(X_0)$ is a compact subset of $X_0$ as long as the following condition holds: Any ideal point of a norm curve in $X_{PSL_2}(M)$ has an associated essential surface with a component $S_0$ having no more than two boundary components. 
\end{cor}

\pf By Proposition \ref{closedinX_0}, it suffices to show that $D^*(X_0)$ is contained in a compact subset of $X_0$. Suppose then that $\{\chi_n \} \subset D^*(X_0)$ is a sequence which converges to an ideal point $x_0$ of $X_0$ and choose $\rho_n \in R_{X_0}$ whose character is $\chi_n$. Fix a component $S_0$ of an essential surface $S$ associated to $x_0$ with $|\partial S_0| \leq 2$. Theorem \ref{torsionfreeconvideal} implies that $\rho(\pi_1(A)))$ is abelian for each $\rho \in R_{X_0}$. Since $X_0$ is a norm curve, $f_{\alpha_0}|X_0$ is non-constant and so there is a Zariski dense subset in $R_{X_0}$ of representations $\rho \in R_{X_0}$ such that $\rho(\alpha_0)$ is loxodromic. Fix such a representation $\rho_0$ and conjugate it so that $\rho_0(\alpha_0)$ is diagonal. Then both $\rho_0(\pi_1(A))$ and $\rho_0(\pi_1(\partial M))$ are contained in ${\cal D}$. 

Let $B$ be the other complementary component of $S_0$. A maximal compression of $\partial B$ in $B$ must yield a family of $2$-spheres as $M$ is small. Thus $B$ is a handlebody and therefore $\pi_1(\partial B) \to \pi_1(B)$ is surjective. Consider a class $\sigma \in \pi_1(M)$ represented by a product of a path in $S_0$ followed by one in $\partial M \cap B$. By hypothesis $|\partial S_0| = 2$ and so $\sigma$ 
is the product of an element of $\pi_1(A)$ and one in $\pi_1(\partial M)$. It follows that $\rho_0(\sigma) \in {\cal D}$. Since $\pi_1(B)$ is generated by such classes and $\pi_1(S_0)$, we see that $\rho_n(\pi_1(B)) \subset {\cal D}$. But then the image of $\rho_0$ is abelian, which is impossible as $R_{X_0}$ contains a Zariski dense subset of such representations. Thus $D^*(X_0)$ contained in a compact subset of $X_0$. 
\qed

Ohtsuki \cite{Oht} has shown that two-bridge knot exteriors satisfy the condition of the previous corollary.  

\begin{cor} \label{2bridge compact} 
Let $X_0$ be a norm curve in the character variety of a two-bridge knot exterior. Then $D^*(X_0)$ and $D_0^*(X_0)$ are compact subsets of $X_0$. 
\qed
\end{cor}

\subsection{Convergent sequences of discrete, co-compact $PSL_2(\mathbb C)$-characters with non-zero volume} \label{convergent} 

Let  $M$ will be a small hyperbolic knot manifold and $X_0$ a non-trivial component of $X_{PSL_2}(M)$. Consider a sequence $\{\chi_{\rho_n}\} \subset D_0^*(X_0)$ of distinct characters which converge to some $\chi_{\rho_0} \in  X_0$. Fix $\rho_n \in R_{X_0}$ whose character is $\chi_n$, set $\Gamma_n = \rho_n(\pi_{1}(M))$, and let $V_n = \mathbb H^3/\Gamma_n$. 

Since $M$ is small, $\|\cdot \|_{X_0} \ne 0$ and so there are only finitely many $n$ such that $\rho_n(\pi_1(\partial M))$ is either $\{\pm I\}$ or contains a paraboic element of $PSL_2(\mathbb C)$. (Otherwise $\rho_n(\pi_1(\partial M))$ would be contained in a paraboic subgroup for infinitely many $n$ and therefore $f_\gamma|X_0 \equiv 0$ for every peripheral $\gamma$). We suppose then that $\rho_n(\pi_1(\partial M))$ contains a loxodromic for each $n$. Since $\Gamma_n$ is torsion free, this implies that $\rho_n(\pi_1(\partial M)) \cong \mathbb Z$ and so there is a unique slope $\alpha_n$ on $\partial M$ which generates $\hbox{kernel}(\rho_n|\pi_{1}(\partial M))$. It follows from Lemma \ref{standardimage} that $V_n$ is a closed hyperbolic $3$-manifold. If $\bar \rho_n \in R_{PSL_2}(M(\alpha_n))$ is the homomorphism induced by $\rho_n$, the proof of part (3) of this lemma shows that  
$$vol(V_n) \leq |vol(\bar \rho_n)| \leq |vol(\rho_n)| \leq vol(M(\alpha_n)) < vol(M).$$

\begin{thm} \label{convord}
Assume that the sequence $\{\chi_n\} \subset D_0^*(X_0)$, as above, converges to a character $\chi_{\rho_0}$ and that the slopes $\alpha_n$ associated to $\rho_n$ are distinct. Then there are 
\begin{description}
\vspace{-.5cm} 
\item[{\rm $\;$ (a)}] a subsequence $\{j\}$ of $\{n\}$ such that $\{V_{j}\}$ converges geometrically to a $1$-cusped hyperbolic $3$-manifold $V$ whose fundamental group containes $\rho_0(\pi_{1}(M))$ as a finite index subgroup. 
\vspace{-.3cm} 
\item[{\rm $\;$ (b)}] a proper non-zero degree map $f_0: M \to V_0$ such that $V_0$ is a compact core of $V$  and if $k_0 : V_0 \to V$ is the inclusion, then $\rho_0 = (k_0)_\# \circ (f_0)_\#$.  
\vspace{-.3cm} 
\item[{\rm $\;$ (c)}] slopes $\beta_j$ on $\partial V_0$ and identifications $V_j \cong V_0(\beta_j)$, such that $(f_0|\partial M)_*(\alpha_j)$ is a multiple of $\beta_j \in H_1(\partial V)$ and if $k_j: V_0 \to V_0(\beta_j)$ is the inclusion, then $\chi_j$ is the character of the composition $(k_j)_\# \circ (f_0)_\#$.
\vspace{-.3cm} 
\item[{\rm $\;$ (d)}] non-zero degree maps $f_j: M(\alpha_n) \to V_0(\beta_j)$ such that the following diagrams are commutative up to homotopy: 
$$\begin{array}{ccc} 
M & \stackrel{f_0}{\longrightarrow} & V_0 \\
\downarrow & & \downarrow \\
M(\alpha_j) & \stackrel{f_j}{\longrightarrow} & V_j \cong V_0(\beta_j) 
\end{array}$$
\end{description}
Moreover, $X_0 = (f_0)_\#^*(Y_0)$ where $Y_0$ is a principal curve for $V_0$. In particular, $X_0$ is a norm curve. 
\end{thm}

\pf After replacing the $\rho_n$ by conjugate representations and passing to a subsequence, we may suppose that $\lim \rho_n = \rho_0$ (see Corollary 2.1 of \cite{CL} for example). Then $\{\Gamma_n\}$ converges algebraically to $\Gamma_0 = \rho_0(\pi_1(M))$. Proposition \ref{convergence}, $\Gamma_0$ is a non-elementary Kleinian group and there are homomorphisms $\theta_n: \Gamma_0 \to \Gamma_{n}$ such that $\rho_{n} = \theta_n \circ \rho_0$ for $n \gg 0$. By passing to a subsequence we may suppose that this is true for $n \geq 1$.  

Since ${\rm kernel}(\rho_n|\pi_1(\partial M)) = \langle \alpha_n \rangle$, $\rho_{n} =  \theta_n \circ \rho_0$, and the slopes $\alpha_n$ are distinct, it follows that $\rho_0|\pi_1(\partial M)$ is injective. By Proposition \ref{convergence} implies that after passing to a subsequence we can suppose that 

\indent (i) $\{\Gamma_n\}$ converges geometrically to a non-elementary Kleinian group $\Gamma$ containing $\Gamma_0$ \\ \indent \hspace{5mm} and $\theta_n$ extends to homomorphisms $\Gamma \to \Gamma_n$ which we still denote by $\theta_n$.

\indent (ii) $\lim V_n = V := \mathbb H^3/\Gamma$ in the sense of Gromov bilipschitz topology. 

\noindent As we noted above, $vol(V_n) < vol(M)$ and therefore $vol(V) = \lim vol(V_n) \leq vol(M)$. It follows that $V$ is a complete, connected, orientable, finite volume hyperbolic $3$-manifold. Further, $V$ has at least one cusp since $\Gamma$ contains $\rho_0(\pi_1(\partial M)) \cong \mathbb Z \oplus \mathbb Z$. Thus $V_n$ is obtained from $V$ by hyperbolic  Dehn filling for large $n$ (cf. \S \ref{convkleinhyp}). 

Let $k_0: V_0 \to V$ be the inclusion of a compact core of $V$. Since $M$ and $V_0$ are $K(\pi,1)$-spaces, there is a map $f_0: M \to V_0$ such that $\rho_0 =  (k_0)_\# \circ (f_0)_\#$. (We have fixed an identification of $\pi_1(V)$ with $\Gamma$ here.) Since $V_0$ is atoroidal, there is a torus $T \subseteq \partial V_0$ such that $\rho_0(\pi_1(\partial M)) \subseteq \pi_1(T)$, at least up to conjugation. Homotope $f_0$ so that $f_0(\partial M) \subseteq  T$. Then $(f_0)\#: \pi_1(\partial M) \to \pi_{1}(T)$ and is injective (by construction), which shows that ${\rm degree}(f_{0}|: \partial M \to T) \ne 0$. On the other hand, if $[\partial M] \in H_2(M)$ and $[T]  \in H_2(V_0)$ are fundamental classes for $\partial M$ and $T$, then $0 = [\partial M]$ so that $0 = (f_0)_*([\partial M]) = {\rm degree}(f_{0}|) [T]$ and therefore $\partial V_0 = T$. 

Recall that $V_n$ is obtained by hyperbolic Dehn filling on $V$. There is some slope $\beta_n$ on $T$ such that $V_n = V_0(\beta_n)$. If $k_n: V_0 \to V_n$ denotes the inclusion, then  $\theta_n \circ (k_0)_\# =  (k_n)_\#$ (cf. Theorem 7.17  of [MT]) and so $\mbox{kernel} (\theta_n) = (k_0)_\#(\langle \langle \beta_n \rangle \rangle_{\pi_1(V_0)})$ where $\langle \langle \beta_n \rangle \rangle_{\pi_1(V_0)}$ is the normal closure in $\pi_1(V_0)$ of the element corresponding to the slope $\beta_n$. Thurston's hyperbolic Dehn filing theorem (see chapter 5 of \cite{Thu} or the appendix to \cite{BoP}) implies that $\pi_1(\partial V_0) \cap \mbox{kernel} (\theta_n) = \langle \beta_n \rangle \cong \mathbb Z$ for large $n$. By passing to a subsequence we can arrange for it to hold for all $n$. 
 
Since $\rho_n = \theta_n \circ \rho_0 = \theta_n \circ  (k_0)_\# \circ (f_0)_\# =  (k_n)_\# \circ (f_0)_\#$, we have $1 = \rho_n(\alpha_n) = (k_n)_\# ((f_0)_\#(\alpha_n))$ and therefore $(f_0)_\#(\alpha_n) \in \pi_1(\partial V_0) \cap \mbox{kernel} (\theta_n) =\langle \beta_n \rangle$. Thus $f_0$ induces a map $f_n: M(\alpha_n) \to V(\beta_n)$ with ${\rm degree}(f_n) = {\rm degree}(f_0)$. If $i_n: M \to M(\alpha_n) $ denotes the inclusion, we have $(f_n)_\# \circ (i_n)_\# = \theta_n \circ (k_0)_\#\circ (f_{0})\# =  (k_n)_\# \circ (f_0)\# = \rho_n$. Hence for large $n$ the following diagrams are commutative up to homotopy  
$$\begin{array}{ccc} 
M & \stackrel{f_0}{\longrightarrow} & V_0 \\
\downarrow & & \downarrow \\
M(\alpha_j) & \stackrel{f_j}{\longrightarrow} & V_j \cong V_0(\beta_j) 
\end{array}$$

To complete the proof, we must show that $X_0 = (f_0)_\#^*(Y_0)$ where $Y_0$ is a principal curve for $V_0$. To that end we note that Thurston's hyperbolic Dehn filling theorem proves that if $Y_0$ is the principal component of $X_{PSL_2}(V_0)$ which contains the character $\chi_0'$ of $\pi_1(V_0) = \Gamma$, then $\chi_0' = \lim_n \chi_n'$ where $\chi_n' \in Y_0$ is the character of our identification $\pi_1(V_0(\beta_n) = \Gamma_n$. By construction $(f_0)_\#^*(\chi_n') = \chi_n$ so that $X_0 \cap (f_0)_\#^*(Y_0)$ is infinite. Lemma \ref{closed} then shows that $X_0 = (f_0)_\#^*(Y_0)$. 
\qed

\begin{cor} \label{discetecharsprincipal} 
Let $M$ be a small hyperbolic knot manifold and $X_0$ a principal component of $X_{PSL_2}(M)$. Then all but finitely many of the elements of $D_0^*(X_0)$ are induced by the complete hyperbolic structure on the interior of $M$ or by Dehn fillings of manifolds finitely covered by $M$.
\end{cor}

\pf We know from Corollary \ref{principalcompact} that $D_0^*(X_0)$ is contained in a compact subset of $X_0$, Thus if the result is false, we could find a convergent sequence $\{\chi_n\} \subset D_0^*(X_0)$ of distinct characters no one of which is induced by a holonomy character of $M$ or one of the Dehn fillings of an oriented manifold it finitely covers. As above we can assume that $\chi_n$ is the character of a representation $\rho_n$ which is peripherally non-trivial. Let $\alpha_n$ be its slope and note that since $\|\cdot\|$ is a norm curve, the function $n \mapsto \alpha_n$ is finite-to one. Thus we can take a subsequence $\{j\}$ of $\{n\}$ for which the $\alpha_j$ are distinct and apply Theorem \ref{convord} to see that there are a non-zero degree map $f_0: M \to N$ where $N$ is hyperbolic, a principal component $Y_0$ of $X_{PSL_2}(N)$ such that $X_0 = (f_0)_\#^*(Y_0)$, and, for infinitely many $j$, $\chi_j$ is the image under  $(f_0)_\#^*$ of the holonomy character of some Dehn filling of $N$. Lemma \ref{kerprincipal}(2) shows that $(f_0)_\#$ is injective and so we can take $f_0: M \to N$ to be a covering map (\cite{Wal2}). But then infinitely many $\chi_j$ are induced from a Dehn filling of an orientable manifold covered by $M$ contrary to our hypotheses. This completes the proof. 
\qed

\begin{cor} \label{characterizingdiscretecharacters}
Let $M$ be a small knot manifold and suppose that there is a norm curve $X_0$ in $X_{PSL_2}(M)$ for which $D_0^*(X_0)$ has an accumulation point in $X_0$. Then there are a hyperbolic manifold $N$, a non-zero degree map $f_0: M \to N$, and a principal component $Z_0$ of $X_{PSL_2}(N)$ such that $X_0 = (f_0)_\#^*(Z_0)$. Further, all but finitely many characters in $D_0^*(X_0)$ are the images under $(f_0)_\#^*$ of the holonomy character of $N$ or one of the Dehn fillings of an oriented manifold finitely covered by $N$.   
\end{cor}

\pf The hypotheses can be used with Theorem \ref{convord} to see that there are a hyperbolic manifold $N_0$, a non-zero degree map $f_0: M \to N_0$, and a principal component $Y_0$ of $X_{PSL_2}(N_0)$ such that $X_0 = f_\#^*(Y_0)$. Let $N \to N_0$ be the cover corresponding to the image of $(f_0)_\#$, $\tilde f_0: M \to N$ a lift of $f_0$, and $Z_0$ the principal curve in $X_{PSL_2}(N)$ obtained by restriction from $Y_0$. Clearly $X_0 = (\tilde f_0)_\#^*(Z_0)$ and final claim of the corollary is a consequence of the previous result applied to $Z_0$. The details are left to the reader.
\qed 

The final results of this section follow immediately from Theorem \ref{torsionfreeconvideal} and Corollaries \ref{lessthan2}, \ref{2bridge compact}, and \ref{characterizingdiscretecharacters}. 

\begin{cor}
Let $M$ be a small knot manifold and suppose that there is a norm curve $X_0$ in $X_{PSL_2}(M)$ for which $D_0^*(X_0)$ has an accumulation point in $X_0$. Then $D_0^*(X_0)$ is compact and has a unique accumulation point corresponding the holonomy character of a hyperbolic knot manifold under a non-zero degree map $M \to N$. 
\qed 
\end{cor} 

\begin{cor} \label{compactconditions} 
Let $M$ be a small knot manifold and $X_0 \subset X_{PSL_2}(M)$ a norm curve. Then $D_0^*(X_0)$ is compact in $X_0$ with at most one accumulation point if for each ideal point $x_0$ of $X_0$ there is a component $S_0$ of an essential surface associated to $x_0$ such that at least one of the following two conditions holds: \\
\indent $(i)$ $\chi|\pi_1(S_0)$ is non-elementary for some $\chi \in X_0$; or \\ 
\indent $(ii)$ $|\partial S_0| \leq 2$.  
\qed
\end{cor}

\begin{thm} \label{twobridgecompact} 
Let $X_0$ be a norm curve in the character variety of a two-bridge knot exterior $M$. Then $D_0^*(X_0)$ is either finite or is a compact subset of $X_0$ with a unique accumulation point. In the latter case there are a two-bridge knot exterior $N$, a non-zero degree map $f: M \to N$, and a principal component $Y_0$ of $X_{PSL_2}(N)$ such that $X_0 = f_\#^*(Y_0)$. 
\end{thm}

\pf The theorem follows from the results cited above and Theorem \ref{2bridgemini} once we note that $M$ must be hyperbolic if $X_{PSL_2}(M)$ is to contain a norm curve. 
\qed

\subsection{Domination and hyperbolic Dehn filling} \label{dominationhyperbolic}

In this section we prove Theorem \ref{thm:domination}. 

Let $M$ be a small knot manifold and $\{\alpha_n\}_{n \geq 1}$ a sequence of distinct slopes on $\partial M$ such that for each $n$ we have a map $f_n: M(\alpha_n) \to V_n$ of degree $d_n \geq 1$ where $V_n$ is hyperbolic. We suppose as well that $\{\alpha_n\}$ does not subconverge projectively to a boundary slope and that there is a slope $\alpha_0$ on $\partial M$ such that $M(\alpha_0)$ does not dominate any hyperbolic $3$-manifold. 

Let $p_n: \tilde V_n \to V_n$ be the finite cover corresponding to $(f_n)_\#(\pi_1(M))$. We can suppose that $p_n$ is a local isometry. Fix a lift $\tilde f_n: M \to \tilde V_n$ of $f_n$ of degree $\tilde d_n \geq 1$ say. If $v_0 > 0$ is the minimal volume for closed, connected, orientable, hyperbolic $3$-manifolds,  then for each $n$ we have 
$vol(M) > vol(M(\alpha_n)) \geq \tilde d_n vol(\tilde V_n) = d_n vol(V_n) \geq d_n v_0 \geq \tilde d_n v_0.$  Thus the $d_n$ and $\tilde d_n$ are bounded so we can assume, after passing to a subsequence, that they are constant, say  $\hbox{degree}(f_n) = d \geq 1, \; \hbox{degree}(\tilde f_n) = \tilde d.$ The degree of each $p_n$ is $d/ \tilde d$. 

We identify $\pi_1(V_n)$ with a subgroup $\Gamma_n$ of $PSL_2(\mathbb C)$ and set $(p_n)_\#(\pi_1(\tilde V_n)) = \tilde \Gamma_n \subseteq \Gamma_n$. Let $i_n: M \to M(\alpha_n)$ be the inclusion and define $\rho_n \in R_{PSL_2}(M)$ to be the composition $\pi_1(M) \stackrel{(i_n)_\#}{\longrightarrow} \pi_1(M(\alpha_n)) \stackrel{(f_n)_\#}{\longrightarrow} \tilde \Gamma_n \subseteq \Gamma_n \subset PSL_2(\mathbb C)$. The character of $\rho_n$ will be denoted by $\chi_n$. These objects combine in the following commutative diagram. 

\vspace{-.05cm}

\begin{diagram}[height=3.4em,w=6em]
& & \pi_1(\tilde V_n)   \\
&  \ruOnto(2,2){ (\tilde f_n)_\#} & \dEmbed_{(p_n)_\#}    \\
 \pi_1(M(\alpha_n))   & \rTo^{\;\;\;\;\; (f_n)_\#} & \Gamma_n =\pi_1(V_n)   \\
\uOnto^{(i_n)_\#} & &   \dEmbed \\
\pi_1(M) & \rTo^{\;\;\;\;\; \rho_n} &  PSL_2(\mathbb C)
\end{diagram}

\vspace{.2cm} 
We claim that each $\rho_n$ is peripherally non-trivial. Otherwise the composition $M \to M(\alpha_n) \stackrel{f_n}{\longrightarrow} V_n$ extends to a map $M(\alpha_0) \to V_n$ of degree $d$, contrary to our hypothesis. Since $\alpha_n \ne \alpha_m$ for $n \ne m$, the characters $\chi_n$ are distinct. After passing to a subsequence we can suppose that they are contained in a non-trivial curve $X_0 \subset X_{PSL_2}(M)$. Proposition \ref{normcondition} shows that $X_0$ is a norm curve. Finally, noting that $vol(\chi_n) = d vol(V_n) \ne 0$ we see that $\chi_n \in D_0^*(X_0)$. Since the slopes $\{\alpha_n\}_{n \geq 1}$ do not projectivly subconverge to a $\partial$-slope, Lemma \ref{projectiveconvergence} shows that there is a subsequence of characters $\{\chi_k\}$ which converge to a character $\chi_{\rho_0} \in D_0^*(X_0)$, and so the conditions of Theorem \ref{convord} are satisfied. 

\noindent{\bf Proof of Theorem \ref{thm:domination}}. By Theorem \ref{convord}, the sequence $\tilde V_k$ converges geometrically to a $1$-cusped hyperbolic $3$-manifold $\tilde V$ for which there are: 
\vspace{-.35cm} 
\begin{description}
\item[{\rm $\;$ (a)}] a proper non-zero degree map $\tilde f_0: M \to \tilde V_0$ such that $\tilde V_0$ is a compact core of $\tilde V$ and if $j_0 : \tilde V_0 \to \tilde V$ is the inclusion map, then $\rho_0 = (j_0)_\# \circ (\tilde f)_\#$; 
\vspace{-.95cm} 
\item\item[{\rm $\;$ (b)}] slopes $\tilde \beta_k$ on $\partial \tilde V_0$ and identifications $\tilde V_k = \tilde V_{0}(\tilde \beta_k)$, such that $(\tilde f_0|\partial M)_*(\alpha_k)$ is a multiple of $\tilde \beta_k \in H_1(\partial \tilde V_0)$ and if $\tilde j_k: \tilde V_0 \to \tilde V_0(\tilde \beta_k)$ is the inclusion, then $\chi_k$ is induced by the composition $(p_k)_\# \circ (\tilde j_k)_\# \circ (\tilde f_0)_\#$.
\vspace{-.95cm} 
\item\item[{\rm $\;$ (c)}] non-zero degree maps $\tilde f'_k: M(\alpha_k) \to \tilde V(\tilde \beta_k)$ such that  the following diagrams are commutative up to homotopy: 
$$\begin{array}{ccc} 
M & \stackrel{\tilde f_0}{\longrightarrow} & \tilde V_0 \\
\downarrow & & \downarrow \\
M(\alpha_k) & \stackrel{\tilde f'_k}{\longrightarrow} & \tilde V_k \cong \tilde V_0(\tilde \beta_k) 
\end{array} \eqno{(4.2.1)}$$
\end{description}
Since  $(p_k)_\# \circ (\tilde f_k)_\# \circ (i_k)_\# = \rho_k = (p_k)_\# \circ (\tilde j_k)_\# \circ (\tilde f_0)_\# = (p_k)_\# \circ (\tilde f'_k)_\# \circ (i_k)_\# $, it follows that $f_k = p_k \circ \tilde f_k$ is homotopic to  $f'_k = p_k \circ \tilde f'_k$. In particular $\hbox{degree}(\tilde f_0) = \hbox{degree}(\tilde f'_k) = \hbox{degree}(\tilde f_k) = \tilde d$ and $\hbox{degree}(f'_k) = \hbox{degree}(f_k) = d$.

Now $\lim_k vol(V_k) = \lim_k (\frac{\tilde d}{d}) vol(\tilde V_k) =  (\frac{\tilde d}{d}) vol(\tilde V)$, and so after passing to a subsequence we may assume that $\{V_k\}$ converges geometrically to a complete hyperbolic $3$-manifold $V$ with finite volume $vol( V) =  (\frac{\tilde d}{d}) vol(\tilde V)$. For $k \gg 0$, $vol(\tilde V) > vol(\tilde V_k)$ and therefore $vol(V) > vol(V_k)$. Thus $V$ has at least one cusp. On the other hand, $p_k$ is a local isometry so for $\tilde x \in \tilde V_k$ we have $\hbox{inj}(\tilde x) \leq (\frac{d}{\tilde d}) \hbox{inj}(p_k(\tilde x))$. Thus if $\mu_0$ is the Margulis constant and $\mu \leq \frac{\tilde d \mu_0}{d}$, we have $p_{k}^{-1}((V_{k})_{(0,\mu]}) \subseteq (\tilde V_k)_{(0,\frac{d \mu}{\tilde d}]}$ ($k \gg 0$). Since there is a sequence $\mu_k \to 0$ such that $(\tilde V_k)_{(0,\frac{d 
\mu_k}{\tilde d}]}$ is a Margulis tube about a geodesic $\tilde \gamma_k$, 
$(V_{k})_{(0,\mu_k]}$ is a Margulis tube about a geodesic $\gamma_k$, because 
a geodesic is unique in its homotopy class. Thus $V$ 
has only one cusp. We note, moreover, that $p_k^{-1}(\gamma_k) = \tilde \gamma_k$ and therefore $\overline{(\tilde V_k)_{(0,\frac{d \mu}{\tilde d}]} \setminus p_{k}^{-1}((V_{k})_{(0,\mu]})} \cong \partial (\tilde V_k)_{(0,\frac{d \mu}{\tilde d}]} \times I$. Thus for large $k$ we can identify $(V_k)_{[\mu, \infty)}$ with a compact core $V_0$ of $V$ and $p_{k}^{-1}((V_k)_{[\mu, \infty)})$ with a compact core $\tilde V_0$ of $\tilde V$. In this way $p_k$ induces a covering map $p_k^0: \tilde V_0 \to V_0$ of degree $d/ \tilde d$. Since $V_0$ and $\tilde V_0$ admit complete finite volume hyperbolic structures on their interiors, after pre-composition by an isotopy of $\tilde V_0$, we can take $p_n^0$ to be a local isometry on the interior of $\tilde V_0$. Now $V_0$ has only finitely many (pointed) covers of degree $d/ \tilde d$ up to equivalence and the isometry group of $\hbox{int}(\tilde V_0)$ is finite, therefore we can restrict to a subsequence and suppose that for all $n, m$, we have $p_n^0 = p_m^0 = p$, say. 

The geometric convergence of $V_k$ to $V$ implies that for large $k$ there are slopes $\beta_k$ on $\partial V_0$ such that $V_k = V_0(\beta_k)$. From the previous paragraph we see that any component of $p^{-1}(\beta_k)$ is isotopic to $\tilde \beta_k$ on $\partial \tilde V_0$. Therefore the following diagrams  are commutative up to homotopy:

$$\begin{array}{ccc} 
\tilde V_0 & \stackrel{p}{\longrightarrow} & V_0 \\
\downarrow & & \downarrow \\
\tilde V_k \cong \tilde V_0(\tilde\beta_k) & \stackrel{p_k}{\longrightarrow} & V_k \cong  V_0(\beta_k) 
\end{array} \eqno{(4.2.2)}$$

\noindent Since $f_k = p_k \circ \tilde f_k$ is homotopic to $f'_k = p_k \circ \tilde f'_k$, by putting together diagrams (4.2.1) and (4.2.2) one deduces that  the  proper map $f = p \circ \tilde f_0: M \to V_0$ of degree $d \geq 1$ makes the following diagrams commute up to homotopy:

$$\begin{array}{ccc} 
M & \stackrel{f}{\longrightarrow} & V_0 \\
\downarrow & & \downarrow \\
M(\alpha_k) & \stackrel{f_k}{\longrightarrow} &  V_k \cong V_0(\beta_k) 
\end{array}$$

\noindent If we assume further that the dominations $f_k : M(\alpha_k) \to V_k \cong  V(\beta_k)$ are strict, then the domination $f_0: M \to V$ must also be strict. Otherwise $V_0$ is homeomorphic to $M$ and $d = \hbox{degree}(f_0) = 1$ (since $vol(M) = vol(V_0)$). Then $(f_0)_\#: \pi_{1} (M) \to  \pi_{1} (V_0)$ is surjective and therefore an isomorphism since $\pi_{1} (M)$ is Hopfian. By Mostow rigidity theorem we can suppose that $f_0$ is a homeomorphism. But then the induced maps $f'_k : M(\alpha_k) \to V(\beta_k)$ are homeomorphisms homotopic to $f_k$, in contradiction with our assumption that $f_k$ is a strict domination. 
\qed

\section{${\cal H}$-minimal Dehn filling}\label{sec:h-minimal}

The goal of this section is to construct collections of infinitely many ${\cal H}$-minimal closed hyperbolic $3$-manifolds by proving Theorem \ref{thm:h-minimal} and Corollaries \ref{cor:integer surgery} and \ref{cor:2-bridge}.
The proofs of these results rely on the following theorem:

\begin{thm} \label{minimal hyperbolic}
Let $M$ be an ${\cal H}$-minimal, small hyperbolic knot manifold and suppose that there is a slope $\alpha_0$ on $\partial M$ such that the Dehn filled manifold $M(\alpha_0)$ does not dominate any closed hyperbolic manifold. \\ 
$(1)$ If $U \subset \mathbb P(H_1(\partial M; \mathbb R))$ is the union of disjoint closed arc neighbourhoods of the finite set of boundary slopes of $M$, then $\mathbb P(H_1(\partial M; \mathbb R)) \setminus U$ contains only finitely many projective classes of slopes $\alpha$ such that $M(\alpha)$ is not ${\cal H}$-minimal. In particular, $M$ admits infinitely many ${\cal H}$-minimal Dehn fillings. \\
$(2)$ If $D_0^*(X_0)$ is a compact subset of $X_0$ for each norm curve in $X_{PSL_2}(M)$, then there are only finitely many slopes $\alpha$ on $\partial M$ such that $M(\alpha)$ is not ${\cal H}$-minimal. In particular, this conclusion holds if for each ideal point $x_0$ of a norm curve $X_0$, there is a component $S_0$ of an essential surface associated to $x_0$ such that at least one of the following two conditions holds: \\
\indent $(i)$ $\chi|\pi_1(S_0)$ is non-elementary for some $\chi \in X_0$; or \\ 
\indent $(ii)$ $|\partial S_0| \leq 2$.  

\end{thm}

\pf (1) Suppose that there are infinitely many projective classes of slopes $\alpha$ in 
$\mathbb P(H_1(\partial M; \mathbb R)) \setminus U$ such that $M(\alpha)$ is not $\cal H$-minimal. Then  there are an infinite sequence of distinct slopes $\alpha_n$ on $\partial M$ which does not subconverge to a $\partial$-slope and strict finite dominations $f_n: M(\alpha_n) \to V_n$, where $V_n$ are closed hyperbolic $3$-manifolds. The sequence $\{\alpha_n\}$ verifies the hypotheses of Theorem \ref{thm:domination}, hence there is a strict domination $f_0: M \to V$, where $V$ is a $1$-cusped, complete, hyperbolic $3$-manifold, contrary to the $\cal H$-minimality of $M$. 

(2) The first assertion follows from the argument in the proof of part (1) while the second follows from Corollary \ref{compactconditions}. 
\qed 

\noindent {\bf Proofs of Theorem \ref{thm:h-minimal} and Corollary \ref{cor:2-bridge}}. Theorem \ref{thm:h-minimal} is the first assertion of Theorem \ref{minimal hyperbolic} while Corollary \ref{cor:integer surgery} follows from second assertion and Corollary \ref{2bridge compact}. 
\qed 

\section{Sets of discrete $PSL_2(\mathbb R)$-characters} \label{psl2rcharacters}

\subsection{Discrete $PSL_2(\mathbb R)$-representations of the fundamental groups of small knot manifolds} \label{realdiscrete} 

In this section we specialize our study to sets of discrete $PSL_2(\mathbb R)$-characters and apply our conclusions to obtain results on $\widetilde{SL_2}$-minimality. This will lead us, for instance, to a proof of Corollary \ref{prop:g-minimalintro} and thus the construction of infinitely many closed minimal manifolds. 

Let $M$ be a small knot manifold and set 
$$D(M; {\mathbb R}) = \{ \chi_\rho \in X_{PSL_2}(M) : \rho \hbox{ is a discrete, non-elementary $PSL_2(\mathbb R)$ representation}\}.$$ 
If $X_0$ is a component of $X_{PSL_2}(M)$ let
$$D(X_0; {\mathbb R}) = \{ \chi_\rho \in X_0 : \rho \hbox{ is a discrete, non-elementary $PSL_2(\mathbb R)$ representation}\}.$$
Thus $D(X_0; \mathbb R) = D(X_0) \cap X_{PSL_2(\mathbb R)}(M)$ and so is closed in $X_0$ (cf. Proposition \ref{closedinX_0}). 

Fix $\rho \in R_{PSL_2}(M)$ such that $\chi_\rho \in D(M; {\mathbb R})$, set $\Delta = \rho(\pi_{1}(M))$, and let ${\cal B} = \mathbb H^2/\Delta$. The underlying surface $|{\cal B}|$ of ${\cal B}$ is orientable and therefore has only cone singularities. 

\begin{lemma} \label{smallfuchsian} $\;$ \\ 
$(1)$ $\Delta$ is either a hyperbolic triangle group or a free product of two finite cyclic groups. \\
$(2)$ $\rho(\pi_1(\partial M)) \cong \mathbb Z/c$ for some $c \geq 0$ and if $c > 0$, $\Delta$ is a hyperbolic triangle group. \\
$(3)$ Suppose that $\chi_{\rho_n} \in D(X_0; \mathbb R), (n \geq 1)$ are distinct and $\rho_n(\pi_1(\partial M)) \cong \mathbb Z / c_n$ with $c_n \geq 1$. Then $\lim_n c_n = \infty$. 

\end{lemma}

\pf (1) If $|{\cal B}|$ is non-compact, then $\Delta = \pi_1({\cal B}) \cong \pi_1(|{\cal B}|) * \mathbb Z/a_1 * \ldots * \mathbb Z/ a_k$ where $a_1, a_2, \ldots , a_k \geq 2$ are the orders of the cone points. On the other hand, we can identify $X_{PSL_2}(\Delta)$ with a closed algebraic subset of $X_{PSL_2}(M)$. Since the latter has complex dimension $1$, either $\pi_1(|{\cal B}|) \cong \{1\}$ and $k \leq 2$ or $\pi_1(|{\cal B}|) \cong \mathbb Z$ and $k = 0$. The latter is impossible since it implies that $\Delta \cong \mathbb Z$. Thus $\Delta$ is a free product of two finite cyclic groups. In this case, if $\alpha \in \hbox{kernel}(\rho|\pi_1(\partial M))$, $X_{PSL_2}(M(\alpha))$ has dimension $1$ and since $M$ is small, $\alpha$ is a boundary slope.  
 
Next suppose that $|{\cal B}|$ is closed. The relation which associates a holonomy representation to a hyperbolic structure determines an embedding of the Teichm\"{u}ller space ${\cal T}({\cal B})$ of ${\cal B}$ in $X_{PSL_2(\mathbb R)}(\pi_1({\cal B})) \subset  X_{PSL_2}(\pi_1({\cal B})) \subset X_{PSL_2}(M)$. Thus ${\cal T}({\cal B})$ has real dimension at most $1$. But this dimension is given by $-3\chi(|{\cal B}|) + 2k$ where $k$ is the number of cone points in ${\cal B}$ (Corollary 13.3.7 \cite{Thu}). Since $|{\cal B}|$ is orientable, the only possibility is for it to be of the form $S^2(a,b,c)$ so that $\Delta$ is a hyperbolic triangle group. 

(2) The first assertion of (2) follows from the elementary observation that an abelian subgroup of $PSL_2(\mathbb R)$ is cyclic. For the second, suppose that $c > 0$ and note that there are infinitely many slopes in $\hbox{kernel}(\rho|\pi_1(\partial M))$. Fix one such slope $\alpha$ and suppose that $\Delta$ is a free product of two finite cyclic groups. There is a principal curve $Y_0 \subset X_{PSL_2}(\Delta) \subset X_{PSL_2}(M(\alpha))$ and so $M(\alpha)$ admits a closed essential surface. Since $M$ is small, $\alpha$ is a boundary slope, and as there are only finitely many such slopes \cite{Hat}, we obtain a contradiction. Thus $\Delta$ is a hyperbolic triangle group. 

(3) Otherwise there is a subsequence $\{j\}$ and $c \geq 1$ such that $c_n = c$ for all $j$. Then for any peripheral class $\gamma$ there are only finitely many possibilities for $f_\gamma(\chi_j)$. Since the $\chi_j$ are distinct this implies that each $f_\gamma$ is constant, which contradicts the smallness of $M$. 
\qed

\begin{lemma} \label{triangleisolated} 
If the image of $\rho \in R_{PSL_2}(M)$ is a discrete hyperbolic triangle group, then $\chi_\rho$ is an isolated point of $D(M; \mathbb R)$. 
\end{lemma}

\pf Suppose that there is a sequence $\{\chi_{\rho_n}\}$ in $D(M; \mathbb R) \setminus \{\chi_\rho\}$ which converges to $\chi_\rho$. By passing to a subsequence and replacing the $\rho_n$ by conjugate representations we may suppose that $\lim_n \rho_n = \rho$ (Lemma 2.1 \cite{CL}) and find homomorphisms $\theta_n: \rho(\pi_1(M)) \to PSL_2(\mathbb C)$ such that $\rho_n = \theta_n \circ \rho$ (Proposition \ref{convergence}). We claim that the $PSL_2(\mathbb C)$ character varieties of triangle groups are finite. Assuming this for the moment, 
by again passing to a subsequence we may find $A_n \in PSL_2(\mathbb C)$ such that $\theta_n = A_n \theta_1 A_n^{-1}$. Then $\rho_n = \theta_n \circ \rho=  A_n (\theta_1 \circ \rho) A_n^{-1}$. Hence $\chi_\rho = \lim_n \chi_{\rho_n} = \lim_n \chi_{\rho_1} = \chi_{\rho_1} \in D(M; \mathbb R) \setminus \{\chi_\rho\}$, which is impossible. Thus $\chi_\rho$ is an isolated point of $D(M; \mathbb R)$.

To see that the character variety of the $(p,q,r)$-triangle group $\Delta(p,q,r) = \langle x,y : x^p = y^q = (xy)^r = 1 \rangle$ is finite, note that there is a natural embedding $X_{PSL_2}(\Delta(p,q,r)) \subset X_{PSL_2}(\mathbb Z/p * \mathbb Z/q)$. Indeed, $X_{PSL_2}(\Delta(p,q,r))$ is contained in the set of points where the regular function $f: X_{PSL_2}(\mathbb Z/p * \mathbb Z/q) \to \mathbb C, \chi_\rho \mapsto \hbox{trace}(\rho(xy))^2$ takes on the value $4 \cos^2(\frac{\pi j}{r})$ for some integer $j$. Now $X_{PSL_2}(\mathbb Z/p * \mathbb Z/q)$ consists of a finite union of curves and isolated points (Example 3.2 \cite{BZ1}) and it is simple to see from the parameterizations given in that example that the restriction of $f$ to any of the curves is non-constant. Thus it takes on the value $4 \cos^2(\frac{\pi j}{r})$ at only finitely many points and therefore $X_{PSL_2}(\Delta(p,q,r))$ is finite. 
\qed 

\begin{lemma} \label{productnotisolated} 
Suppose that the image $\Delta$ of $\rho \in R_{PSL_2}(M)$ is isomorphic to $\mathbb Z / p * \mathbb Z / q$. Then $\chi_\rho$ is an accumulation point of $D(X_0; \mathbb R)$ where $X_0 = \rho^*(Y_0)$ for some principal component $Y_0$ of $X_{PSL_2}(\Delta)$. Further, $D(X_0; \mathbb R)$ is non-compact in $X_0$ and there is a compact subset $K \subset X_0$ such that \\ 
\indent $($a$)$ $\hbox{int}(K)$ contains all characters in $D(X_0; \mathbb R)$ of representations whose images are \\ \indent \hspace{5mm} hyperbolic triangle groups, and \\ 
\indent $($b$)$ $\overline{(X_0 \setminus K)} \cap D(X_0; \mathbb R)$ contains all characters in $D(X_0; \mathbb R)$ of representations whose \\ \indent \hspace{5mm}  images are $\mathbb Z / p * \mathbb Z / q$. 
\end{lemma}

\pf The inclusion $\Delta \to PSL_2(\mathbb C)$ is contained in a unique curve $Y_0 \subset X_{PSL_2}(\Delta)$ (cf. Example 3.2 \cite{BZ1}). Set $X_0 = \rho^*(Y_0)$. The remaining assertions of the lemma are a consequence of the discussion in Remark \ref{pqnoncompact} and Theorem 2.3 of \cite{Kn}.
\qed

\begin{lemma} \label{smallproduct} 
Let $\{\chi_n\} \subset D(X_0; \mathbb R)$ be a sequence of distinct characters of representations $\rho_n$ with image a free product of two finite cyclic groups. Then there are an epimorphism $\pi_1(M) \to \mathbb Z/p * \mathbb Z / q$ $(2 \leq p, q)$ and a principal curve $Y_0 \subset X_{PSL_2}(\mathbb Z/p * \mathbb Z/q)$ which maps bijectively to $X_0$ under the inclusion $X_{PSL_2}(\mathbb Z/p * \mathbb Z/q) \subset X_{PSL_2}(M)$. In particular, $X_0$ is an $\alpha_0$-curve for some slope $\alpha_0$ on $\partial M$ and $\rho_n(\pi_1(M)) \cong \mathbb Z/p * \mathbb Z / q$ for all $n$. 
\end{lemma}

\pf Choose $n \gg 0$ such that $\chi_n$ is a simple point of $X_{PSL_2}(M)$. By hypothesis, the image $\Delta$ of $\rho_n$ is isomorphic to $\mathbb Z/p * \mathbb Z / q$ for some $2 \leq p, q$. There is a principal curve $Y_0 \subset X_{PSL_2}(\Delta)$ containing the inclusion $\Delta \to PSL_2(\mathbb C)$ and since $\chi_n$ is a simple point, its image in $X_{PSL_2}(M)$ is $X_0$. Lemma \ref{smallfuchsian}(2) implies that $X_0$ is an $\alpha_0$-curve for some slope $\alpha_0$. Finally, for each $n$, there is an epimorphism $\mathbb Z / p * \mathbb Z / q \cong \Delta \to \rho_n(\pi_1(M)) \cong \mathbb Z / r * \mathbb Z / s$ for some $r, s \geq 2$. It follows from Example 3.2 of \cite{BZ1} that the induced homomorphisms $\mathbb Z / p, \mathbb Z / q \to \mathbb Z / r * \mathbb Z / s$ are injective. Further, since these images conjugate into one of $\mathbb Z / r,  \mathbb Z / s$ and generate $\mathbb Z / r * \mathbb Z / s$, we must have $\mathbb Z / r * \mathbb Z / s \cong \mathbb Z / p * \mathbb Z / q$. 
\qed

\subsection{Unbounded sequences of discrete $PSL_2(\mathbb R)$-characters} \label{realunbounded} 

Let $M$ be a small knot manifold and $X_0$ a component of $X_{PSL_2}(M)$. Consider a sequence $\{\chi_n\}$ in $D(X_0; \mathbb R)$ which converges to an ideal point $x_0$ of $X_0$. The following lemma is a consequence of Lemmas \ref{projectiveconvergence}(1) and \ref{smallfuchsian}. 

\begin{lemma} \label{projectiveconvergencereal} 
Let $M$ be a small knot manifold, $X_0$ a curve component of $X_{PSL_2}(M)$, and $\{\chi_n\} \subset D(X_0; \mathbb R)$ a sequence which converges to an ideal point $x_0$ of $X_0$. Fix $\rho_n \in R_{X_0}$ such that $\chi_n = \chi_{\rho_n}$ and let $\alpha_0$ be  the $\partial$-slope associated to $x_0$.  For $n \gg 0$, $\hbox{kernel}(\rho_n|\pi_1(\partial M)) \cong \mathbb Z$ and $\rho_n(\pi_1(\partial M)) \cong \mathbb Z$ where the $\mathbb Z$ factor is generated by a loxodromic. 
\qed
\end{lemma}

There is no subgroup of $PSL_2(\mathbb R)$ isomorphic to the tetrahedral, octahedral, or icosahedral group. Thus the next result follows directly from Theorem \ref{elementaryimage}. 

\begin{thm} \label{convidealreal}
Let $M$ be a small knot manifold, $X_0$ a component of $X_{PSL_2}(M)$, and $\{\chi_n\} \subset D(X_0; \mathbb R)$ a sequence which converges to an ideal point $x_0$ of $X_0$. If $S_0$ is a component of an essential surface associated to $x_0$, then for $n \gg 0$, the image of $X_0$ in $X_{PSL_2}(S_0)$ is contained in $X_{{\cal N}}(S_0)$.  
\qed 
\end{thm}

\subsection{Convergent sequences of discrete $PSL_2(\mathbb R)$-characters} \label{realconvergent} 

Let  $M$ be a small knot manifold and $X_0$ a non-trivial component of $X_{PSL_2}(M)$. We are interested in the accumulation points of $D(X_0; \mathbb R)$ in $X_0$. 

\begin{thm} \label{convordreal}
Let  $M$ be a small knot manifold, $X_0$ a non-trivial component of $X_{PSL_2}(M)$, and $\{\chi_n\} \subset D(X_0; \mathbb R)$ a sequence of distinct characters which converge to some $\chi_{\rho_0} \in  X_0$.  Then \\ 
$(1)$ $\rho_0(\pi_1(M)) = \Delta_0$ is discrete, non-elementary, and isomorphic to $\mathbb Z / p * \mathbb Z /q$ for some integers $2 \leq p, q$; \\ 
$(2)$ there is a principal component $Y_0 \subset X_{PSL_2}( \Delta_0)$ such that $X_0 = \rho_0^*(Y_0)$; \\ 
$(3)$ there is a unique slope $\alpha_0$ on $\partial M$ such that $\rho_0(\alpha_0) = \pm I$ and $X_0$ is an $\alpha_0$-curve; \\
$(4)$ if $\rho_n(\pi_1(\partial M))$ is finite for infinitely many $n$, $\rho_0(\pi_1(\partial M)) \cong \mathbb Z$  is generated by a parabolic. 
\end{thm}

\pf Fix $\rho_n \in R_{X_0}$ whose character is $\chi_n$. After replacing the $\rho_n$ by conjugate representations (over $PSL_2(\mathbb R)$) and passing to a subsequence, we may suppose that $\lim \rho_n = \rho_0$. Let $\Delta_n$ the image of $\rho_n$ ($n \geq 0$) and ${\cal B}_n = \mathbb H^2/ \Delta_n$ ($n \geq 1$). Lemma \ref{smallfuchsian} shows that for $n \geq 1$, $\Delta_n$ is either a free product of two finite cyclic groups or a hyperbolic triangle group. Thus the topological orbifold type of ${\cal B}_n$ is either $\mathbb R^2(p,q)$ or $S^2(p,q,r)$. Since $\{\Delta_n\}$ converges algebraically to $\Delta_0$, Proposition \ref{convergence} implies that $\Delta_0$ is a non-elementary Kleinian group and after passing to a subsequence we may suppose that $\{\Delta_n\}$ converges geometrically to a Fuchsian group $\Delta$ containing $\Delta_0$. Further, there are homomorphisms $\theta_n: \Delta \to \Delta_n$ such that $\rho_n = \theta_n \circ \rho_0$ ($n \geq 1$) and $\lim_n \theta_n$ is the inclusion $\Delta \to PSL_2(\mathbb C)$. 

Assume first that the image of $\rho_n$ is a free product of finite cyclic groups for infinitely many $n$. By Lemma \ref{smallproduct} there are an integer $n \gg 0$, integers $p, q \geq 2$, and a principal component $Z_0$ of $X_{PSL_2}(\Delta_n)$ such that $\Delta_n \cong \mathbb Z / p * \mathbb Z /q$ and $X_0 = \rho_n^*(Z_0) \subset \rho_n^*(X_{PSL_2}(\Delta_n))$. Hence as $\rho_0$ is irreducible, we have $\rho_0 = \psi \circ \rho_n$ for some $\psi \in R_{PSL_2}(\Delta_n)$. It follows that we have surjective homomorphisms $\Delta_n \stackrel{\psi}{\longrightarrow} \Delta_0$ and $\Delta_0 \stackrel{\theta_n}{\longrightarrow} \Delta_n$. Since $\Delta_0$ and $\Delta_n$ are Hopfian, $\theta_n$ is an isomorphism. Thus $Y_0 = \theta_n^*(Z_0)$ is a principal component of $X_{PSL_2}(\Delta_0)$ and $X_0 = \rho_n^*(Z_0) = \rho_0^*(Y_0)$. Lemma \ref{smallfuchsian} shows that the remaining conclusions (3) and (4) of the proposition hold. 

Next assume that $\Delta_n$ is isomorphic to the $(a_n, b_n, c_n)$ triangle group $\Delta(a_n, b_n, c_n)$ where $2 \leq a_n \leq b_n \leq c_n$. Then ${\cal B}_n = S^2(a_n, b_n, c_n)$. We know that $\rho_0(\pi_1(\partial M)) \cong \mathbb Z / d$ for some $d \geq 0$. If $d > 0$, then $\rho_n(\pi_1(\partial M)) = \theta_n(\rho_0(\pi_1(\partial M)))$ is a quotient of  the finite group $\mathbb Z / d$ for all $n$, which contradicts Lemma \ref{smallfuchsian}. Thus $\rho_0(\pi_1(\partial M)) \cong \mathbb Z$. Let $\alpha_0$ be the unique slope such that $\rho_0(\alpha_0) = \pm I$. 

\begin{claim} \label{limitpq} 
The sequence $\{c_n\}$ tends to infinity and after passing to a subsequence we can find integers $2 \leq p \leq q$ such that $a_n = p, b_n = q$ for all $n$. Further, $\Delta \cong \mathbb Z/p * \mathbb Z/q$ and $\Delta_0$ has index at most $2$ in $\Delta$. If it has index $2$, then $\Delta \cong \mathbb Z/2 * \mathbb Z/q, \Delta_0 \cong \mathbb Z/q * \mathbb Z/q$, and $c_n$ is odd.  
\end{claim}

\pf  If $\{c_n\}$ is a bounded sequence, then so are $\{a_n\}, \{b_n\}$ and so after passing to a subsequence we can suppose that they are constants $a, b, c$. We know that $\mathbb H^2/ \Delta = \lim_n \mathbb H^2/ \Delta_n = \lim_n {\cal B}_n= S^2(a,b,c)$. Thus $\Delta \cong \Delta(a,b,c)$ and as this group is Hopfian, it follows that $\theta_n: \Delta \to \Delta_n$ is an isomorphism for all $n$. Since the groups $\Delta_n$ are conjugate in $PSL_2(\mathbb R)$ and the outer automorphism group of  $\Delta(a,b,c)$ is finite, it follows that there are only finitely many conjugacy classes among the representations $\rho_n = \theta_n \circ \rho_0$, which contradicts our assumptions. Thus after passing to a subsequence we may suppose that $\lim_n c_n = \infty$. 

If $\{a_n\}$ is not bounded, then up to passing to a subsequence we may suppose that $\lim_n a_n = \infty$. It follows that $\lim_n b_n = \infty$ and therefore $\mathbb H^2/ \Delta = \lim_n {\cal B}_n$ is a thrice-punctured sphere. But then $\Delta$ is a free group on two generators, and therefore the non-abelian group $\Delta_0$ is free. Thus the dimension of $X_{PSL_2}(\Delta_0)$ is at least $3$. But this is impossible as $\rho_0^*: X_{PSL_2}(\Delta_0) \to X_{PSL_2}(M)$ is injective. Thus $\{a_n\}$ is bounded so that after passing to a subsequence we may suppose that $a_n = p \geq 2$ for all $n$. 

A similar argument shows that if $\{b_n\}$ is unbounded, then $\Delta_0$ is a non-abelian subgroup of $\Delta \cong \mathbb Z/p * \mathbb Z$. Then $\Delta_0$ is a free product of at least two cyclic groups, each of which is either free or has order dividing $p$. If there are either three such factors or two with one of them free, a contradiction is obtained as in the previous paragraph. On the other hand if $\Delta_0 \cong \mathbb Z/r * \mathbb Z/s$ where $r$ and $s$ divide $p$, then $\Delta(p, b_n, c_n) = \theta_n(\Delta_0)$ is generated by two elements of order dividing $p$. Knapp \cite{Kn} studied when two elliptics can generate a triangle group and determined necessary and sufficient conditions on their orders and the coefficients of the triangle group for this to occur. It follows from Theorem 2.3 of \cite{Kn} (and its proof) that if $\Delta(p, b_n, c_n)$ is generated by elements of bounded order, then $\{b_n\}$ is a bounded sequence, contrary to our assumptions.  Thus by passing to a subsequence we may suppose that $b_n = q \geq p$ for all $n$. 

The work above shows that $\mathbb H^2/ \Delta = \lim_n S^2(p,q, c_n) = \mathbb R^2(p,q)$ so that $\Delta \cong \mathbb Z / p * \mathbb Z / q$. Hence $\Delta_0 \subset \Delta$ is a free product of cyclic groups. and the smallness of $M$ implies that it must be of the form $\mathbb Z/r * \mathbb Z/s$ where each of $r, s$ divides at least one of $p, q$. It follows that $\Delta(p, q, c_n)$ is generated by two elements whose orders divide $r, s$ respectively. Given our constraints on $c_n$ and $r, s$, Theorem 2.3 of \cite{Kn} shows that the conclusion of the claim holds. 
\qed (Claim \ref{limitpq})

There is a principal component $Y_0$ of $X_{PSL_2}(\Delta_0)$ which contains the character of the inclusion $\Delta_0 \to PSL_2(\mathbb C)$. Since $\lim_n \theta_n$ is this inclusion and the algebraic components of $X_{PSL_2}(\Delta_0)$ are topological components (see Example 3.2 of \cite{BZ1}), if $n \gg 0$, $\chi_{\theta_n} \in Y_0$. On the other hand, $\chi_n$ is a simple point of $X_{PSL_2}(M)$ for $n \gg 0$.  Since $\chi_{\rho_0} = \lim_n \chi_n = \lim_n \rho_0^*(\chi_{\theta_n})$ it follows that $X_0 = \rho_0^*(Y_0)$. This proves (1) and (2) while (3) is a consequence of the (1), (2),  and Lemma \ref{smallfuchsian}. Finally, to prove (4), note that if $\alpha_1 \ne \alpha_0$ is a slope, then $|\hbox{trace}(\rho_0)(\alpha_1)| = \lim_n |\hbox{trace}(\rho_n)(\alpha_1)| \leq 2$. On the other hand if $\rho_0(\alpha_1)$ is elliptic, then $\rho_n(\alpha_1)$ is elliptic of the same order for $n \gg 0$. This contradicts Lemma \ref{smallfuchsian}. Thus $\rho_0(\alpha_1)$ is parabolic. 
\qed 

\begin{cor} \label{drintersectionfinite}
Suppose that $M$ is a small knot manifold. \\ 
$(1)$ If $X_0$ is a norm curve component of $X_{PSL_2}(M)$, then the intersection of $D(X_0; \mathbb R)$ with any compact subset of $X_0$ is finite. \\ 
$(2)$ If $\pi_1(M)$ does not surject onto a free product of non-trivial cyclic groups. Then the intersection of $D(M; \mathbb R)$ with any compact subset of $X_{PSL_2}(M)$ is finite. 
\qed
\end{cor}

\begin{exa} \label{propertyq} 
{\rm The character variety of a knot manifold $M$ whose fundamental group admits a discrete epimorphism onto a free product of finite cyclic groups contains an $\alpha_0$-curve for some slope $\alpha_0$. Hence Example \ref{twistpretzel} gives many examples for which this does not occur. In particular Corollaries \ref{principalcompact} and \ref{dmrfinite} show that if $M$ is the exterior of a hyperbolic twist knot or a $(-2,3,n)$ pretzel knot with $n \not \equiv 0$ (mod $3$), then $D(M; \mathbb R)$ is finite. Gonz\`alez-Acu$\tilde{n}$a and Ramirez \cite{GR1}, \cite{GR2} have studied the problem of when the fundamental group of the exterior $M$ of a knot in the $3$-sphere admits an epimorphism onto a free product $\mathbb Z / p * \mathbb Z /q$ for some integers $p, q \geq 2$. It is simple to see that in this case $p, q$ are relatively prime. Hartley and Murasugi showed \cite{HM} that  the epimorphism factors through a homomorphism $\pi_1(M) \to \pi_1(M_{p,q})$ whose image is normal with cokernel finite cyclict. This implies that the Alexander polynomial of $M$ is divisible by that of $M_{p,q}$. These conclusions hold more generally for manifolds $M$ with $H_1(M) \cong \mathbb Z$ (cf. the proof of Theorem \ref{pqdomination}). Gonz\`alez-Acuna and Ramirez \cite{GR1} have given an algorithm which determines which two-bridge knot exteriors have fundamental groups which admit such a representation. This work easily shows that the fundamental group of the exterior of the $\frac{p}{q}$ two-bridge knot, $p$ prime, admits no such representation. }
\end{exa}

\begin{cor} \label{dmrfinite}
$\;$  \\ 
$(1)$ If $M$ is a small knot manifold and $X_0$ is a non-trivial component of $X_{PSL_2}(M)$ such that for each connected, essential surface $S_0$ in $M$ there is a character $\chi \in X_0$ such that $\chi|\pi_1(S_0)$ is strictly irreducible, then $D(X_0; \mathbb R)$ is finite. \\
$(2)$ Let $M$ and $N$ be small hyperbolic knot manifolds and suppose that $\varphi: \pi_1(M) \to \pi_1(N)$ is a virtual epimorphism. Then if $Y_0$ is a principal component of $X_{PSL_2}(N)$ and $X_0 = \varphi^*(Y_0)$, then $D(X_0; \mathbb R)$ is finite. 
\end{cor}

\pf (1) By Theorem \ref{convidealreal} we deduce that $D(X_0; \mathbb R)$ is compact. If it has an accumulation point then Theorem \ref{convordreal} implies that there is a surjection $\rho: \pi_1(M) \to \mathbb Z / p * \mathbb Z / q$ and a principal component $Y_0$ of $X_{PSL_2}(\mathbb Z / p * \mathbb Z / q)$ such that $X_0 = \rho^*(Y_0)$. But then Lemma \ref{productnotisolated} shows that $D(X_0; \mathbb R)$ is not compact. Thus $D(X_0; \mathbb R)$ is finite. 

(2) Corollary \ref{principalcompact} implies that $D(X_0; \mathbb R)$ is compact. If it has an accumulation point then Theorem \ref{convordreal} implies that for each irreducible $\chi_\rho \in X_0$, $\rho(\pi_1(M))$ is generated by two torsion elements (cf. Remark \ref{pqnoncompact}). But this is clearly not the case for the image by $\rho^*$ of the discrete faithful character of $\pi_1(N)$. Thus $D(X_0; \mathbb R)$ is finite. 
\qed 

\subsection{Discrete $PSL_2(\mathbb R)$-representations and domination} \label{psl2domination}

\begin{thm} \label{pqdomination} 
Let $M$ be a knot manifold with $H_1(M) \cong \mathbb Z$ and suppose that there is a homomorphism $\rho_0 \in R_{PSL_2}(M)$ with discrete, non-elementary image $\Delta_0 \cong \mathbb Z / p * \mathbb Z / q$. Suppose further that $\rho_0(\lambda_M)$ is parabolic for any longitudinal class $\lambda_M \in \pi_1(\partial M)$. Then there are a Seifert fibred manifold $N$ whose interior has base orbifold $\mathbb H^2 / \Delta_0 \cong \mathbb R^2(p,q)$ and a domination $f: (M, \partial M) \to (N, \partial N)$ such that the composition $\pi_1(M) \stackrel{f_\#}{\longrightarrow} \pi_1(N) \to \Delta_0$ is conjugate to $\rho_0$. 
\end{thm}

\pf Consider the central  extension 
$$1 \to K \to \mbox{Isom}_0(\widetilde{SL_2}) \stackrel{\psi}{\longrightarrow} PSL_2(\mathbb R) \to 1$$ 
where $\mbox{Isom}_0(\widetilde{SL_2})$ is the component of the identity in $\mbox{Isom}(\widetilde{SL_2})$ and $K \cong \mathbb R$ (cf. pp. 464-465 of \cite{Sc}). It is simple to see that for each torsion element $x \in PSL_2(\mathbb R)$, there is a unique torsion element $\tilde A \in \psi^{-1}(A) \subset \mbox{Isom}_0(\widetilde{SL_2})$. Thus $\rho_0$ lifts to a representation $\tilde \rho_0: \pi_1(M) \to \mbox{Isom}_0(\widetilde{SL_2})$ whose image is isomorphic to $\Delta_0$. Fix a non-zero homomorphism $\phi: \pi_1(M) \to K$ and note that 
$$\tilde \rho: \pi_1(M) \to \mbox{Isom}_0(\widetilde{SL_2}), \gamma \mapsto \phi(\gamma) \tilde \rho_0(\gamma)$$ 
is another homomorphism which lifts $\rho_0$. Set $\tilde \Delta_\phi = \tilde \rho(\pi_1(M))$.

\begin{claim} \label{disretelift} 
$\tilde \Delta_\phi$ is discrete, torsion free, and is the fundamental group of a Seifert manifold $N$ with base orbifold $D^2(p,q)$.  
\end{claim}

\pf Since $\Delta_0$ is discrete in $PSL_2(\mathbb R)$, $\tilde \Delta_\phi$ is discrete in $\mbox{Isom}_0(\widetilde{SL_2})$ if and only if it intersects the central subgroup $K$ of $\mbox{Isom}_0(\widetilde{SL_2})$ in a discrete subgroup. This intersection is precisely $\tilde \rho(\mbox{kernel}(\rho_0)) = \phi(\mbox{kernel}(\rho_0)) \subset \phi(\pi_1(M)) \subset K$. The latter group is isomorphic to $\mathbb Z$ by construction, and so is discrete. Thus $\tilde \Delta_\phi$ is discrete. 

Suppose that $\gamma \in \pi_1(M)$ and $\tilde \rho(\gamma)^n = 1$ for some positive $n$. Then up to conjugation, $\rho_0(\gamma)^n = \pm I$ is also torsion and therefore $\tilde \rho_0(\gamma)^n = 1$ as well. But then $1 =  \tilde \rho(\gamma)^n = \phi(\gamma)^n \tilde \rho_0(\gamma)^n = \phi(\gamma)^n$. Since $K$ is torsion free we conclude that $\gamma \in \hbox{kernel}(\phi)$, and since $H_1(M) \cong \mathbb Z$ and $\phi \ne 0$, $\hbox{kernel}(\phi) = [\pi_1(M), \pi_1(M)]$. Hence $\gamma \in  [\pi_1(M), \pi_1(M)]$ and therefore the image of $\rho_0(\gamma)$ in $H_1(\Delta_0)$ is zero.  But $\Delta_0 \cong \mathbb Z / p * \mathbb Z / q$ so that $\Delta_0 \to H_1(\Delta_0)$ is injective on torsion elements. Thus $\rho_0(\gamma) = 1$ and therefore $\tilde \rho(\gamma) = \phi(\gamma) \tilde \rho_0(\gamma) = 1$. This proves that $\tilde \Delta_\phi$ is torsion free. 

The conclusions of the two previous paragraphs imply that $\tilde \Delta_\phi$ acts freely and properly discontinuously on $\widetilde{SL_2}$. Let $W = \widetilde{SL_2}/\tilde \Delta_0^\phi$ be  the quotient manifold. Now $\tilde \Delta_\phi \cap K \ne \{0\}$ as otherwise $\psi | \tilde \Delta_\phi \to \Delta_0$ would be an isomorphism, which contradicts the result of the last paragraph. Thus $\tilde \Delta_\phi \cap K \cong \mathbb Z$ and so $K/(\tilde \Delta_\phi \cap K) \cong S^1$. On the other hand, $\mathbb H^2/ \Delta_0 \cong \mathbb R^2(p,q)$. Thus there is an orbifold bundle $S^1 \to W \to \mathbb R^2(p,q)$ so that $W$ admits a compactification $N$ with boundary a torus. Further, $N$ admits a Seifert fibering with base orbifold $D^2(p,q)$. This completes the proof of the claim. 
\qed (Claim \ref{disretelift})

To complete the proof of the proposition we must show that there is a domination $M \geq N$. To that end, fix a map $f: M \to N$ which realizes $\tilde \rho: \pi_1(M) \to \tilde \Delta_\phi = \pi_1(N)$. We must show that $f_\#|\pi_1(\partial M)$ is injective and has image contained in a peripheral subgroup of $\pi_1(N)$. 

By hypothesis $\rho_0(\lambda_M)$ is parabolic. In particular it has infinite order and is distinct from a primitive element $\alpha_0 \in \hbox{kernel}(\rho_0|\pi_1(\partial M))$ (cf. Theorem \ref{convordreal}(3)). It follows that $1 \ne \phi(\alpha_0) \in K$. It is easy to see that the restriction of $\tilde \rho$ to $\langle \lambda_M, \alpha_0 \rangle \cong \mathbb Z^2$ is injective and since $\tilde \Delta_\phi$ is torsion free, the same holds for its restriction to $\pi_1(\partial M)$. 

Finally, to show that $\tilde \rho(\pi_1(\partial M))$ is peripheral, it suffices to see that $\rho_0(\pi_1(\partial M))$ is peripheral in $\Delta_0 = \pi_1(D^2(p,q))$. But this is clear since it is a parabolic subgroup of $\Delta_0$. This completes the proof. 
\qed 

\begin{rem}
{\rm The condition that $H_1(M) \cong \mathbb Z$ was used to guarantee that $\tilde \Delta$ is torsion free. Without this condition we can still construct a proper non-zero degree map from $M$ to a $3$-dimensional Seifert orbifold, but the underlying space of the orbifold might be $S^1 \times D^2$. }
\end{rem}

\begin{cor} \label{seifertdomination}
Let  $M$ be a small hyperbolic knot manifold with $H_1(M) \cong \mathbb Z$, $X_0$ a non-trivial component of $X_{PSL_2}(M)$, and $\{\chi_{\rho_n}\} \subset D(X_0; \mathbb R)$ a sequence of distinct characters which converges to $\chi_{\rho_0} \in X_0$. Suppose further that for each $n$, $\rho_n(\pi_1(\partial M))$ is finite. Then $\rho_0$ has discrete, non-elementary image isomorphic to a free product of two finite cyclic groups and $\rho_0(\pi_1(\partial M))$ is parabolic. If $\rho_n(\lambda_M) \ne \pm I$ for infinitely many $n$, there is a strict domination $M \geq N$ for some Seifert manifold $N$ with incompressible boundary. 
\end{cor}

\pf  Suppose that $\lim_n \chi_{\rho_n} = \chi_{\rho_0}$ and set $\Delta_0 = \rho_0(\pi_1(M))$. By Theorem \ref{convordreal}, $\Delta_0 \subset PSL_2(\mathbb R)$ is discrete and isomorphic to $\mathbb Z / p * \mathbb Z / q$ for some integers $2 \leq p, q$. Our hypotheses imply that $\rho_0(\lambda_M)$ is parabolic (cf. Proposition \ref{elemordiscr}(2) and Theorem \ref{convordreal}(4)). Thus Theorem \ref{pqdomination} implies the desired conclusion. 
\qed 

\begin{cor} \label{cor:s-domination} 
Let  $M$ will be a small hyperbolic knot manifold, $\{\alpha_n\}$ a sequence of distinct slopes on $\partial M$, and $\{\chi_n\} \subset D(M; \mathbb R)$ a sequence of characters of representations $\rho_n$ such that $\rho_n(\alpha_n) = \pm I$ for all $n$. If there are infinitely many distinct $\chi_n$ and the sequence $\{\chi_n\}$ subconverges to a character $\chi_{\rho_0}$, then \\
$(1)$ the image of $\rho_0$ is isomorphic to a discrete, non-elementary free product of two finite cyclic groups; \\
$(2)$ $\rho_n(\pi_1(\partial M))$ is finite for infinitely many $n$ and $\rho_0(\pi_1(\partial M))$ is parabolic. \\
$(3)$ if $H_1(M) \cong \mathbb Z$ and $\rho_0(\lambda_M) \ne \pm I$, $M$ strictly dominates a Seifert manifold with incompressible boundary. 
\end{cor} 

\pf After passing to a subsequence we can assume that the $\chi_n$ are distinct. Part (3) of Theorem \ref{convordreal} shows that there is a unique slope $\alpha_0$ on $\partial M$ such that $\rho_0(\alpha_0) = \pm I$. Then for $n \gg 0$ we have $\rho_n(\alpha_0) = \pm I$ (Proposition \ref{elemordiscr}). Since the $\alpha_n$ are distinct, this implies that  for large $n$, $\rho_n(\pi_1(\partial M))$ is a finite cyclic group of order dividing $\Delta(\alpha_0, \alpha_n)$. Corollary \ref{seifertdomination} then yields a strict domination $f: M \to N$ where $N$ is a Seifert manifold with incompressible boundary. 
\qed 

\section{Minimal Dehn fillings}  \label{gminimal}

In the section we use the results of the paper to construct various infinite families of minimal closed $3$-manifold. 

\begin{lemma} \label{finiteredhaken} $\;$ \\ 
$(1)$ If $M$ is a small knot manifold, there are only finitely many slopes $\alpha$ on $\partial M$ such that $M(\alpha)$ is either reducible or Haken. \\ 
$(2)$ A closed, connected, orientable manifold with infinite fundamental group is either reducible, Haken, or admits a geometric structure modelled on $Nil, \mathbb H^3$, or $\widetilde{SL_2}$. 
\end{lemma}

\pf (1) If $M(\alpha)$ contains an essential surface $S$ and we isotope $S$ so as to minimize $|S \cap \partial M|$, then $S_0 := S \cap M$ is an essential surface in $M$. Since $M$ is small, $\partial S_0 \ne \emptyset$ and has slope $\alpha$. Thus $\alpha$ is a boundary slope. By \cite{Hat}, there are at most finitely many such $\alpha$. Thus (1) holds.

(2) By the geometrization theorem of Perelman we see that a closed, connected, orientable manifold $W$ which is irreducible though not Haken admits a  geometric structure. 
If the structure is $Sol$, $W$ is Haken since it is irreducible and contains an essential torus (\cite{Sc}). If it is $\mathbb S^2 \times \mathbb R, \mathbb E^3$ or $\mathbb H^2 \times \mathbb R$, then $W$ amits a Seifert fibre structure with zero Euler number and therefore it is either reducible or Haken (\cite{Sc}). If it is $\mathbb S^3$, $\pi_1(W)$ is finite. This proves (2). 
\qed

\begin{thm} \label{finitesl2domination}
Suppose that $M$ is a small ${\cal H}$-minimal hyperbolic knot manifold which has the following properties:\begin{description} 
\vspace{-.45cm} \item[\hspace{1.5mm} {\rm(a)}]  There is a slope $\alpha_0$ on $\partial M$ such that $M(\alpha_0)$ is ${\cal H}$-minimal.
\vspace{-.35cm}  \item[\hspace{1.5mm} {\rm(b)}]  For each norm curve  $X_0 \subset X_{PSL_2}(M)$ and for each essential surface $S$ associated to an ideal point of $X_0$ there is a character $\chi_\rho \in X_0$ which restricts to a strictly irreducible character on $\pi_1(S)$. 
\vspace{-.35cm} \item[\hspace{1.5mm} {\rm(c)}]  There is no surjective homomorphism from $\pi_1(M)$ onto a Euclidean triangle group.
\vspace{-.35cm} \item[\hspace{1.5mm} {\rm(d)}]  There is no epimorphism $\rho: \pi_1(M) \to \Delta(p,q,r) \subset PSL_2(\mathbb R)$ such that $\rho(\pi_{1} (\partial M))$ is elliptic or trivial.
\end{description}
\vspace{-.45cm}  Then all but finitely many Dehn fillings $M(\alpha)$ yield a minimal manifold.
\end{thm}

\pf By Theorem \ref{minimal hyperbolic}(2) and Lemma \ref{finiteredhaken}, we need only show that $M(\alpha)$ is $Nil$-minimal and $\widetilde{SL_2}$-minimal for all but finitely many slopes $\alpha$ on $\partial M$. 

Suppose that there is a slope $\alpha$ and a domination $f$ from $M(\alpha)$ to a closed $Nil$-manifold $V$ with base orbifold ${\cal B}_n$. By passing to a cover of $V$ we may suppose that $f_\#$ is surjective. We can also suppose that $\alpha$ is not a boundary slope so that $M(\alpha)$ is not Haken. Since ${\cal B}$ is Euclidean, the only possibility is that ${\cal B} \cong S^2(a,b,c)$ for some Euclidean triple $(a,b,c)$. But then we would have an epimorphism $\pi_1(M) \to \pi_1(M(\alpha)) \to \pi_1(V) \to \pi_1(S^2(a,b,c)) \cong \Delta(a,b,c)$, which contradicts (c). Thus $M(\alpha)$ is $Nil$-minimal for all but finitely many $\alpha$. 

Suppose next that there are a sequence of distinct slopes $\alpha_n$ and dominations $f_n$ from $M(\alpha_n)$ to a closed $\widetilde{SL_2}$-manifold $V_n$ with base orbifold ${\cal B}_n$. By passing to a cover of $V_n$, we may suppose that $(f_n)_\#$ is surjective for all $n$. Let $\rho_n$ be the composition $\pi_1(M) \to \pi_1(M(\alpha_n)) \stackrel{(f_n)_\#}{\longrightarrow} \pi_1(V_n) \to \pi_1({\cal B}_n) \subset PSL_2(\mathbb R) \subset PSL_2(\mathbb C)$. By passing to a subsequence we may suppose that $\chi_n \in X_0$ for some non-trivial curve $X_0$. Hypothesis (d) implies that $\hbox{kernel}(\rho_n|\pi_1(\partial M)) = \langle \alpha_n \rangle$ so that there are infinitely many distinct $\chi_n$ and $\rho_n(\pi_1(\partial M))$ is infinite. It follows that $\rho_n(\pi_1(\partial M))$ contains loxodromics. Thus Corollary \ref{normcondition} implies that $X_0$ is a norm curve. But then hypothesis (b) and Corollary \ref{dmrfinite}(1) imply that there are only finitely many $\chi_n$, contrary to the construction. Thus there is no sequence $\{\alpha_n\}$ as above and so $M(\alpha)$ is $\widetilde{SL_2}$-minimal for all but finitely many $\alpha$. 
\qed 

\begin{cor}\label{cor:s-minimal2}
If $M$ is the exterior of a hyperbolic twist knot, then all but finitely many Dehn fillings $M(\alpha)$ yield a minimal manifold. 
\qed
\end{cor}

\pf Hypothesis (a) of Theorem \ref{finitesl2domination} clearly holds for $M$, and since the only non-trivial curve in $X_{PSL_2}(M)$ is a principal curve (\cite{Bu}), hypothesis (b) holds as well (cf. the proof of Corollary \ref{principalcompact}). Finally, hypotheses (c) and (d) are true for $M$ by Proposition \ref{nil} and Proposition \ref{prop:epi}. 
\qed   

Condition (d) of Theorem \ref{finitesl2domination} is difficult to verify in general. Nevertheless, the following results show that we can still construct infinite families of minimal Dehn fillings in quite general situations. First we need to prove an elementary lemma. 

\begin{lemma} \label{infinitelymanyslopes}
Let $\alpha_0, \alpha_1, \ldots, \alpha_n$ be projectively distinct primitive elements of $\mathbb Z^2$ and suppose that $L_1, L_2, \ldots , L_m$ are subgroups of $\mathbb Z^2$ none of which contains $\alpha_0$. For each $i = 1, 2, \ldots , n$, let $U_i$ be an arc neighbourhood of $[\alpha_i] \in \mathbb P(\mathbb R^2)$ and suppose that $\overline{U_i} \cap \overline{U_j} = \emptyset$ for $i \ne j$. Then there are infinitely many primitive $\alpha \in \mathbb Z^2$ such that $\alpha \not \in L_1 \cup \ldots \cup L_m$ and $[\alpha] \not \in U_1 \cup \ldots \cup U_n$. 
\end{lemma}

\pf Since each $L_j$ is contained in a rank $2$ subgroup of $\mathbb Z^2$ in the complement of $\alpha_0$, we can assume, without loss of generality, that each $L_j$ has rank $2$. Let $\beta_0 \in \mathbb Z^2$ be dual to $\alpha_0$ and fix coprime integers $a, b$ such that $\delta_0 := a \alpha_0 + b \beta_0 \ne \alpha_i$, $0 \leq i \leq n$. Set $L_0 = \{\alpha_0 + n \delta_0 : n \in \mathbb Z\}$ and note that from the definition of $U = U_1 \cup \ldots \cup U_n$ and choice of $\delta_0$, there is some $k_0 > 0$ such that if $|k| \geq k_0$, $[\alpha_0 + k \delta_0] \not \in U$. Define $d \geq 1$ to be the index of $L_1 \cap L_2 \cap \ldots \cap L_m$ in $\mathbb Z^2$ and note that for each $k \in \mathbb Z$, the class $\alpha_k = \alpha_0 + dkb \delta_0 \not \in (L_1 \cup L_2 \cup \ldots \cup L_m)$. The proof is completed by observing that $\alpha_k = (1 + abkd) \alpha_0 + b^2kd \beta_0$ is primitive and $[\alpha_k] \not \in U$ for $|k| \geq k_0$.  
\qed

\begin{thm} \label{uinfinite} 
Let $M$ be an ${\cal H}$-minimal, small, hyperbolic knot manifold and suppose that  $H_1(M) \cong \mathbb Z \oplus T$ where \\ 
\indent $(a)$ $H_1(\partial M) \to H_1(M)/T \cong \mathbb Z$ is surjective, and  \\
\indent $(b)$ $\mathbb Z / a \oplus \mathbb Z / b$ is not a quotient of $T$ for $(a,b) = (2,3), (2,4), (3,3)$. \\
Suppose as well that  \\ 
\indent $(c)$ there is no discrete, non-elementary representation $\rho \in R_{PSL_2(\mathbb R)}(M)$ such that \\ \indent \hspace{5mm} $\rho(\pi_1(M))$ is isomorphic to a free product of two non-trivial cyclic groups and $\rho(\pi_1(\partial M))$ \\ \indent \hspace{5mm} is parabolic; \\
\indent $(d)$ there is a slope $\alpha_0$ on $\partial M$ such that $\pi_1(M(\alpha_0))$ admits no homomorphism onto a \\ \indent \hspace{5mm} non-elementary Kleinian group or a Euclidean triangle group. \\ 
If $U \subset \mathbb P(H_1(\partial M; \mathbb R))$ is the union of disjoint closed arc neighbourhoods of the finite set of boundary slopes of $M$, then there are infinitely many slopes $\alpha$ such that $[\alpha] \in \mathbb P(H_1(\partial M; \mathbb R)) \setminus U$ and $M(\alpha)$ is minimal.  
\end{thm}

\pf By Lemma \ref{finiteredhaken} and Theorem \ref{minimal hyperbolic}(1), it suffices to show that there are infinitely many slopes $\alpha$ such that $[\alpha] \in \mathbb P(H_1(\partial M; \mathbb R)) \setminus U$ and $M(\alpha)$ is both $Nil$-minimal and $\widetilde{SL_2}$-minimal. As we argued in the proof of Theorem \ref{finitesl2domination}, if $\alpha$ is not a boundary slope and there is a domination $M(\alpha) \to V$ where $V$ is a $Nil$ or $\widetilde{SL_2}$ manifold, there is an epimorphism $\rho: \pi_1(M) \to \Delta(a,b,c)$ which can suppose lies in $D(M; \mathbb R)$ if $(a,b,c)$ is a hyperbolic triple. 

Suppose first of all that $\rho: \pi_1(M) \to \Delta(a,b,c)$ is surjective and $(a,b,c)$ is a Euclidean triple with $a \leq b \leq c$. There is an epimorphism $\mathbb Z \oplus T = H_1(M) \to H_1(\Delta(a,b,c)) \cong \mathbb Z / a \oplus \mathbb Z / b$ where $(a,b)$ is one of the pairs $(2,3), (2,4), (3,3)$. Hence our hypotheses imply that there is some $\gamma \in \pi_1(\partial M)$ which is sent to a non-zero element of $H_1(\Delta(a,b,c))$ under the composition $\pi_1(M) \stackrel{\rho}{\longrightarrow} \Delta(a,b,c)) \to  H_1(\Delta(a,b,c))$. It is a simple exercise to then show that $\rho(\gamma)$ has non-trivial finite order in $\Delta(a,b,c)$. (For instance, use the fact that $\Delta(a,b,c)$ can be considered a subgroup of the upper-triangular matrices in $PSL_2(\mathbb C)$.) Since an abelian subgroup of an infinite triangle group is cyclic, $\rho(\pi_1(\partial M)) \cong \mathbb Z / d$ where $d \in \{2,3,4, 6\}$. Thus there are only finitely many possibilities for $\hbox{kernel}(\rho|\pi_1(\partial M))$, say $L_1, \ldots , L_k$. By hypothesis none of them contain $\alpha_0$. Further, if $\alpha \not \in L_1 \cup \ldots \cup L_k$, then $\pi_1(M(\alpha))$ admits no homomorphism onto a Euclidean triangle group. 

Next set $D_U(M; \mathbb R) := \{\chi_\rho \in D(M; \mathbb R) : \rho(\alpha) = \pm I \hbox{ for some slope } \alpha \hbox{ such that } [\alpha] \not \in U\}$ and suppose it is infinite. If $\alpha$ is a slope such that $[\alpha] \not \in U$, then $\alpha$ is not a boundary slope and so there are only finitely many $\chi_\rho \in D_U(M; \mathbb R)$ such that $\rho(\alpha) = \pm I$ (Corollary \ref{smallcharactervariety}). Hence we can find a sequence of distinct slopes $\alpha_n$, a sequence of distinct characters $\chi_{\rho_n} \in D_U(M; \mathbb R)$, and a component $X_0$ of $X_{PSL_2}(M)$ such that $[\alpha_n] \not \in U$,  $\rho_n(\alpha_n) = \pm I$, and $\chi_{\rho_n} \in X_0$. Lemma \ref{projectiveconvergence} shows that $\{\chi_{\rho_n}\}$ does not accumulate to an ideal point of $X_0$. Thus we can suppose that it converges to some $\chi_{\rho_0} \in X_0$. Theorem \ref{convordreal} implies that $\rho_0(\pi_1(M))$ is a free product of two finite cyclic groups and $\rho_0(\pi_1(\partial M))$ is parabolic. But this contradicts our hypotheses. Thus $D_U(M; \mathbb R)$ is finite, say $D_U(M; \mathbb R) = \{\chi_{\rho_1}, \chi_{\rho_2} , \ldots , \chi_{\rho_l}\}$. Set $L_j' = \hbox{kernel}(\rho_j|\pi_1(\partial M))$ ($1 \leq j \leq l$). Then $\alpha_0 \not \in (L'_1 \cup L'_2 \cup \ldots \cup L'_m)$ and if $\alpha \not \in (L'_1 \cup L'_2 \cup \ldots \cup L'_m)$ is a slope such that $[\alpha] \not \in U$, $\pi_1(M(\alpha))$ admits no homomorphism onto a hyperbolic triangle group. The proof is completed by applying Lemma \ref{infinitelymanyslopes} to the subgroups $L_1, \ldots , L_k, L'_1, \ldots , L'_l$. 
\qed  

\begin{exa} \label{lensspace} 
{\rm The theorem applies to the exterior of many knots in lens spaces. For instance, it follows from work of Indurskis \cite{In} that if $M_m$ is the manifold obtained by $m$-Dehn filling on one component of the right-hand Whitehead link, then $X_{PSL_2}(M_m)$ has exactly one non-trivial component and is therefore minimal. For $|m| > 4$, $M_m$ is hyperbolic. If $\mu$ is the slope on $\partial M_m$ corresponding to  a meridian of the Whitehead link, $M_m(\mu) \cong L(m, 1)$. Since $H_1(M_m) \cong \mathbb Z \oplus \mathbb Z / m$ and $H_1(\partial M_m) \to H_1(M_m) / T_1(M_m)$ is onto, the hypotheses of Theorem \ref{uinfinite} are satisfied as long as $m \not \equiv 0$ {\rm (mod $6$)}. For such $m$, $M_m(\alpha)$ is minimal for infinitely many slopes $\alpha$.}
\end{exa}

\begin{thm} \label{minginfinite}
Suppose that $M$ is a minimal small hyperbolic knot manifold such that $H_1(M)$ $\cong \mathbb Z$ and that \\ 
\indent $(a)$ there is no homomorphism $\rho: \pi_1(M(\lambda_M)) \to PSL_2(\mathbb R)$ such that $\rho(\pi_1(M(\lambda_M)))$ is \\ \indent \hspace{5mm} a free product of two non-trivial cyclic groups and $\rho(\pi_1(\partial M))$ parabolic. \\ 
\indent $(b)$ there is a slope $\alpha_0$ on $\partial M$ such that $\pi_1(M(\alpha_0))$ admits no homomorphism onto a \\ \indent \hspace{5mm} non-elementary Kleinian group or a Euclidean triangle group. \\ 
Then there are infinitely many slopes $\alpha$ on $\partial M$ such that $M(\alpha)$ is minimal.  
\end{thm} 

\pf The proof is similar to that of Theorem \ref{uinfinite}. As before it suffices to show that there are subgroups $L_1, L_2, \ldots , L_m$ of $H_1(\partial M)$, none of which contain $\alpha_0$, such that if $\alpha \not \in  L_1 \cup \ldots \cup L_m$ is a slope, though not a boundary slope, then $\pi_1(M(\alpha))$ admits no surjective homomorphism onto an infinite triangle group. Since $H_1(M) \cong \mathbb Z$, the homological conditions (a) and (b) from the statement of Theorem \ref{uinfinite} hold and so there are subgroups $L_1, \ldots , L_k$ of $H_1(\partial M)$, none of which contain $\alpha_0$, such that if $\alpha \not \in L_1 \cup \ldots \cup L_k$, then $\pi_1(M(\alpha))$ admits no homomorphism onto a Euclidean triangle group. 

To derive a similar conclusion for hyperbolic triangle groups, the proof of Theorem \ref{uinfinite} shows that it suffices to fix a disjoint union $U \subset \mathbb P(H_1(\partial M; \mathbb R))$ of closed arc neighbourhoods of the finite set of boundary slopes of $M$ and prove that $D_U(M; \mathbb R) := \{\chi_\rho \in D(M; \mathbb R) : \rho(\alpha) = \pm I \hbox{ for some slope } \alpha \hbox{ such that } [\alpha] \not \in U\}$ is finite. Suppose otherwise and note that as in the proof of Theorem \ref{uinfinite}, we can find a representation $\rho_0: \pi_1(M) \to PSL_2(\mathbb R)$ with discrete image isomorphic to a free product of non-trivial cyclic groups such that $\rho_0(\pi_1(\partial M))$ is parabolic. Hypothesis (a) implies that $\rho_0(\lambda_M) \ne \pm I$ and so Theorem \ref{pqdomination} implies that $M$ strictly dominates some Seifert manifold $N$ with incompressible boundary. This contradicts the minimality of $M$. Thus $D_U(M; \mathbb R)$ is finite and the proof proceeds as in that of Theorem \ref{uinfinite}. 
\qed

\begin{cor} \label{s3ginfinite}
Let $M$ be a minimal, small, hyperbolic $3$-manifold which is the exterior of a knot $K$ in the $3$-sphere. If there is no homomorphism $\rho: \pi_1(M(\lambda_M)) \to PSL_2(\mathbb R)$ such that $\rho(\pi_1(M(\lambda_M)))$ is a free product of two non-trivial cyclic groups and $\rho(\pi_1(\partial M))$ parabolic, then there are infinitely many slopes $\alpha$ on $\partial M$ such that $M(\alpha)$ is minimal. 
\qed 
\end{cor}

\begin{exa} \label{s3ginfiniteegs}
{\rm If the Alexander polynomial of a knot $K \subset S^3$ with exterior $M$ is not divisible by the Alexander polynomial of a non-trivial torus knot, there is no homomorphism of $\pi_1(M)$ onto the free product of two non-trivial finite cyclic groups (cf. Remark \ref{propertyq}). Thus if its exterior is minimal, small, and hyperbolic, there are infinitely many slopes $\alpha$ on $\partial M$ such that $M(\alpha)$ is minimal. This provides many examples. For others, let $K$ be the $(-2, 3, n)$ pretzel where $n \not \equiv 0$ (mod $3$). We noted in Example \ref{twistpretzel} that there is a unique non-trivial component of $X_{PSL_2}(M)$ and it is principal and used this to deduce that $M$ is minimal. It also implies that there is no homomorphism $\rho: \pi_1(M(\lambda_M)) \to PSL_2(\mathbb R)$ such that $\rho(\pi_1(M(\lambda_M)))$ is a free product of two non-trivial cyclic groups (Lemma \ref{productnotisolated}). Thus Corollary \ref{s3ginfinite} may be applied to see that there are infinitely many slopes $\alpha$ on $\partial M$ such that $M(\alpha)$ is minimal. As a final example, Riley has shown that if $K$ is a two-bridge knot and $\rho \in R_{PSL_2}(M)$ is irreducible with $\rho(\pi_1(\partial M))$ parabolic, then $\rho(\lambda_M) \ne \pm I$ (Lemma 1 \cite{Ri}). In particular, if $M$ is the exterior of a $\frac{p}{q}$ two-bridge knot, it is minimal (Corollary \ref{2bridgemini2}), small, and hyperbolic if $p$ is prime and $q \not \equiv \pm 1$ (mod $p$). Thus the corollary implies that there are infinitely many slopes $\alpha$ on $\partial M$ such that $M(\alpha)$ is minimal.}
\end{exa}

\appendix

\section{On the smoothness of dihedral characters} \label{2-bridge}

One goal of this appendix is to prove that if $\mu$ is a meridinal class of the $p/q$ two-bridge knot, then $d_{M_{p/q}}(\mu) = \frac{p-1}{2}$. In order to do this, we determine a useful criterion for the smoothness of dihedral characters. 

\subsection{A cohomological calculation}\label{sec:cohocalcul}

Let $\Gamma$ be a finitely generated group, $V$ is a complex vector space, and $\theta: \Gamma \to GL(V)$ a homomorphism. We use $b_1(\Gamma; \theta)$ to denote the complex dimension of $H^1(\Gamma; V_\theta)$. For instance if $\rho \in R_{PSL_2}(\Gamma)$, the induced action of $\Gamma$ on $sl_2(\mathbb C)$
given by the composition $\Gamma \stackrel{\rho}{\longrightarrow} PSL_2(\mathbb C)
\stackrel{Ad}{\longrightarrow} Aut(sl_2(\mathbb C))$ gives rise to the cohomology group $H^1(\Gamma; sl_2(\mathbb C)_{Ad\rho})$ whose dimension is $b_1(\Gamma; Ad \rho)$. 

Identify the dihedral group of $2n$ elements $D_n$ with the subgroup of ${\cal N}$ generated by the matrices $\pm \left(
\begin{array}{cc} 0 & 1 \\ -1 & 0 \end{array}
\right)$ and $\pm \left( \begin{array}{cc} \zeta & 0 \\ 0 & \zeta^{-1}
\end{array} \right)$ where $\zeta =
\exp(\frac{2 \pi i}{2n})$. Any subgroup of $PSL_2(\mathbb C)$ abstractly
isomorphic to $D_n$ is conjugate in
$PSL_2(\mathbb C)$ to $D_n$.

For each divisor $d \geq 1$ of $n$ we have surjections $\theta_{n,d}: D_n
\to D_d$ given by
$$\theta_{n,d}(\pm \left(\begin{array}{cc} 0 & 1 \\ -1 & 0 \end{array}
\right)) =
\pm \left( \begin{array}{cc} 0 & 1 \\ -1 & 0 \end{array} \right) \mbox{ and
} \theta_{n,d}(\pm \left(
\begin{array}{cc} \zeta & 0 \\ 0 & \zeta^{-1} \end{array} \right)) = \pm
\left( \begin{array}{cc}
\zeta^{\frac{n}{d}} & 0 \\ 0 & \zeta^{-\frac{n}{d}}\end{array} \right).$$

\begin{lemma} \label{b1Adrho}
Let $\rho: \Gamma \to PSL_2(\mathbb C)$ be a representation whose image is
$D_n$, $n > 1$.
For each divisor $d \geq 1$ of $n$ let $\rho_d$ be the composition of
$\rho$ with $\theta_{n,d}$ and set
$\Gamma_{2d} = ker(\rho_d)$. Then
$$b_1(\Gamma; Ad \rho) = b_1(\Gamma_2) - b_1(\Gamma) + \frac{1}{\phi(n)}
\sum_{d|n}
\mu(\frac{n}{d}) b_1(\Gamma_{2d})$$
where $\phi$ is Euler's $\phi$-function and $\mu$ is the M\"{o}bius function.
\end{lemma}

\pf Consider the real basis
$$e_1 = \left( \begin{array}{cc} 1 & 0 \\ 0 & -1 \end{array} \right),
\;\;\;\;\;
e_2 = \left( \begin{array}{cc} 0 & 1 \\ 0 & 0 \end{array} \right), \;\;\;\;\;
e_3 = \left( \begin{array}{cc} 0 & 0 \\ 1 & 0 \end{array} \right)$$
of $sl_2(\mathbb C)$. Let $\Theta = Ad|: {\cal N} \to Aut(sl_2(\mathbb
C))$. The reader
will verify that the span $\langle e_1 \rangle \cong \mathbb C$ of $e_1$ is
invariant under
$\Theta({\cal N})$ as is $\langle e_2, e_3 \rangle \cong \mathbb C^2$. Thus
$$sl_2(\mathbb C)_{\Theta} = \mathbb C_{\Theta_1} \oplus \mathbb
C^2_{\Theta_2} \eqno{({\rm A}.1.1)}$$
where $\Theta_1: {\cal N} \to  GL_1(\mathbb C)$ is given by
$$\Theta_1(A) = \left\{\begin{array}{rl} 1_{\mathbb C} & \mbox{if } A \in
{\cal D} \\ -1_{\mathbb C} &
\mbox{if } A \in {\cal N} \setminus{\cal D} \end{array} \right.$$
and, in terms of the ordered basis $\{e_2, e_3\}$, $\Theta_2: {\cal N} \to
GL_2(\mathbb C)$ is given by
$$\Theta_2(\pm \left(\begin{array}{cc} u & 0 \\ 0 & u^{-1} \end{array}
\right)) =
\left(\begin{array}{cc} u^2 & 0 \\ 0 & u^{-2, } \end{array} \right),
\Theta_2(\pm \left(\begin{array}{cc} 0 & v
\\ -v^{-1} & 0 \end{array} \right)) = \left(\begin{array}{cc} 0 & -v^2 \\
-v^{-2} & 0 \end{array} \right).$$

Without loss of generality we may suppose that the image of $\rho$ lies in
${\cal N}$.
Then (A.1.1) yields the decomposition $sl_2(\mathbb C)_{Ad\rho} = \mathbb
C_{\theta_1} \oplus \mathbb
C^2_{\theta_2}$ where $\theta_j = \Theta_j \circ \rho$. Hence
$$b_1(\Gamma; Ad\rho) = b_1(\Gamma; \theta_1) +
b_1(\Gamma; \theta_2).$$
The proof of the lemma now follows from Claims (1) and (2) below.

\noindent {\bf Claim 1}. $\;$  $b_1(\Gamma; \theta_1) =
b_1(\Gamma_2) - b_1(\Gamma)$.

\pf The $\Gamma$-module $\mathbb C[\Gamma/\Gamma_2] = \mathbb C[\mathbb
Z/2]_\Gamma$
splits into two $1$-dimensional modules $$\mathbb C[\mathbb Z/2]_\Gamma =
\mathbb C_1 \oplus \mathbb
C_{\theta_1}$$ where $\mathbb C_1$ is the trivial $\Gamma$-module. Then
$H^1(\Gamma_2; \mathbb C) \cong
H^1(\Gamma;\mathbb C[\Gamma/\Gamma_2]) = H^1(\Gamma; \mathbb C[\mathbb
Z/2]_\Gamma) =  H^1(\Gamma; \mathbb C_1)
\oplus H^1(\Gamma; \mathbb C_{\theta_1})$. Thus $b_1(\Gamma_2) =
b_1(\Gamma) + b_1(\Gamma; \theta_1)$.
\qed (of Claim 1)

\noindent {\bf Claim 2}. $\;$ $b_1(\Gamma; \theta_2) =
\frac{1}{\phi(n)} \sum_{d|n}
\mu(\frac{n}{d}) b_1(\Gamma_{2d})$ {\it where $\phi$ is Euler's
$\phi$-function and $\mu$ is
the M\"{o}bius function.}

\pf Fix a divisor $d \geq 1$ of $n$ and observe that the $\Gamma$-module
$\mathbb C[\Gamma/\Gamma_{2d}] \cong \mathbb C[D_d]$ splits as a sum
$$\mathbb C[\Gamma/\Gamma_{2d}] = \mathbb C_1 \oplus \mathbb C_{\theta_1}
\oplus \bigoplus_{r=1}^d \mathbb C^2_{\delta_r} \eqno{({\rm A}.1.2)}$$
where $\delta_r = \Theta \circ \delta_r^0 \circ \rho_d: \Gamma \to
GL_2(\mathbb C)$ with
$\delta_r^0: D_d \to D_d$ given by
$$\delta_r^0(\pm \left(\begin{array}{cc} 0 & 1 \\ -1 & 0 \end{array}
\right)) =
\pm \left( \begin{array}{cc} 0 & 1 \\ -1 & 0 \end{array} \right) \mbox{ and
} \delta_r^0(\pm \left(
\begin{array}{cc} \zeta & 0 \\ 0 & \zeta^{-1} \end{array} \right)) = \pm
\left( \begin{array}{cc}
\zeta^{\frac{rn}{d}} & 0 \\ 0 & \zeta^{-\frac{rn}{d}}\end{array} \right).$$
(see \S 5.3 of \cite{Se} for example). If $r$ divides $d$, then $\delta_r^0 \circ \rho_d =
\rho_{\frac{d}{r}}$. Moreover, if $r_1$
and $r_2$ have the same order in $\mathbb Z/d$, there is a $\Gamma$-module
isomorphism between $\mathbb
C^2_{\delta_{r_1}}$ and $\mathbb C^2_{\delta_{r_2}}$. (For under this
condition there is an automorphism $\psi$
of the group $\theta_{r_1}(\Gamma)$ such that $\delta_{r_2} = \psi \circ
\delta_{r_1}$.) Combining these
observations with (A.1.2) and Claim 1 shows that
$$b_1(\Gamma_{2d}) - b_1(\Gamma_2) =  \sum_{e|d} \phi(e) b_1(\Gamma; \Theta \circ \rho_e).$$
This formula holds for each $d$ which divides $n$, so
the M\"{o}bius inversion formula (see \S 16.4 of \cite{HW} for example) yields
$$b_1(\Gamma; \theta_2) = b_1(\Gamma; \Theta \circ
\rho_n)=
\frac{1}{n} \sum_{d|n} \mu(\frac{n}{d}) (b_1(\Gamma_{2d}) -
b_1(\Gamma_2))  = \frac{1}{\phi(n)}  \sum_{d|n} \mu(\frac{n}{d})
b_1(\Gamma_{2d}),$$
as $n > 1$.
\qed (of Claim 2)

This completes the proof of Lemma \ref{b1Adrho}.
\qed

\subsection{A criterion for the smoothness of dihedral characters of knot groups} \label{criterion}

For a knot $K$ in a $\mathbb Z$-homology $3$-sphere $W$ we use  $\widehat{W}_2(K) \to W$ to denote the $2$-fold cover of $W$ branched over $K$. It is well-known that any irreducible representation of the fundamental group of the exterior of $K$ with values in ${\cal N}$ has image $D_n$ for some odd $n$. Moreover, Klassen observed that if $\Delta_K$ is the Alexander polynomial of $K$,  there are exactly $\frac{|\Delta_K(-1)| - 1}{2}$ characters of such representations (compare Theorem 10 of \cite{Kl}).

A {\it simple point} of a complex affine algebraic set $V$ is a point of $V$ which is contained in a unique algebraic component of $V$ and is a smooth point of that component. 

\begin{lemma} \label{b11}
Let $M$ be the exterior of a knot $K$ in a $\mathbb Z$-homology $3$-sphere
$W$. Suppose that $\rho: \pi_1(M) \to PSL_2(\mathbb C)$ has image $D_n$ where $n > 1$. Then the associated $2n$-fold cover $\widetilde{M}_\rho \to M$ extends to a branched cover $p: \widehat{W}_\rho(K) \to W$, branched over $K$. Moreover \\
$(1)$ $p$ factors through an $n$-fold cyclic (unbranched) cover $\widehat{W}_\rho(K) \to \widehat{W}_2(K)$ and the $2$-fold branched cyclic cover $\widehat{W}_2(K) \to W$. \\
$(2)$ if $b_1(\widehat{W}_\rho(K)) = 0$, then $H^1(M; Ad\rho) \cong \mathbb C$ and $\chi_\rho$ is a simple point of $X_{PSL_2}(M)$. 
\end{lemma} 

\pf (1) Fix meridinal and longitudinal classes $\mu$ and $\lambda$ in $\pi_1(\partial M) \subset \pi_1(M)$. Denote by $\widetilde{M}_2 \to M$ the $2$-fold cover of $M$. Since $W$ is a $\mathbb Z$-homology $3$-sphere we have $b_1(M) = b_1(\widetilde{M}_2) = 1$. 

The subgroup $\rho^{-1}({\cal D})$ has index $2$ in $\pi_1(M)$ and so equals $\pi_1(\widetilde{M}_2)$. Hence $\rho|\pi_1(\widetilde{M}_2)$ has image $D_n \cap {\cal D} \cong \mathbb Z/n$. Since $\pi_1(M)$ is generated by $\mu$ and $\pi_1(\widetilde{M}_2)$ we see that $\rho(\mu) \in {\cal N} \setminus {\cal D}$ and
therefore has order $2$. Further, since $\lambda$ is a double commutator, $\pi_1(\partial \widetilde{M}_2) \subset \hbox{kernel}(\rho)$. In particular if $\widetilde{M}_\rho \to \widetilde{M}_2$ is the regular cover associated to $\rho|\pi_1(\widetilde{M}_2)$, then $|\partial \widetilde{M}_\rho| = n$.

Since $\rho(\mu^2) = \pm I$, $\rho|\pi_1(\widetilde{M}_2)$ factors through $\pi_1(\widehat{W}_2(K))$ and defines an $n$-fold cyclic cover $\widehat{W}_\rho(K) \to \widehat{W}_2(K)$ which composes with $\widehat{W}_2(K) \to W$ to produce the desired  cover of $W$ branched over $K$. It is clear that $\widehat{W}_\rho(K) \to \widehat{W}_2(K)$ is obtained from $\widetilde{M}_\rho \to \widetilde{M}_2$ by equivariant Dehn filling. 

(2) Suppose now that $b_1(\widehat{W}_\rho(K)) = 0$. For each $d \geq 1$ which divides $n$, let $\rho_d = \theta_{n,d} \circ \rho$ and
$\widetilde{M}_{2d} \to M$ the associated cover. The second paragraph of the proof of (1) shows that
$|\partial \widetilde{M}_{2d}| = d$ and since
each boundary component is a torus, $b_1(\widetilde{M}_{2d}) \geq d$. On the other
hand, the third paragraph shows that there is a
Dehn filling of $\widetilde{M}_{2d}$ which yields $\widehat{W}_{\rho_d}(K)$. Now by
construction, $\widehat{W}_{\rho_d}(K)$ is covered
by $\widehat{W}_{\rho}(K)$ and therefore $b_1(\widehat{W}_{\rho_d}(K)) \leq b_1(\widehat{W}_{\rho}(K)) = 0$. Thus $d \leq
b_1(\widetilde{M}_{2d}) \leq b_1(\widehat{W}_{\rho_d}(K)) + |\partial \widetilde{M}_{2d}| \leq d$. Plugging $b_1(\widetilde{M}_{2d}) = d$ into the
conclusion of Lemma \ref{b1Adrho} shows that
$$b_1(M; Ad \rho) = b_1(\widetilde{M}_2) - b_1(M) + \frac{1}{\phi(n)}
\sum_{d|n} \mu(\frac{n}{d}) d = \frac{1}{\phi(n)}
\sum_{d|n} \mu(\frac{n}{d}) d.$$
It is well-known that for $n > 1$ we have $\frac{1}{\phi(n)} \sum_{d|n} \mu(\frac{n}{d}) d = 1$ (see Identity 16.3.1 of \cite{HW} for example) and thus, $b_1(M; Ad \rho) = 1$. Theorem 3 of [BZ3] now shows that $\chi_\rho$ is a simple point of $X_{PSL_2}(M)$. This completes the proof.
\qed

\begin{cor} \label{2br1}
If $K$ is a two-bridge knot with exterior $M$ and $\rho \in R_{PSL_2}(M)$ has image $D_n$ where $n > 1$, then $H^1(M; Ad\rho) \cong \mathbb C$ and $\chi_\rho$ is a simple point of $X_{PSL_2}(M)$.  
\end{cor}

\pf Since the $2$-fold branched cyclic cover of $W = S^3$ over $K$ is a lens space, Proposition \ref{b11}(1) implies that $b_1(\widehat{W}_\rho(K)) = 0$. The desired conclusion now follows from conclusion (2) of Proposition \ref{b11}.
\qed

\subsection{Proof of Proposition \ref{meridiandegree}} \label{pf510} 

Let $p \geq 1,q$ be relatively prime integers where $p$ is odd. 
We observed in \S  \ref{criterion} that given a knot $K \subset S^3$ with exterior $M$, the image of any  homomorphism $\rho: \pi_1(M) \to {\cal N}$ with non-abelian image is $D_n$ for some odd $n \geq 3$. Moreover, the number of characters of such representations is exactly $\frac{|\Delta_K(-1)| - 1}{2} = \frac{|H_1(\hat S^3(K))| - 1}{2} < \infty$. For  $K = k_{p/q}$ we have $|\Delta_K(-1)| = |H_1(L(p,q))| = p$. This discussion yields our next lemma. 

\begin{lemma} \label{strnontriv} 
Every non-trivial curve in the $PSL_2(\mathbb C)$ character variety of the exterior of a knot in $S^3$ is strictly non-trivial. 
\qed
\end{lemma}

Consider a non-trivial curve $X_0 \subset X_{PSL_2}(M_{p/q})$. It is shown in \cite{HT} that the meridinal slope $\mu$ of $k_{p/q}$ is not a boundary slope.  Since $M_{p/q}$ is small, Propositions \ref{idealvalue} shows that for each ideal point $x$ of $\tilde X_0$, $\Pi_x(\tilde f_{\mu}) > 0$. Thus $\Pi_x(\tilde f_{\mu^2})  =  \Pi_x(\tilde f_{\mu}(\tilde f_{\mu} + 4)) > 0$ as well.
Then $Z_x(\tilde f_{\mu}) = Z_x(\tilde f_{\mu^2}) = 0$ and
so by Identity (2.4.1) we have
$$0 < d_M(\mu) = d_M(\mu^2) - d_M(\mu) =
\sum_{\stackrel{non-trivial}{X_0}} \; \sum_{x \in X_0^\nu} (Z_x(\tilde f_{\mu^2})  - Z_x(\tilde
f_{\mu})). \eqno{({\rm A}.3.1)}$$
It follows from Proposition 1.1.3 of \cite{CGLS} that $Z_x(\tilde f_{\mu}) \leq Z_x(\tilde f_{\mu^2})$ for each $x \in X_0^\nu$.
Moreover, Proposition 1.5.4 of that paper shows that if $Z_x(\tilde
f_{\mu}) < Z_x(\tilde f_{\mu^2})$ for some $x \in X_0^\nu$ and $\nu(x) = \chi_\rho$, then
$\rho(\mu^2) = \pm I$. In
particular, the restriction of $\rho$ to the fundamental group of the
$2$-fold cover of $M_{p/q}$
factors through the fundamental group of $L(p,q)$, the $2$-fold cover of
$S^3$ branched over $k_{p/q}$. Thus the image of $\rho$ is finite and as 
we can suppose that it is not cyclic (Proposition 1.5.5 of \cite{CGLS}), it must be 
a non-abelian dihedral group. 

Now $\rho(\mu)$ is neither parabolic nor $\pm I$, so that 
$Z_x(\tilde f_\mu(x)) = 0$. We know that $\chi_\rho$ is a simple point of $X_{PSL_2}(M_{p/q})$ by Corollary \ref{2br1} and we claim (see Lemma \ref{simple} below)
that $Z_x(\tilde f_{\mu^2}) = 1$. 
Note that these two facts and Identity (A.3.1) show that
$d_M(\mu)$ equals the number
of irreducible, dihedral characters which lie on some non-trivial curve in
$X_{PSL_2}(M_{p/q})$. But by
Proposition \ref{posdim}, every such character lies on such a curve, and
since there are $\frac{p-1}{2}$
irreducible, dihedral characters of $\pi_1(M_{p/q})$, we have 
$d_M(\mu) = \frac{p-1}{2}$,  which is what we set out to
prove.
\qed

\begin{lemma} \label{simple}
Let $\chi_\rho$ be an irreducible, dihedral character of $\pi_1(M_{p/q})$, $X_0$ the
unique curve in $X_{PSL_2}(M_{p/q})$ which contains it, and $x$ the unique point of $X_0^\nu$ such that $\nu(x) = \chi_\rho$. Then $Z_x(\tilde f_{\mu^2}) = 1$.
\end{lemma}

\pf The proof that $Z_x(\tilde f_{\mu^2}) = 1$ is essentially identical to the proof of Theorem 2.1 (2) of \cite{BB}, though with with some slight modifications as $2\mu$ is not a
primitive class in $H_1(\partial M)$. These modifications are simple and we describe them next.

Let $M_{p/q}(2\mu)$ be the space obtained by attaching a solid torus to
$M_{p/q}$ by a covering map which maps $S^1 \times \{1\}$ homeomorphically
to $\lambda_{M_{p/q}}$ and is a
$2$-fold cover of $\{*\} \times \partial D^2$ to $\mu$. Note that there is
a unique $2$-fold cover of
$M_{p/q}(2\mu)$ and its total space is $L(p,q)$, the $2$-fold cover of
$S^3$ branched over
$k_{p/q}$. Note as well that $\rho$ factors through $\pi_1(M_{p/q}(2\mu))$.
Lemma \ref{b1Adrho}
applied to this situation shows that $H^1(M_{p/q}(2\mu)) = 0$. The proof of
Lemma 1.8 of [BB] shows that
if $u \in Z^1(\pi_1(M_{p/q}); Ad\rho)$ represents a non-zero class in
$H^1(\pi_1(M_{p/q});
Ad\rho)$ then $u(\mu^2) \ne 0$. The calculation $Z_x(\tilde f_{\mu^2}) = 1$
now follows in a similar fashion to the calculation in the proof of Theorem
2.1 (2) of [BB].   
\qed

\section{Peripheral values of homomorphisms of twist knot groups} \label{twistknotperipheral}

In this appendix we show that the trefoil knot is the only twist knot whose group admits a homomorphism onto an infinite triangle group such that the image of the peripheral subgroup is finite. 

After Hoste and Shanahan, we identify the $n$-twist knot $K_n$ with the knot $J(2, 2n)$ of \cite{HS1}. When $n = -1, 0, 1$, $K_n$ is the figure $8$ knot, the trivial knot, and the trefoil knot respectively. 

The fundamental group of the exterior $M_n$ of $K_n$ admits a presentation 
$$\pi_1(M_n) = \langle a, b : a (ab^{-1}a^{-1}b)^n = (ab^{-1}a^{-1}b)^n b \rangle$$ 
where $a$ and $b$ are meridinal classes (cf. Proposition 1 of \cite{HS1}). Set $w = ab^{-1}a^{-1}b$ so that the relation becomes $a w^n = w^n b$.

\begin{prop}\label{nil}
If there is a surjective homomorphism $\rho: \pi_1(M_n) \to \Delta(p,q,r)$ where $(p,q,r)$ is a Euclidean triple, then $n = 1$ and hence $K_n$ is the trefoil knot. Further, $(p,q,r) = (2,3,6)$ up to permutation.
\end{prop}

\pf First we observe that any two elements of $\Delta(p,q,r)$ which are of infinite order commute since they correspond to translations under the natural embedding $\Delta(p,q,r) \to \hbox{Isom}_+(\mathbb E^2)$. Thus $\rho(a)$ and $\rho(b) = \rho(w)^{-n}\rho(a) \rho(w)^n$ must be elliptic. Since $H_1(\Delta(p,q,r))$ is necessarily cyclic, we have $(p,q,r) = (2,3,6)$ up to permutation. Further, $\rho(a)$ generates $H_1(\Delta(2,3,6)) \cong \mathbb Z/6$ so that $\rho(a)$ and $\rho(b)$ have order $6$. Then up to replacing $\rho$ by a conjugate representation we may suppose that 
$$\rho(a) = xy \in \langle x, y : x^2 = y^ 3 = (xy)^6 = 1 \rangle = \Delta(2,3,6).$$
We claim that up to conjugating $\rho$ by a power of $xy$, we can suppose that $\rho(b) = yx$. To see this, fix a tessellation ${\cal T}$ of $\mathbb E^2$ by triangles with angles $\frac{\pi}{2}, \frac{\pi}{3}, \frac{\pi}{6}$, and identify $\Delta(2,3,6) \subset \hbox{Isom}(\mathbb E^2)$ with the set of orientation preserving symmetries of ${\cal T}$.  
the elements $xy, yx$ are conjugate elements of order $6$ in $\Delta(2,3,6)$ and form a generating set. Moreover, the tessellation ${\cal T}$ can be described as follows. Let $A$ be the fixed point of $xy$, $B \ne A$ that of $yx$, and let $C$ be the midpoint of $[A, B]$. Denote by $L$ the line through $C$ which is orthogonal to $[A, B]$ and let $T(A, D, E)$ be the triangle with vertices $A, D, E$ where $D, E \in L$ are equidistant to $C$ and the angles at $A, D, E$ are $\frac{\pi}{3}, \frac{\pi}{3}, \frac{\pi}{6}$ respectively. The triangle 
$T(A,C,D)$ is a face of ${\cal T}$ and so the tessellation is its orbit  under the action of $\Delta(2,3,6)$. Moreover, $T(A, D, E)$ is the union of the two adjacent faces $T(A,C,D)$ and $T(A,C,E)$ and so is a fundamental domain for $\Delta(2,3,6)$. Since $\mathbb E^2$ admits self-similarities of arbitrary scale factor, it is clear that any two elements of order $6$ in $\hbox{Isom}(\mathbb E^2)$ with distinct fixed points generate a subgroup isomorphic to $\Delta(2,3,6)$ with invariant tessellation and fundamental domain constructed  as above. In particular this is the case for $\rho(a), \rho(b)$. Let ${\cal T}'$ and $T(A, C', E')$ be the associated tessellation and fundamental domain. Since $\rho(a), \rho(b)$ generate $\Delta(2,3,6)$, ${\cal T} = {\cal T}'$ and so $T(A, C, E)$ can be obtained from $T(A, C', E')$ by a rotation about $A$ of angle $\frac{2 \pi j}{3}$ for some integer $j$. This rotation is given by $(xy)^{\epsilon j}, \epsilon \in \{\pm 1\}$, so if we replace $\rho$ by $(xy)^{-\epsilon j} \rho (xy)^{\epsilon j}$, the new fixed point of $\rho(b)$ is $B$. Since $\rho(a)$ and $\rho(b)$ are conjugate, it follows that $\rho(b) = yx$. 

With these calculations in hand, we see that $v := [a, b^{-1}] = x(yx)^3$ is a product of two elements of order $2$ with distinct fixed points. Thus $v$ is a translation which leaves ${\cal T}$ invariant. Since $\rho(b) = v^{-n} \rho(a) v^n$, $B$, the fixed point of $\rho(b)$, equals $v^{-n}(A)$. But examination of ${\cal T}$ shows the only way this is possible is for $n = 1$. 
\qed

\begin{prop}\label{prop:epi}
If there is a surjective homomorphism $\rho: \pi_1(M_n) \to \Delta(p,q,r) \subset PSL_2(\mathbb R)$ such that $(p,q,r)$ is a hyperbolic triple and $\rho(a)$ is elliptic, then $n = 1$, so that $K_n$ is the trefoil knot. Further, $(p,q,r) = (2,3,r), \, r \geq 7$ up to permutation.
\end{prop}

\pf Suppose that there is a surjective homomorphism $\rho: \pi_1(M_n) \to \Delta(p,q,r) \subset PSL_2(\mathbb R)$ such that $(p,q,r)$ is a hyperbolic triple and $\rho(a)$ is elliptic. Clearly $n \ne 0$ and $\Delta(p,q,r)$ is generated by two conjugate elliptics. Theorem 2.3 of \cite{Kn} implies that up to permuting $p,q,r$, one of the following two scenarios arises. 
\vspace{-.4cm}  
\begin{description}
\item[{\rm (a)}] $(p,q,r) = (2,q, r)$ where $\rho(a)$ has order $q$ and $r \geq 3$ is odd. Further, there is an integer $s$ relatively prime to $q$ such that in the standard presentation $\Delta(2, q,r) = \langle x, y, z : x^2, y^q, z^r \rangle$ we have $\rho(a) = y^s, \rho(b) = xy^sx^{-1}$. 
\item[{\rm (b)}] $(p,q,r) = (2, 3, r)$ where where $\rho(a)$ has odd order $r \geq 7$. Further,  there is an integer $s$ relatively prime to $r$ such that in the standard presentation $\Delta(2, 3, r) = \langle x, y, z : x^2, y^3, z^r \rangle$ we have $\rho(a) = z^s, \rho(b) = yxy^{-1}z^syxy^{-1}$. 
\end{description}
Set $v = \rho(w)$ so that  
$$v = [\rho(a), \rho(b^{-1})] = \left\{ \begin{array}{ll} (y^sxy^{-s})((xy^{-s})x(xy^{-s})^{-1}) & \hbox{in scenario (a)} \\ ((z^sy)x(z^sy)^{-1})((yxy^{-1}z^{-s}y)x(yxy^{-1}z^{-s}y)^{-1}) & \hbox{in scenario (b)}. \end{array} \right.$$ 
In either case $v = uu'$ where $u, u'$ are of order $2$. Denote by $R, R'$ the fixed points of $u, u'$ and observe that if $R = R'$, then $u = u'$ and therefore $v = uu' = 1$. Then $\rho(b) = \rho(w^{-n}aw^n) = v^{-n} \rho(a) v^n = \rho(a)$, which is impossible. Thus $R \ne R'$ and it is easy to see that if $\gamma$ denotes the geodesic in $\mathbb H^2$ which contains both $R$ and $R'$, then $v$ is a hyperbolic element of $PSL_2(\mathbb R)$ with invariant geodesic $\gamma$ and translation length $2d_{\mathbb H^2}(R, R')$. The proof of the proposition is similar in the two possible scenarios. We analyse them separately. 

Assume first that we are in scenario (a). There is a fundamental domain for $\Delta(2,q,r)$ in $\mathbb H^2$ which is a geodesic triangle $T = T(A,B,C)$ having vertices $A = \hbox{Fix}(y), B = \hbox{Fix}(xyx^{-1}), C = \hbox{Fix}(z)$ and the midpoint $P$ of $[A, B]$ is $\hbox{Fix}(x)$ (so $x(A) = B$). The angles of $T$ at $A, B, C$ are $\frac{\pi}{q}, \frac{\pi}{q}, \frac{2 \pi}{r}$ respectively. The hyperbolicity of $v$ implies that for $l \ne 0$, $v^l(C), v^l(P) \not \in T$ and so if $T \cap v^l(T) \ne \emptyset$, then up to replacing $l$ by its negative we have $T \cap v^l(T) = \{A\}$ and $v^l(B) = A$. 

Since $xy^sx = \rho(b) = \rho(w^{-n}aw^n) = v^{-n} y^s v^n$, there is an integer $m$ such that $v^n = y^mx$. Then $v^n(B) = y^mx(B) = A$ so that $v^n(T) \cap T= \{A\}$. It now follows from the previous paragraph that if for some $l \ne 0$ we have $T \cap v^l(T) \ne \emptyset$, then up to replacing $l$ by its negative we have $v^l(B) = A = v^n(B)$. Thus $l = \pm n$ and so for $d \ne e$, 
$$v^{dn}(T) \cap v^{en}(T) = \left\{\begin{array}{cl} 
v^{dn}(B) & e  = d - 1 \\
v^{(d+1)n}(B) & e = d + 1 \\
\emptyset & e \ne d \pm 1 \end{array} \right.$$ 
and the reader will verify that $\Gamma_0 = \cup_d v^{dn}(T)$ is an infinite chain of geodesic triangles which is closed, properly embedded, and separating in $\mathbb H^2$.  
It follows that $n = \pm 1$ as otherwise $v(\Gamma_0) \cap \Gamma_0 = \emptyset$ and so the side of $\Gamma_0$ in $\mathbb H^2$ containing $v(\Gamma_0)$ is invariant under $v$. Thus $v^l(\Gamma_0) \cap \Gamma_0 = \emptyset$ for all $l > 0$, contrary to the fact that $v^{|n|}(\Gamma_0) = \Gamma_0$. 

If $n = -1$, then $y^{m}x = v^{-1}$ so that $xy^{m}x^{-1} = (y^{-s}xy^{s})x(y^{s}xy^{-s})$. In particular, $(y^{-s}xy^{s})x(y^{s}xy^{-s})$ fixes $x(A) = B$. But $(y^{-s}xy^{s})x(y^{s}xy^{-s})$ is a product of three conjugates of $x$ with fixed points $P, y^s(P), y^{-s}(P)$ and the reader will verify that this is impossible because such a configuration of order 2 elliptics cannot fix $B$. (Alternately, we refer the reader to the proof of Theorem 11.5.2 of [Beardon] where the fixed points of a product of three order 2 elliptics are analysed. The analysis implies that if $(y^{-s}xy^{s})x(y^{s}xy^{-s})$ has a fixed point, then this fixed point and $B$ lie on opposite sides of the geodesic through $y^s(P)$ and $y^{-s}(P)$.)  

Finally, if $n = 1$, $K_n$ is the trefoil knot and we have $y^{m}x = v = y^sxy^{-s}xy^{-s}xy^{s}x$ so that $y^{m-3s} = (xy^{-s})^3$. If $(xy^{-s})^3 \ne 1$ then $xy^{-s}$ fixes $A$, which is impossible. Thus $y^{m-3s} = (xy^{-s})^3 = 1$. We will show that $r = 3$ to complete this part of the proof. Let $D$ be the fixed point of $xy^{-s}$. Since $xy^{-s}(A) = B$, $D$ lies on the perpendicular $L$ to $[A,B]$ through $P$. Now $D \ne P$ as otherwise $y^s = 1$. On the other hand, $C$ and $x(C)$ are the closest points to $P$ of the given tessellation of $\mathbb H^2$ which lie on $L$. Further, since $z(B) = (xz^{-1})(B) = A$ we see that $\frac{2\pi}{3}$, the absolute value of the angle of rotation of $xy^{-s}$ is bounded above by that of $z$ at $C$ or $xz^{-1}x$ at $x(C)$, which is $\frac{2\pi}{r}$, with equality if and only if $D \in \{C, x(C)\}$. Thus $\frac{2\pi}{3} \leq \frac{2\pi}{r}$ which implies that $r=3$ and we are done.  

Now suppose that we are in scenario (b) and consider the geodesic triangle $T_0$ in $\mathbb H^2$ with vertices $A, B, C$ such that $A = \hbox{Fix}(x), B = \hbox{Fix}(y), C = \hbox{Fix}(z)$. The angles of $T_0$ at $A, B, C$ are $\frac{\pi}{2}, \frac{\pi}{3}, \frac{\pi}{r}$ respectively where $r \geq 7$ is odd. Recall that 
$$v = [(z^sy)x(z^sy)^{-1}][(yxy^{-1}z^{-s}y)x(yxy^{-1}z^{-s}y)^{-1}]$$ 
is a product of two conjugates of $x$ and observe that that form of $\rho(b)$ given in scenario (b) implies that 
$$v^n = z^m yxy^{-1}$$ 
for some integer $m$. 

Consider the geodesic triangle $T = T(C,D, E)$ containing $T_0$ with vertices $C, D = (yxy^{-1})(C) = \hbox{Fix}((yxy^{-1})z(yxy^{-1})), E = y(C) = \hbox{Fix}(yzy^{-1})$ of angles $\frac{\pi}{r}, \frac{\pi}{r}, \frac{4\pi}{r}$ respectively. 

\begin{claim} \label{intersection}
For any integer $l \ne 0$, $v^l(T) \cap T$ is either empty or one of the vertices $C, D$.
\end{claim}

\pf (of Claim \ref{intersection}) First we show that $v^l(\hbox{int}(T)) \cap \hbox{int}(T) = \emptyset$. Decompose $T$ as $T_1 \cup T_2 \cup T_3$ where $T_1 = T(B,C,E), T_2 = T(B, E, B')$ where $B' = (yxy^{-1})(B)$, and $T_3 = T(B', D, E)$ and note that each of $T_1, T_2, T_3$ is a fundamental domain for the action of $\Delta(2,3,r)$ on $\mathbb H^2$. Since $T_1$ is sent to the geodesic triangle $T_3$ by the elliptic element $(yxy^{-1})(xzx^{-1})(yxy^{-1})^{-1}$, $v^l(\hbox{int}(T_1)) \cap \hbox{int}(T_3) = \emptyset$. Similarly $v^l(\hbox{int}(T_1)) \cap \hbox{int}(T_2) = \emptyset$ so that $v^l(\hbox{int}(T_1)) \cap \hbox{int}(T) = \emptyset$. In the same way we see that $v^l(\hbox{int}(T_2)) \cap \hbox{int}(T) = v^l(\hbox{int}(T_3)) \cap \hbox{int}(T) = \emptyset$, which is what we needed to prove. 

Second we claim that $v^l(E) \not \in T$. If this is false we have $v^l(E) \in \{C, D\}$ (i.e. the only valency $2r$ vertices in $T$ are $C,D,E$ and $v^l(E) \ne E$ as it is hyperbolic), say $v^l(E) = C$. Then the axis of $v$ is perpendicular to the perpendicular bisector $L$ of $[C, E]$. On the other hand, $v^n(D) = z^m yxy^{-1}(D) = C$, so the axis of $v$ is also perpendicular to the perpendicular bisector $L'$ of either $[C, D]$. But this is impossible since $L \cap L' = \{y(E)\} \ne \emptyset$. Hence $v^l(E) \ne C$ and a similar argument shows it does not equal $D$. 

Third we observe that $v^l(T) \cap T$ contains no edge of the tessellation. By the first paragraph, such an edge would have to lie in $\partial T$ and by the second it could not contain $E$. Since $v^l$ preserves the combinatorial type of the vertices it is now easy to use the method of the first paragraph to obtain a contradiction. 

These observations imply that if $v^l(T) \cap T$ is non-empty, then it is a vertex of $T$. This proves Claim \ref{intersection} 
\qed

Since $v^n(D) = C$, the claim implies that if $v^l(T) \cap T \ne \emptyset$, then after possibly replacing $l$ by its negative we have $v^l(D) = C$. Thus $l = \pm n$ and the intersection is $C$ if $l = n$ and $D$ otherwise. Arguing as in scenario (a) we have $|n| = 1$. If $n = -1$, we have $z^myxy^{-1} = v^{-1}$ and therefore 
$$(yxy^{-1})z^m(yxy^{-1}) = [(z^{-s}y)x(z^{-s}y)^{-1}](yxy^{-1})[(z^sy)x(z^sy)^{-1}].$$
Hence $D$ is fixed by the product $[(z^{-s}y)x(z^{-s}y)^{-1}](yxy^{-1})[(z^sy)x(z^sy)^{-1}]$ of three conjugates of $x$ with fixed points $y(P), (z^{s}y)(P), (z^{-s}y)(P)$. As in the analysis of the case $n = -1$ in scenario (a), it can be verified that this cannot occur. Thus $n = 1$ and $K_n$ is the trefoil knot. 
\qed  

\section{Bending} \label{bending} 

Let $\Gamma$ be a finitely generated group which splits over a subgroup $\Gamma_0$ and $\rho: \Gamma \to PSL_2(\mathbb C)$ a homomorphism such that $\rho(\Gamma_0)$ is abelian but not isomorphic to $\mathbb Z/2 \oplus \mathbb Z/2$. Under this condition we can perform a deformation operation on $\chi_\rho$ known as bending. The details of the  construction depend on whether the splitting is a free product with amalgamation or an HNN extension and are dealt with in \S \ref{fpa} and \S \ref{hnn} respectively. 

Recall the subgroups ${\cal D, N}$ of $PSL_2(\mathbb C)$ defined in \S \ref{generalities} and set 
$${\cal P}_+ = \{\pm \left( \begin{array}{cc} 1 & z \\ 0 & 1 \end{array}
\right) \; | \; z \in {\mathbb C}\} \subset {\cal T}_+ = \{\pm \left( \begin{array}{cc} z & w \\ 0 & z^{-1} \end{array}
\right)  \; | \; z, w \in {\mathbb C}, z \ne 0\}$$ 
Under the natural action of $PSL_2(\mathbb C)$ on $\mathbb CP^1$, the fixed point sets of ${\cal T}_+$ and ${\cal P}_+$ coincide and consist of a single line $L_+$. That of ${\cal D}$ consists of two lines $\{L_+, L_-\}$. 

The centraliser of a subset $E$ of $PSL_2(\mathbb C)$ wil be denoted by $Z_{PSL_2}(E)$ and the component of the identity of $Z_{PSL_2}(E)$ will be denoted 
$Z_{PSL_2}^0(E)$. For $\pm I \ne A \in PSL_2(\mathbb C)$ we have 
$$Z_{PSL_2}^0(A) \mbox{ is conjugate to } 
\left\{ \begin{array}{ll}
{\cal D} & \mbox{if $A$ is diagonalisable} \\
{\cal P}_+ & \mbox{if  $A $ is parabolic}
\end{array} \right. $$
Thus when $E \ne \{\pm I\}$, $Z_{PSL_2}^0(E)$ is abelian and reducible. 
For a group $\pi$ and representation $\rho \in R_{PSL_2}(\pi)$, we use $Z_{PSL_2}(\rho), Z_{PSL_2}^0(\rho)$ to denote, respectively, the centraliser and  the component of the identity of the centraliser of $\rho(\pi)$.

\subsection{$\Gamma = \Gamma_1 *_{\Gamma_0} \Gamma_2$} \label{fpa}

Fix $\rho: \Gamma \to PSL_2(\mathbb C)$ a homomorphism such that $\rho(\Gamma_0)$ is abelian but not isomorphic to $\mathbb Z/2 \oplus \mathbb Z/2$. Denote by $\rho_j$ the restriction of $\rho$ to $\Gamma_j $. For each $S \in Z_{SL_2}(\rho_0)$ define $\rho_S: \Gamma \to PSL_2(\mathbb C)$ to be the homomorphism determined by the push-out diagram 

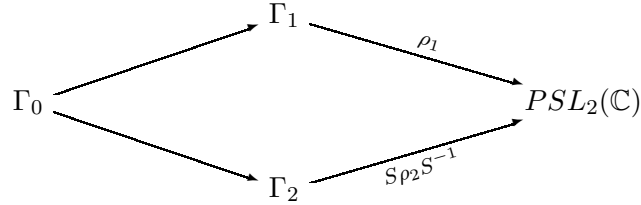
\vspace{-.4cm}
\begin{diagram}[height=1.5em,w=4em]
      & &  \Gamma_1 & &      \\
      &\;\;\;\;\;\;\;\;\;\;\;\;\;\; \ruTo \;\;\;\;\;\;\;\;\;\;\;\;\;\;&    &  
\;\;\;\;\;\;\;\;\;\;\;\;\;\;\rdTo^{\;\;\;\rho_1} 
\;\;\;\;\;\;\;\;\;\;\;\;\;\; &              \\
 \Gamma_0  &   &  &   & PSL_2(\mathbb C)                 \\
 & \rdTo   &             &   \ruTo_{S \rho_2 S^{-1}}                  \\
 & & \Gamma_2 &   
\end{diagram}

\vspace{.2cm} 
\noindent We say that the character $\chi_{\rho_S}$ is obtained by {\it bending} $\chi_\rho$ by $S$. The {\it bending function} 
$$\beta_\rho: Z^0_{PSL_2}(\rho_0) \to X_{PSL_2}(\Gamma), \; S \mapsto \chi_{\rho_S}.$$
We say that $\rho$ can be {\it bent non-trivially} if $\beta_\rho$ is non-constant. Our next result determines necessary and sufficient conditions for this to occur. 

\begin{lemma} \label{constbendsep}
Suppose that $\rho \in R_{PSL_2}(\Gamma$ is such that $\rho(\Gamma_0)$ is abelian but not isomorphic to $\mathbb Z/2 \oplus \mathbb Z/2$. The bending function $\beta_\rho: Z^0_{PSL_2}(\rho_0) \to X_{PSL_2}(\Gamma)$ is constant if and only if one of the following two situations arises: \\ 
\indent $(a)$ $\rho_0(\Gamma_0) = \{\pm I\}$ and either $\rho_1(\Gamma_1) = \{\pm I\}$ or $\rho_2(\Gamma_2) = \{\pm I\}$. \\ 
\indent $(b)$ $\rho_0(\Gamma_0) \ne \{\pm I\}$ and either $\rho_1(\Gamma_1)$ is abelian and reducible, or $\rho_2(\Gamma_2)$  is abelian and \\ \indent \hspace{5mm} reducible, or $\rho$ is reducible.  
\end{lemma}

\pf  Suppose that the correspondence $S \mapsto \chi_{\rho_S}$ is constant. We leave the justification of the following claim to the reader. 

\begin{claim} \label{constconj}
{\rm Let $A, B, C \in SL_2(\mathbb C)$. Then $\mbox{trace}(ASBS^{-1}) =\mbox{trace}(AB)$ for all $S \in Z_{SL_2}(C)$ if and only if one of the following two situations arise: \\ 
\indent $(a)$  $C = \pm I$ and either $A = \pm I$ or $B = \pm I$. \\ 
\indent $(b)$ $C \ne \pm I$ and either $[A,C] = I$ or
$[B,C] = I$ or $A,B$ and $C$  have a common \\ \indent \hspace{5mm} eigenvector.  } 
\qed 
\end{claim}

\noindent  Set $G_j = \mbox{ image}(\rho_j)$ for $j = 0,1,2$. The claim shows that if $G_0 = \{\pm I\}$, then either $G_1 = \{\pm I\}$ or $G_2 = \{\pm I\}$, and we are done.  Assume then that $G_0 \ne \{\pm I\}$ and that both $G_1$ and $G_2$ are either irreducible or non-abelian subgroups of $PSL_2(\mathbb C)$. It follows that neither $G_1 \subset Z^0_{PSL_2}(G_0)$ nor $G_2 \subset Z^0_{PSL_2}(G_0)$. We will show that $\rho$ is reducible.  

Fix $C \in G_0 \setminus \{\pm I\}$ and observe that our hypotheses show 
$$Z^0_{PSL_2}(C) = Z^0_{PSL_2}(G_0) = \left\{\begin{array}{ll} {\cal D} & \mbox{if } G_0 \subset {\cal D} \\  
{\cal P}_+ & \mbox{if } G_0 \subset {\cal P}_+ \end{array} \right.$$
If there is some $A_0 \in G_1 \setminus Z^0_{PSL_2}(G_0)$, then Claim \ref{constconj} implies that for each $B \in G_2$, 
either $B \in Z^0_{PSL_2}(G_0) \subset {\cal T}_+$ or $A_0, B$ and $C$ have a common fixed point in $\mathbb CP^1$. 
It follows that each element of $G_2$ fixes at least one of $L_+$ and $L_-$ and it is simple to deduce from this that 
$G_2$ fixes one of these lines. A similar argument shows that $G_1$ fixes one of them as well. 
If $G_1$ and $G_2$ have a common fixed point, then $\rho$ is reducible, so suppose that they do not. One of them, say $G_1$,
fixes $L_+$ and not $L_-$, while $G_2$ fixes $L_-$ and not $L_+$. Since $G_0 \subset G_1 \cap G_2$, it fixes both $L_+$
and $L_-$ and therefore we must have $G_0 \subset {\cal D}$. By choice, $A_0 \in G_1 \setminus {\cal D}$ and so its fixed
point set is $L_+$. On the other hand we have assumed that there is some  
$B_0 \in G_2 \setminus Z^0_{PSL_2}(G_0) = G_2 \setminus {\cal D}$. Its fixed point set is $L_-$. 
But this is impossible as Claim \ref{constconj} implies that $A_0, B_0$ and $C$ have a common fixed point in $\mathbb CP^1$. 
Thus $G_1$ and $G_2$ do have a common fixed point, and so $\rho$ is reducible. 

Conversely if either $G_1 \subset Z^0_{PSL_2}(G_0)$, or
$G_2 \subset Z^0_{PSL_2}(G_0)$, or $G_0 \ne \{\pm I\}$ and $\rho$ is reducible, 
then Claim \ref{constconj} implies that the correspondence $S \mapsto \chi_{\rho_S}$, 
where $S \in Z^0_{PSL_2}(G_0)$, is constant. This completes the proof of the lemma.  
\qed 

\subsection{$\Gamma = (\Gamma_1)_{\Gamma_0}$}  \label{hnn} 

In this case there is an injective homomorphism $\varphi: \Gamma_0 \to \Gamma$ such that 
$$\Gamma = \langle \Gamma_1, \mu : \mu \gamma \mu^{-1} = \varphi(\gamma), \gamma \in \Gamma_0 \rangle.$$
Set $\Gamma_0' = \varphi(\Gamma_0)$ and for $\rho \in R_{PSL_2}(\Gamma)$ we take $\rho_1, \rho_0, \rho_0'$ to be its restriction to $\Gamma_1, \Gamma_0, \Gamma_0'$ respectively. The correspondence $\rho \in R_{PSL_2}(\Gamma) \mapsto (\rho_1, \rho(\mu)) \in R_{PSL_2}(\Gamma_1) \times PSL_2(\mathbb C)$ determines an identification 
$$R_{PSL_2}(\Gamma) = \{(\rho_1, A) \in R_{PSL_2}(\Gamma_1) \times PSL_2(\mathbb C) : A \rho_0(\gamma) A^{-1}  = \rho_0'(\varphi(\gamma)) \mbox{ for all } \gamma \in \Gamma_0 \}.$$
Note that $(\rho_1, A), (\rho_1, B) \in R_{PSL_2}(\Gamma)$ if and only if $B = AS$ for some $S \in Z_{PSL_2}(\rho(\Gamma_0))$. In particular, if $\rho(\Gamma_0)$ is abelian but not $\mathbb Z/2 \oplus \mathbb Z/2$, we have a bending function 
$$\beta_{(\rho_1, A)}: Z^0_{PSL_2}(\rho(\Gamma_0)) \to X_{PSL_2}(\Gamma_0), S \mapsto \chi_{(\rho, AS)}.$$ 

\begin{lemma} \label{constbendnonsep}
Suppose that $(\rho_1, A) = \rho \in R_{PSL_2}(\Gamma)$ is a representation with $\rho(\Gamma_0)$ is abelian but not $\mathbb Z / 2 \oplus \mathbb Z /2$. The bending function $\beta_{(\rho_1, A)}$ is constant if and only if $\rho(\Gamma_0) \ne \{\pm I\}$ and, after a possible conjugation, one of the following two situations arises:  \\
\indent $(a)$ $\rho(\Gamma_1) \subset {\cal D}$ and $A = \pm \left[ \begin{array}{cc}  0 & 1 \\ -1 & 0 \end{array} \right]$. \\
\indent $(b)$  $\rho(\Gamma_1) \subset {\cal T}_+$ and $A \in {\cal T}_+$. \\ 
In particular, $\rho_1$ is reducible and $\rho$ is either reducible or conjugate into ${\cal N}$.  
\end{lemma}

\pf Let $G_1, G_0, G_0'$ denote the images of $\rho_1, \rho_0, \rho_0'$ respectively. 
After possibly replacing $\rho = (\rho_1, A)$ by a conjugate representation, we may suppose that either $G_0 = \{\pm I\}$, or  $\{\pm I\} \ne G_0  \subset {\cal D}$, or $\{\pm I\} \ne G_0 \subset {\cal P}$. We consider these three cases separately.

\noindent {\bf Case 1}. $G_0  = \{\pm I\}$. 
  
\noindent Then $Z_{PSL_2}^0(\rho_0) = PSL_2(\mathbb C)$ and so in general, $\beta_{(\rho_1, A)}(\mu) =  \pm \mbox{trace}(A) \ne \pm \mbox{trace}(AS)= \beta_{(\rho_1, AS)}(\mu)$ for $S \in Z_{PSL_2}^0(\Gamma_0)$, $\beta_{(\rho_1, A)}$ is not constant. 

\noindent {\bf Case 2}. $\{\pm I\} \ne G_0  \subset {\cal D}$.

\noindent Then $Z_{PSL_2}^0(\rho_0) = {\cal D}$. If $\beta_{(\rho_1, A)}$ is constant, then $\pm \mbox{trace}(\rho(\gamma)A) = \beta_{(\rho_1, A)}(\gamma\mu) = \beta_{(\rho, AS)}(\gamma\mu) = \pm \mbox{trace}(\rho(\gamma)AS)$ for each $\gamma \in \Gamma_1$ and $S \in  {\cal D}$. It follows that $\rho(\gamma)A \in {\cal N} \setminus {\cal D}$ for each $\gamma \in \Gamma_1$. Hence $A \in {\cal N} \setminus {\cal D}$ and $\rho(\Gamma_1) \subset {\cal D}$. After a further conjugation we may suppose that $A = \pm \left[ \begin{array}{cc}  0 & 1 \\ -1 & 0 \end{array} \right]$. 

Conversely suppose that $A = \pm \left[ \begin{array}{cc}  0 & 1 \\ -1 & 0 \end{array} \right]$ and $\rho(\Gamma_1) \subset {\cal D}$. Consider a word $w = \Pi_j \mu^{a_j} x_j$ where $a_j \in \mathbb Z$ and $x_j \in \Gamma_1$. Set $D_j = \rho(x_j) \in {\cal D}$. Then for any $S \in {\cal D}$ we have 
$$(\rho, AS): w \mapsto \Pi_{j} (AS)^{a_j} D_j = (AS)^{a_1 + a_2 + \ldots + a_n} 
\Pi_j D_j^{(-1)^{a_{j+1} + a_{j+2} + \ldots + a_n}}.$$ 
The trace of the right-hand side of this identity is independent of $S$, so that $\beta_{(\rho_1, A)}$ is constant. 
 
\noindent {\bf Case 3}. $\{\pm I\} \ne G_0  \subset {\cal P}$.

\noindent Then $Z_{PSL_2}^0(\rho_0) = {\cal P}$. If $\beta_{(\rho_1, A)}$ is constant, then $\mbox{trace}(\rho(\gamma)A) = \beta_{(\rho_1, A)}(\gamma \mu) = \beta_{(\rho_1, AS)}(\gamma \mu) = \mbox{trace}(\rho(\gamma)AS)$ for each $\gamma \in \Gamma_1$ and $S \in  {\cal P}$. It follows that $\rho(\gamma)A \in {\cal T}_+$ for each 
$\gamma \in \Gamma_1$. Hence $A \in {\cal T}_+$ and $\rho(\Gamma_1) \subset {\cal T}_+$. 

Conversely suppose that $A \in {\cal T}_+$ and $\rho(\Gamma_1) \subset {\cal T}_+$. Consider a word $w = \Pi_{j} \mu^{a_j} x_j$ where $a_j \in \mathbb Z$ and $x_j \in \Gamma_1$. Set $U_j = \rho(x_j) \in {\cal T}_+$. Then for any $S \in {\cal P}$ we have 
$$(\rho, AS): w \mapsto \Pi_{j} (AS)^{a_j} U_j.$$ 
The trace of the right-hand side of this identity is independent of $S$, so that $\beta_{(\rho_1, A)}$ is constant. 

This completes the proof of the lemma. 
\qed

\medskip\noindent

\noindent Michel Boileau, Laboratoire \'Emile Picard, Universit\'e Paul Sabatier, Toulouse Cedex 4, France \\  
e-mail: boileau@picard.ups-tlse.fr 

\noindent
Steven Boyer, D\'ept. de math., UQAM, P. O. Box 8888, Centre-ville, Montr\'eal, Qc, H3C 3P8, Canada
\newline\noindent
e-mail: boyer@math.uqam.ca


\footnotesize
\begin{thebibliography}{XXXXX}

\bibitem[Ad]{Ad} C. Adams, \emph{Unknotting tunnels in hyperbolic 3-manifolds}, Math. Ann.
\textbf{302} (1995), 177--195. 

\bibitem[BMS]{BMS} G. Baumslag, J. Morgan and P. Shalen, \emph{Generalized triangle groups}, Math. Proc. Camb. Phil. Soc. {\bf 102} (1987), 25--31. 


\bibitem[B]{B} L. Ben~Abdelghani.
\emph{Espace des repr{\'e}sentations du groupe d'un n{\oe}ud dans un
 groupe de {L}ie}. Th\`{e}se, Universit\'{e} de Bourgogne, 1998.

\bibitem[BB]{BB} L. Ben~Abdelghani and S. Boyer, \emph{A calculation of the
Culler-Shalen seminorms
associated to small Seifert Dehn fillings}, Proc. London Math. Soc.
\textbf{83} (2001), 235--256.

\bibitem[BeP]{BeP} R. Benedetti and C. Petronio.
 \emph{Lectures on hyperbolic geometry}, Universitext, Springer, Berlin, 1992.


\bibitem[BoP]{BoP} M. Boileau and J. Porti, \emph{Geometrization of
3-orbifolds of cyclic type}, Ast\'eris\-que Monograph, \textbf{272}, 2001.

\bibitem[BWa1]{BWa1} M. Boileau and S.C. Wang, \emph{Non-zero degree maps and
surface bundles over $S^1$}, J. Diff. Geom. \textbf{43} (1996), 789--806.

\bibitem[BWa2]{BWa2} M. Boileau and S.C. Wang, \emph{Degree one maps
between small
3-manifolds and Heegaard genus}, arXiv:math.GT/0402121.

\bibitem[BWe]{BWe} M. Boileau and R. Weidmann, \emph{The structure of
3-manifolds with 2-generated fundamental group}, Topology \textbf{44} (2005), 283--320.

\bibitem[BolZ]{BolZ} M. Boileau and B. Zimmermann, \emph{ Symmetries of
nonelliptic Montesinos links}. Math. Ann. \textbf{277} (1987), no. 3,
563--584.

\bibitem[Boy]{Boy} S. Boyer, \emph{On the local structure of ${\rm
SL}(2,{\mathbb C})$-character varieties at
reducible characters}, Topology Appl. \textbf{121} (2002), 383--413.

\bibitem[BMZ]{BMZ} S. Boyer, T. Mattman and X. Zhang,
\emph{The fundamental polygons of twisted knots and the (-2,3,7)-pretzel knot}, Knots'96, 
World Scientific Publishing Co. PteLtd, 1997, 159--172.


\bibitem[BLZ]{BLZ} S. Boyer, E. Luft and X. Zhang,
\emph{On the algebraic components of the ${\rm SL}(2,\mathbb C)$
character varieties of knot exteriors}, Topology \textbf{41} (2002), 
667--694.

\bibitem[BZ1]{BZ1} S. Boyer, X. Zhang, \emph{On Culler-Shalen seminorms and
Dehn filling}, Ann. Math.
\text\bf{148} (1998), 737--801.

\bibitem[BZ2]{BZ2} S. Boyer, X. Zhang, \emph{A proof of the finite filling
conjecture}, J. Differential Geom.
\textbf{59} (2001), 87--176.

\bibitem[BZ3]{BZ3} S. Boyer, X. Zhang, \emph{On simple points of the
character varieties of 3-manifolds},
Knots in Hellas '98 (Delphi), Ser. Knots Everything \textbf{24}, World Sci.
Publishing,
River Edge, NJ, 2000, 27--35.

\bibitem[Br]{Br} W. Browder, \emph{Surgery on simple-connected manifolds},
 Berlin-Heidelberg-New York, Springer, 1972.

\bibitem[Brn]{Brn} K. Brown \emph{cohomology of groups},
 Berlin-Heidelberg-New York, Springer GTM, 1982.
 
 \bibitem[Bu]{Bu} G. Burde, \emph{$SU(2)$-representation spaces for two-bridge knot groups}, Math. Ann. {\bf 288} (1990), 103--119.

\bibitem[CJ]{CJ} A.~Casson and D.~Jungreis.
\emph{Convergence groups and Seifert fibered $3$-manifolds}, Invent. Math.
\textbf{118} (1994), 441--456.

\bibitem[CCGLS]{CCGLS}
D. Cooper, M. Culler, H. Gillet, D. D. Long, and P. Shalen,
{\it Plane curves associated to character varieties of $3$-manifolds}, Invent.
Math. {\bf 118} (1994), 47--84.

\bibitem[CL]{CL} D. Cooper and D. Long, \emph{Roots of unity and the character variety of a knot complement},  J. Aus. Math. Soc. {\bf 55} (1993), 90--94.	

\bibitem[CGLS]{CGLS} M. Culler, C. Gordon, J. Luecke, and P. Shalen,
\emph{Dehn surgery on knots}, Ann. of Math., \textbf{125} (1987), 237--300.


\bibitem[CJR]{CJR} M. Culler, W. Jaco and H. Rubistrein,
\emph{Incompressible surfaces in once-punctured torus
bundles}, Proc. London Math. Soc. (3) \textbf{45}, (1982), 385--419.

\bibitem[CS]{CS} M. Culler and P. B. Shalen.
\emph{Varieties of group representations and splittings of 3-manifolds},
Ann. of Math. \textbf{117} (1983), 109--146.


\bibitem[Dun]{Dun} N. Dunfield, \emph{Cyclic surgery, degrees of maps of
characters curves and volume rigidity
for hyperbolic manifolds}, Invent. Math. \textbf{136} (1999), 623--657.


\bibitem[FH]{FH} W. Floyd and A. Hatcher, \emph{Incompressible surfaces in
punctured-torus bundles}, Topology and its application \textbf{13} (1982),
263--282.

\bibitem[Fra]{Fra} S. Francaviglia, \emph{Hyperbolic volume of representations of fundamental groups of cusped $3$-manifolds}, arXiv:math.GT/0305275v1

\bibitem[Ga]{Ga} D.~Gabai,
\emph{Convergence groups are Fuchsian groups}, Ann. of Math. \textbf{136} (1992), 447--510.

\bibitem[GLM]{GLM} C.McA. Gordon, R. Litherland and K. Murasugi, {\em
Signatures of covering links}, Can. J. Math. \textbf{33}, 1981, 331-394.

\bibitem[GR1]{GR1} F. Gonz\`alez-Acuna and A. Ramirez \emph{Two-bridge knots with property Q}, Quart. J. Math. {\bf 52} (2001), 447--454. 

\bibitem[GR2]{GR2} F. Gonz\`alez-Acuna and A. Ramirez \emph{Epimorphisms of knot groups onto free products}, Topology {\bf 42} (2003), 1205--1227. 

\bibitem[GW]{GW} F. Gonz\`alez-Acuna and W. Whitten \emph{Imbeddings of three-manifold groups}, Mem. Amer. Math. Soc. {\bf 99}, 1992.  

\bibitem[Gre]{Gre} L. Greenberg, \emph{Discrete subgroups of the Lorentz group}, Math. Scand. \textbf{10} (1962), 85--107.

\bibitem[Gro]{Gro} M. Gromov, 
 \emph{Metric structure for riemannian and non-riemannian spaces}, Progress in Math. \textbf{110},
Birkh\"auser 1999.



\bibitem[HW]{HW} G. Hardy and G. Wright, {\it An Introduction to the Theory of Numbers}, 4th ed.,  Oxford University Press, 1975. 

\bibitem[Ha]{Ha} R. Hartley, \emph{Knots with free period}, Canad. J. Math. \textbf{33} (1981), 91--102

\bibitem[HM]{HM} R. Hartley and K. Murasugi. \emph{Homology invariants}, Canad. J. Math. \textbf{30} (1978), 655--670.



\bibitem[Hat]{Hat} A. Hatcher, \emph{On the boundary curves of
incompressible surfaces}, Pac. J. Math.
\textbf{99} (1982), 373--377.


\bibitem[HT]{HT} A. Hatcher and W. Thurston. \emph{Incompressible surfaces in
 2-bridge knot complements}, Invent. Math. \textbf{79} (1985), 225--246.

\bibitem[HWZ]{HWZ} C. Hayat-Legrand, S.C. Wang and H. Zieschang,
\emph{Minimal Seifert manifolds},
 Math. Ann. \textbf{308} (1997), 673--700.
 
\bibitem[HS1]{HS1} J. Hoste and P.D. Shanahan,
\emph{A formula for the A-polynomial of twist knots},
 J. Knot Theory and Ramifications \textbf{13} (2004), 193--209. 

\bibitem[Hu]{Hu} H. Huang, \emph{Branched coverings and non-zero degree
maps between Seifert manifolds},
Proc. Amer. Math. Soc. \textbf{130} (2002), 2443--2449.

\bibitem[In]{In}  G. Indurskis, \emph{Remplissages d'une composante du bord de l'ext\'erieur de l'entrelacs de Whitehead}, UQAM doctoral dissertation, 2005. 

\bibitem[Ja]{Ja} W.H.Jaco
{\em Lectures on three-manifold topology}, AMS, 1977.



\bibitem[JM]{JM} T. J\o rgensen, J. Marsden, \emph{Algebraic and geometric convergence of Kleinian groups}, Math. Scand. {\bf 66} (1990), 47--72.

\bibitem[Ki]{Ki} R. Kirby, \emph{Problems in low dimensional topology},
Geometric Topology,
Edited by H. Kazez, AMS/IP 1997.

\bibitem[Kl]{Kl} E. Klassen, \emph{Representations of knots groups in $SU(2)$},
Trans. Amer. Math. Soc. {\bf 326} (1991), 795--828.
 

\bibitem[Kn]{Kn} A.W. Knapp, \emph{Doubly generated Fuchsian groups}, Mich. Math. J. {\bf 15}
(1968), 289--304. 

\bibitem[LM]{LM} A.\ Lubotzky and A.\ R.\ Magid,
\emph{Varieties of representations of finitely generated groups},
Mem. Amer. Math. Soc. \textbf{58}, 1985.

\bibitem[MT]{MT} K.~Matsuzaki and M.~Taniguchi, \emph{Hyperbolic manifolds
and Kleinian groups},
Oxford Math. Monographs, Oxford, 1998.

\bibitem[Mat]{Mat}  T. Mattman, \emph{The Culler-Shalen seminorms of the $(-2,3,n)$ pretzel knot}, J. Knot Theory and Ramifications \textbf{11} (2002), 1251--1289. 

\bibitem[Mey]{Mey}  R. Meyerhoff, \emph{The cusped hyperbolic 3-orbifold of minimum volume}, Bull. Amer. Math. Soc. \textbf{13} (1985), 154--156. 

\bibitem[MS1]{MS1}  J. Morgan and P. Shalen, \emph{Valuations, trees, and degenerations of hyperbolic structures, I}, Ann. Math. \textbf{120} (1984), 401--476. 


\bibitem[MS3]{MS3}  J. Morgan and P. Shalen, \emph{Valuations, trees, and degenerations of hyperbolic structures, III}, Ann. Math. \textbf{127} (1988), 467--519. 

\bibitem[Ne]{Ne} P. E. Newstead, \emph{Introduction to moduli problems and orbit spaces}, 
Tata Institute Lecture Notes, Narosa Publishing House, New Delhi, 1978. 

\bibitem[Oht]{Oht}  T. Ohtsuki, \emph{Ideal points and incompressible surfaces in two-bridge knot complements},   J. Math. Soc. Japan \textbf{46} (1994), 51--87.

\bibitem[ORS]{ORS} T. Ohtsuki, R. Riley and M Sakuma,
\emph{Epimorphisms between 2-bridge link groups}, preprint, June
2006.

\bibitem[ReW]{ReW} A.W. Reid and S. Wang,  \emph{Non-Haken $3$-manifolds are not large with respect to mappings of non-zero degree}, Comm. Anal. Geom, \textbf{7} (1999), 105--132.

\bibitem[RWZ]{RWZ} A.W. Reid, S. Wang and Q. Zhou \emph{Generalized hopfian property, a minimal Haken manifold, and epimorphisms between 3-manifold groups}, Acta Math. Sinica, English Ser. 
\textbf{18} (2002), 157--172.

\bibitem[Re]{Re} A. Reznikov, \emph{Analytic topology}, European Congress of Mathematics, Vol. I (Barcelona, 2000), 519--532, Progr. Math. {\bf 201}, BirkhŠuser, Basel, 2001

\bibitem[Ri]{Ri} Y.R. Riley, \emph{Knots with parabolic property P}, Quart. J. Math.
\textbf{35} (1974), 273--283.

\bibitem[Ron1]{Ron1} Y. Rong, \emph{Maps between Seifert spaces of infinite
$\pi_1$}, Pacific J. Math.
\textbf{160} (1993), 143--154.

\bibitem[Ron2]{Ron2} Y. Rong, \emph{Degree one maps of Seifert manifolds
and a note on Seifert volume},
Top. Appl.. \textbf{64} (1995), 191--200.


\bibitem[RoW]{RoW} Y. Rong and S.C.Wang, \emph{The preimage of submanifolds},
Math. Proc. Camb. Phil. Soc. \textbf{112} (1992), 271--279.


\bibitem[Sak]{Sak} M. Sakuma, \emph{The geometries of spherical Montesinos links},
Kobe J. Math. \textbf{7} (1990), 167--190.

\bibitem[Sc]{Sc} P. Scott, \emph{The geometries of $3$-manifolds},
Bull. Lond. Math. Soc. \textbf{15} (1983), 401--487.

\bibitem[Se]{Se} J.-P. Serre, \emph{Repr\'esentations lin\'eaires des groupes finis},
Hermann, Paris, 3rd ed, 1978.

\bibitem[So1]{So1} T. Soma, \emph{Degree-one maps between hyperbolic
3-manifolds with the same volume limit},
Trans. Amer. Math. Soc. \textbf{353} (2001),  2753--2772.


\bibitem[So2]{So2} T. Soma, \emph{Epimorphism sequences between hyperbolic
3-manifold groups}, Proc. Amer. Math.
Soc. \textbf{130}, (2002), 1221--1223.

\bibitem[Tak]{Tak} M. Takahashi, \emph{Two-bridge knots have property P}, Mem. Amer. Math. Soc. \textbf{29}, 1981.

\bibitem[Tan]{Tan} D. Tanguay,  {\it Chirurgie finie et noeuds rationnels}, Ph. D. thesis,
UQAM, 1995.



\bibitem[Thu]{Thu} W. Thurston, \emph{Topology and Geometry of 3-manifolds},
Princeton Lecture Notes, 1978.


\bibitem[Wan1]{Wan1} S.C.Wang, \emph{The existence of maps of nonzero
degree between aspherical 3-manifolds},
Math. Z. \textbf{208} (1991), 147--160.


\bibitem[Wan2]{Wan2} S.C. Wang, \emph{Non-zero degree one maps between
3-manifolds}, Proceedings of ICM 2002, Beijing, 457--468.

\bibitem[Wal1]{Wal1} F.~Waldhausen.
 \emph{Gruppen mit Zentrum un 3-dimensionale Mannigfaltigkeiten},
Topology \textbf{6} (1967), 505--517.


\bibitem[Wal2]{Wal2} F.~Waldhausen.
 \emph{On irreducible $3$-manifolds which are sufficiently large},
 Ann. of Math. \textbf{87} (1968), 56--88.
 
 \bibitem[Wed]{Wed} R.~Weidmann, 
 \emph{On the rank of amalgamated products and product knot groups},
Math. Ann. \textbf{312} (1998), 761--771.




\end{thebibliography}
\end{document}